\documentclass{amsbook}
\usepackage{epsfig}
\usepackage[arrow,cmtip,line,matrix]{xypic}

\newcommand{\bd}{\partial \kern -.1em}
\newcommand{\inter}[1]{{\mathaccent'27 #1}}
\newcommand{\tinter}{{\/\hbox{\rm Int}\thinspace}}
\newcommand{\hint}{\psi}

\newcommand{\UU}{\mathbf u}
\newcommand{\VV}{\mathbf v}
\newcommand{\WW}{\mathbf w}
\newcommand{\x}{\times}

\newcommand{\sub}{\subseteq}
\newcommand{\into}{\rightarrow}
\newcommand{\cv}[2]{\left(\begin{matrix} #1 \\ #2 \end{matrix} \right)}
\newcommand{\mtbt}[4]
  {\left(\begin{matrix} #1 & #2 \\ #3 & #4 \end{matrix} \right)}
\newcommand{\dtbt}[4]
    {\left|\begin{matrix} #1 & #2 \\ #3 & #4 \end{matrix} \right|}
\newcommand{\Mob}{{M\"obius band\relax{}}}
\newcommand{\defit}[1]{{\itshape "#1"}}
\newcommand{\comps}{\hbox{\bf C}}
\newcommand{\reals}{\hbox{\bf R}}

\newcommand{\ints}{\hbox{\bf Z}}
\newcommand{\FIGURE}[2]{\begin{figure}\vskip 10pt
  \psfig{figure=./#1.eps}
  \caption{#2\label{#1}}
  \end{figure}
  \vskip 10pt
}
\newcommand{\strutdepth}{\dp\strutbox}
\newcommand{\marginalnote}[1]
   {\strut\vadjust{\kern-\strutdepth\domarginalnote{#1}}}
\newcommand{\domarginalnote}[1]{\vtop to \strutdepth{
  \baselineskip\strutdepth
   \vss\llap{ #1\ \ }\null}}  %get sevenpoint font here
\newcommand{\showlabels}{\global\let\showlabelflag=T}
\newcommand{\hidelabels}{\global\let\showlabelflag=F}
\let\showlabelflag=F
\newcommand{\mylabel}[1]{\ifx T\showlabelflag\marginalnote{#1}\fi\label{#1}}

\newcommand{\BlackBox}{\hfill\vrule height8.0902pt width5pt depth0pt}
\newenvironment{demo}[1]{\textsc{#1}}{\BlackBox\par\vspace{3pt}}
\newenvironment{remark}{\textsc{Remark}}{\par\vspace{3pt}}

\newtheorem{dummy}{realdumb}[section]
\newtheorem{thm}[dummy]{Theorem}
\newtheorem{lemma}[dummy]{Lemma}

\newtheorem{cor}{Corollary}[dummy]

\title{Seifert Fibered Spaces \\[5pt] Notes for a course given in
the Spring of 1993}

\author{Matthew G. Brin}

\address{Binghamton University \\ Binghamton, NY 13902-6000 \\ USA}

\address{\begin{center}\copyright 1993 All rights
reserved\end{center}}

\thanks{As an exercise in learning \LaTeX, these notes were reset in
January of 1996.  The appearance has changed, but the content
remains the same}

\begin{document}

\setlength{\baselineskip}{18pt plus 1pt minus 1pt}

\maketitle

\newpage

\vbox to 1 truein{}

\begin{center}{\bfseries\LARGE Preface to the arXiv
edition}\end{center}

These notes have resided on my web page for as long as I have had a
web page.  They have been referred to at least once in a published
article.  It has taken me this long (1993--2007) to figure out that
it would make sense for them to reside in the arXiv instead.  I
don't know that the arXiv has been in existence since 1993, but it
has been around long enough for me to know better.  From now on my
web page will have a link to the copy in the arXiv and I will get
out of the business of serving these notes to the public.

The notes are not complete.  They contain a description of compact
three dimensional Seifert fibered spaces and a classification up to
homeomorphism of compact three dimensional Seifert fibered spaces
with non-empty boundary.  It would have been nicer to go all the way
to a classification of all compact three dimensional Seifert fibered
spaces both with and without boundary, but time ran out and I never
added to the notes beyond what was covered in the course.

\vspace{10pt}

{}\hfill{\itshape Binghamton, 2007}

\tableofcontents
\chapter{Basics}

\section{Introduction}

Most of this chapter is taken from the first 8 of the 15 sections of
Seifert's original paper on fibered spaces: H. Seifert, ``Topologie
dreidimensionaler gefaserte R\"aume,'' Acta Math. 60 (1932),
147--238.  An English translation by W. Heil: ``Topology of
3-dimensional fibered spaces,'' appears starting on page 359 of
``Seifert and Threlfall: A textbook in topology'' (yes that is the
full title) by (of course) H. Seifert and W. Threlfall, an English
translation by M. Goldman of Seifert and Threlfall's 1934 book
``Lehrbuch der Topologie.''  The English volume in question is
edited by J. Birman and J. Eisner and was published in 1980 by
Academic Press, New York, as volume 89 in the series on Pure and
Applied Mathematics.  The examples in this chapter were mostly taken
from the book by W. Jaco that is identified more completely in the
introductory paragraph to Chapter 2.

Seifert fibered spaces are 3-manifolds that are unions of pairwise
disjoint circles.  There will be restrictions on how the circles fit
together.  These restrictions will be discussed shortly.  If \(M\)
is a Seifert fibered space, then the circles just referred to are
the \defit{fibers} of \(M\).  If \(M\) and \(N\) are Seifert fibered
spaces, then a \defit{fiber preserving homeomorphism} from \(M\) to
\(N\) is a homeomorphism from \(M\) to \(N\) that takes each fiber
of \(M\) onto a fiber of \(N\).  If such a homeomorphism exists,
then we say that \(M\) and \(N\) are \defit{fiber homeomorphic.} Our
first tasks will be to define Seifert fibered spaces, to classify
Seifert fibered spaces up to fiber preserving homeomorphism, and
then to classify Seifert fibered spaces up to arbitrary
homeomorphism.  [The course ended before this last part could be
accomplished.  Compact Seifert fibered spaces with non-empty
boundary were classified.]  One expects the last classification to
be somewhat looser than the first since more homeomorphisms are
allowed.  For the most part this is not the case, and the situations
in which classes combine because of the extra homeomorphisms are few
in number and easy to describe.  To save words, we define the
\defit{fiber type} of a Seifert fibered space to be the class of all
Seifert fiber spaces that are fiber homeomorphic to it.

\section{Definition of Seifert fibered spaces}

Many spaces are defined using local models.  An \(n\)-manifold
without boundary is defined to be a separable metric space in which
every point has an open neighborhood homeomorphic to the fixed local
model \(\reals^n\).  An \(n\)-manifold with boundary is a separable
metric space in which every point has a closed neighborhood
homeomorphic to the fixed local model \({\hbox{\bf B}}^n\).  Seifert
fibered spaces are defined using an infinite collection of local
models, and these models are not models for neighborhoods of points
but for neighborhoods of fibers.
\subsection{\mylabel{LocalModel}Structure of the local model}

\FIGURE{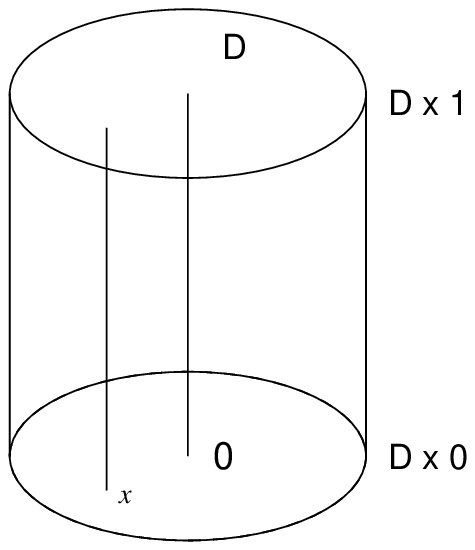}{\(D\x I\)}

Let \(D\) be the unit disk in the complex plane \(\comps\) and let
\(D\x I\) be fibered by the intervals \(\{x\}\x I, x\in D\) as shown
in Figure \ref{DxI}.

Let \(\nu\) and \(\mu\) be integers with \(\mu\ne0\).  Define
\(\rho:D\into D\) by \[\rho(x)=xe^{2\pi i(\nu/\mu)}.\] That is,
\(\rho\) rotates \(D\) by \(\nu/\mu\) of a full circle.  The
homeomorphism \(\rho\) is completely determined by the rational
number \(\nu/\mu\mod1\).  Thus we assume that \(\nu\) and \(\mu\)
are relatively prime, that \(\mu>0\) and that \(0\le\nu<\mu\).  Let
\(T\) be the quotient of \(D\x I\) obtained by identifying \((x,0)\)
with \((\rho(x),1)\) for each \(x\in D\).  Since \(\rho\) is
isotopic to the identity map, \(T\) is homeomorphic to \(D\x S^1\).

For a given \(x\in D\), there are a finite number of images of \(x\)
in \(D\) under powers of \(\rho\).  Since \(\nu\) and \(\mu\) are
relatively prime, this number of images is exactly \(\mu\) if \(x\)
is not the center of \(D\), and there is exactly one image if
\(x=0\), the center of \(D\).  We can now describe \(T\) as a union
of circles.  One circle is the image of \(\{0\}\x I\) under the
quotient map.  The other circles are each the union of images of
\(\mu\) intervals of the form \(\{\rho^i(x)\}\x I\) where \(i\) runs
from 0 to \(\mu-1\).  We call these circles the \defit{fibers} of
\(T\).  We use \(T\) to refer to both the underlying solid torus
together with the added structure consisting of the set of fibers of
\(T\).  We call such a \(T\) a \defit{fibered solid torus}.  We call
the fiber of \(T\) that is the image of \(\{0\}\x I\) the
\defit{centerline} of \(T\).  We refer to \(T\) as the fibered solid
torus determined by the rational number \(\nu/\mu\mod1\).

If \(\nu=0\), then the fibering of \(T\) is the product fibering of
\(D\x S^1\).  In this case we call \(T\) an \defit{ordinary solid
torus}.

There are many ordinary solid tori in a fibered solid torus.  Let
\(x\) be a point of \(D\) other than the center and let \(T\) be the
fibered solid torus determined by the rational number \(\nu/\mu\) as
a quotient of \(D\x I\).  Let \(H\) be the fiber of \(T\) containing
the image of \(\{x\}\x I\).  There is a disk \(E\) in \(D\) that
contains a neighborhood \(x\) and does not contain the center of
\(D\).  If \(x\) is in \(\bd D\), then \(x\) is in \(\bd E\), and if
\(x\) is in \(\inter D\), then \(x\) is in \(\inter E\).  By
choosing \(E\) small enough, we can keep \(E\) disjoint from its
\(\mu\) images under powers of \(\rho\).  Now the union of the
images of \(\rho^i(E)\x I\), \(0\le i<\mu\) is a solid torus \(T'\)
in \(T\) fibered as an ordinary solid torus by the fibers of \(T\)
that it contains.  If \(x\) is in \(\inter E\), then we can regard
\(x\) as the center point of \(E\).  Thus \(H\) can be regarded as
the centerline of \(T'\).  So every fiber in \(T\) that is not the
centerline of \(T\) is contained in a union of fibers of \(T\) that
forms an ordinary solid torus in \(T\).  Further, if the fiber is in
\(\inter T\), then it is the centerline of this ordinary solid
torus.

A fibered solid torus \(T\) is orientable.  If \(D\x I\) is given an
orientation, then \(T\) can inherit the orientation from \(D\x I\).
We will show below that the choice of orientation is important in
determining certain invariants.  We will therefore choose consistent
orientations for the various local models.  We will do this by
picking one orientation for \(D\x I\) and letting the models inherit
that orientation.  When we discuss below invariants of fibered solid
tori, we will then describe our choice of orientation for \(D\x I\).

If \(T\) is the fibered solid torus determined by \(\nu/\mu\) as a
quotient of \(D\x I\), then there is a fiber preserving
homeomoprphism from \(T\) to the image of \(D'\x I\) where \(D'\) is
the half unit disk in \(\comps\).  The homeomporhism is induced by
radial dilation by one half.

\subsection{The definition}

Let a 3-manifold \(M\) be a pairwise disjoint union of simple closed
curves called the \defit{fibers} of \(M\).  We will use \(M\) to
refer to the underlying 3-manifold together with the added structure
consisting of the set of fibers of \(M\).  A subspace of \(M\) is
\defit{saturated} if it is a union of fibers.  Such a subspace is
given the extra structure of its set of fibers.  We say that \(M\)
is a \defit{Seifert fibered space} if each fiber \(H\) of \(M\) has
a saturated, closed neighborhood \(N\) with a fiber preserving
homeomorphism from \(N\) to a fibered solid torus \(T\).  If \(M\)
is oriented, then we require that the fiber preserving homeomorphism
from \(N\) to \(T\) be orientation preserving.

We pick out certain saturated torus neighborhoods of fibers for
special consideration.  If a fiber \(H\) has a point in \(\inter
M\), then \(H\sub\inter N\) and \(H\sub\inter M\).  Thus a fiber of
\(M\) lies either entirely in \(\bd M\) or entirely in \(\inter M\).
If \(H\) lies in \(\bd M\), then a saturated neighborhood \(N\) of
\(H\) that is fiber homeomorphic to a fibered solid torus must have
\(H\) in \(\bd N\).  By previous remarks, a possibly smaller
saturated neighborhood of \(H\) will be fiber homeomorphic to an
ordinary solid torus.  That this might not be all of \(N\) is seen
by letting \(M\) be a non-ordinary fibered solid torus and letting
\(N=M\).

If \(H\) lies in \(\inter M\), then there is a saturated
neighborhood \(N\) of \(H\) that has a homeomorphism to a fibered
solid torus \(T\).  The homeomorphism carries \(H\) into \(\inter
T\), and may or may not carry \(H\) to the centerline of \(T\).  If
\(H\) is not carried to the centerline of \(T\), then previous
remarks imply that there is a smaller saturated neighborhood of
\(H\) that is fiber homeomorphic to an ordinary solid torus \(T'\)
under a homeomorphism that carries \(H\) to the centerline of
\(T'\).  By narrowing the neighborhood, we can preserve all that we
have and insist that the neighborhood lie in \(\inter M\).

Thus every fiber in \(M\) is either in \(\bd M\) and contained in
the boundary of an ordinary solid torus, or is in \(\inter M\) and
is the centerline of a fibered (perhaps ordinary) solid torus in
\(\inter M\).  For a fiber \(H\) in \(M\), we say that \(N\) is a
\defit{fibered solid torus neighborhood} of \(H\) if \(H\) is in
\(\bd M\) and \(N\) is a saturated neighborhood of \(H\) that is
fiber homeomorphic to an ordinary solid torus, or if \(H\) is in
\(\inter M\) and \(N\) lies in \(\inter M\) and has a fiber
preserving homeomorphism to a fibered solid torus \(T\) that carries
\(H\) to the centerline of \(T\).  We will see later that \(H\)
determines its fibered solid torus neighborhoods in a strong way.

Since a fiber of \(M\) that intersects that boundary of \(M\) must
lie in the boundary of \(M\) and have an orientable neighborhood in
\(\bd M\), it follows that boundary components of \(M\) are
saturated tori, Klein bottles and open annuli.  Klein bottles can
arise because it is possible to represent a Klein bottle as a union
of annuli.

\section{Fibered solid tori}

Before we look at Seifert fibered spaces, it will be useful to study
the local model.

\subsection{Background}

We need some information about solid tori and their boundary tori.
We will take the facts in this item as given.  Let \(T\) be a solid
torus and let \(\bd T\) be its boundary torus.  Both \(\pi_1(T)\)
and \(H_1(T)\) are isomorphic to \(\ints\) and both \(\pi_1(\bd T)\)
and \(H_1(\bd T)\) are isomorphic to \(\ints\oplus\ints\).

We concentrate for a while on \(H_1(\bd T)\).  Once an ordered pair
of generators is chosen for \(H_1(\bd T)\), we can think of elements
of \(H_1(\bd T)\) as column vectors \(\cv{i}{j}\) with integer
entries where \(i\) is the coefficient of the first generator and
\(j\) is the coefficient of the second.  Automorphisms of \(H_1(\bd
T)\) are then represented by two by two integer matrices
\(\mtbt{a}{b}{c}{d}\) with \(\dtbt{a}{b}{c}{d}=\pm1\).  If a
different pair of generators of \(H_1(\bd T)\) is chosen, then the
matrix representing an automorphism will change, but its determinant
will not.  Thus the determinant of an automorphism of \(H_1(\bd T)\)
is well defined.

Given a pair of generators \((\UU,\VV)\) of \(H_1(\bd T)\), there is
a unique anti-symmetric pairing (bilinear form) \(\hint(\,\,,\,)\)
on \(H_1(\bd T)\) for which \(\hint(\UU,\VV)=+1\).  It is defined by
\(\hint\left(\cv{a}{b},\cv{c}{d}\right)=\dtbt{a}{c}{b}{d}\).  The
representing matrix of the pairing is \(\mtbt{0}{1}{-1}{0}\).  For
an element \(\WW=a\UU+b\VV\) of \(H_1(\bd T)\), we have
\(\hint(\WW,\UU)=-b\), \(\hint(\WW,\VV)=a\), and
\(\WW=\hint(\WW,\VV)\UU-\hint(\WW,\UU)\VV=
\hint(\WW,\VV)\UU+\hint(\UU,\WW)\VV\).  This is useful since it
allows coefficients to be read off from the pairing.  The pairing is
preserved by automorphisms of determinant \(+1\) and negated by
automorphisms of determinant \(-1\).  This pairing and its negative
are the only two anti-symmetric pairings on \(H_1(\bd T)\) whose
image is all of \(\ints\).  We can think of a choice of orientation
of \(\bd T\) either as a choice of one of the two classes of ordered
pairs of generators of \(H_1(\bd T)\) connected by automorphisms of
determinant \(+1\), or as a choice of one of the two anti-symmetric
pairings on \(H_1(\bd T)\) whose image is all of \(\ints\).  A
choice \((\UU,\VV)\) of generating pair and a choice of pairing
\(\hint(\,\,,\,)\) is consistent if \(\hint(\UU,\VV)=+1\).  We can
abuse notation and think of this as a choice of orientation for
\(H_1(\bd T)\) as well and label automorphisms of \(H_1(\bd T)\) as
orientation preserving or orientation reversing as given by the
determinant.

Let a pair of generators \((\UU,\VV)\) of \(H_1(\bd T)\) be given
and let \(\hint(\,\,,\,)\) be the anti-symmetric pairing on
\(H_1(\bd T)\) for which \(\hint(\UU,\VV)=+1\).  Assume that
\[\hint\left(\cv{a}{b},\cv{c}{d}\right)=\epsilon=\pm1\] for elements
\(\cv{a}{b}\) and \(\cv{c}{d}\) in \(H_1(\bd T)\).  Then the matrix
\(\mtbt{a}{c}{b}{d}\) represents an automorphism of \(H_1(\bd T)\)
which preserves orientation depending on \(\epsilon\) and which
carries \(\UU\) to \(\cv{a}{b}\) and \(\VV\) to \(\cv{c}{d}\).  Thus
\(\left(\cv{a}{b},\cv{c}{d}\right)\) is a generating pair for
\(H_1(\bd T)\) and gives the same orientation as \((\UU,\VV)\) if
and only if \(\epsilon=+1\).  Thus the pairing picks out generating
pairs for \(H_1(\bd T)\).

Another view of an orientation for \(\bd T\) is that it gives a
consistent way of assigning one of \(\pm1\) to an isolated,
transverse intersection of an ordered pair of oriented open arcs in
\(\bd T\).  This can be extended by summing to a definition of
intersection number for ordered pairs of oriented closed curves that
can be shown to depend only on the homology classes that the curves
represent and to be bilinear.  This then gives a well defined notion
of intersection numbers on ordered pairs of homology classes in
\(H_1(\bd T)\).  By considering the ``obvious'' generators of
\(H_1(\bd T)\) --- two coordinate slices of \(S^1\x S^1\) --- it can
be seen that the intersection numbers give an anti-symmetric pairing
whose image is all of \(\ints\).  Thus the pairing is one of the
``orientation pairings'' discussed above.  We can now refer to the
``orientation pairings'' as intersecton pairings.

We can use the intersection pairing to gather information about
simple closed curves on \(\bd T\).  A simple closed curve on \(\bd
T\) separates if and only if it is null homologous in \(\bd T\).  If
a simple closed curve \(J\) represents a non-trivial element of
\(H_1(\bd T)\), then it does not separate, and it is transverse to
some simple closed curve \(K\) that it intersects exactly once.
Thus \(\hint(J,K)=\pm1\) and \((J,K)\) represents a generating pair
for \(H_1(\bd T)\).  Further, the homology class \([J]\) that \(J\)
represents is indivisible in \(H_1(\bd T)\) --- it is an integer
multiple of no element other than itself or its inverse --- and if
\([J]=\cv{a}{b}\) with respect to an arbitrary pair of generators
for \(H_1(\bd T)\), then \(a\) and \(b\) are relatively prime.  A
converse is true.  If a class in \(H_1(\bd T)\) is indivisible, then
it is represented by a simple closed curve.  This can be seen by
using the ``obvious'' generating pair, noting that the coefficients
of the element must be relatively prime, and then plotting a lift of
the desired curve in the universal cover as a straight line with
slope the quotient of the coefficients.

We need more information about curves.  (Is a good reference for
this paragraph the paper by Epstein on curves on surfaces and
isotopies?)  If two simple closed curves on \(\bd T\) have
intersection pairing \(n\), then one curve can be altered by an
isotopy so that the two curves intersect transversely and intersect
in exactly \(|n|\) points.  From this it follows that every
automorphism of \(H_1(\bd T)\) is realized by some homeomorphism of
\(\bd T\) because the new generators can now be represented by a
pair of transverse simple closed curves that intersect in one point
and it is easy to construct a homeomorphism carrying the ``obvious''
generators to the new generating curves.  A second fact classifies
homeomorphisms.  Two simple closed curves on \(\bd T\) are
homologous if and only if they are homotopic in \(\bd T\), in which
case they are also isotopic in \(\bd T\).  (The last part
generalizes to all surfaces.  Two embedded curves are homotopic if
and only if they are isotopic where the embeddings, homotopies and
isotopies are all required to preserve boundaries in that inverse
images of boundaries are always boundaries.)  From this information
and the Alexander Lemma it can be shown that a homeomorphism that
induces the identity on \(H_1(\bd T)\) is isotopic to the identity
on \(\bd T\).  Thus self homeomorphisms of \(\bd T\) are classified
up to isotopy by the induced automorphisms on \(H_1(\bd T)\) and,
once an ordered pair of generators for \(H_1(\bd T)\) is chosen, by
two by two integer matrices of determinant \(\pm1\).  Chosing
ordered pairs of generators for the first homologies of a pair of
tori allows the classification of homeomorphisms, up to isotopy,
between the tori by two by two, invertible, integer matrices.

We call a simple closed curve on \(\bd T\) that is null homotopic in
\(T\) but not in \(\bd T\) a \defit{meridian} of \(T\).  There are
only two classes in \(H_1(\bd T)\) that contain meridians of \(T\).
These classes are the negatives of each other so that there is only
one meridian of \(T\) up to isotopy and reversal of orientation.  A
self homeomorphism of \(\bd T\) extends to all of \(T\) if and only
if the homology class of the meridian is preserved up to sign.  We
call any simple closed curve \(\ell\) on \(\bd T\) a \defit{longitude}
for \(T\) if there is a meridian for \(T\) that intersects \(\ell\)
exactly once and which is transverse to \(\ell\).  A longitude for
\(T\) represents a generator of \(\pi_1(T)=H_1(T)=\ints\) and a
meridian-longitude pair represents a generating pair for \(H_1(\bd
T)\).  Any meridian-longitude pair for \(T\) can be carried to any
other meridian-longitude pair for \(T\) by a self homeomorphism of
\(T\).

We will be interested in understanding fiber preserving
homeomorphisms of fibered solid tori.  We will use this information
to obtain numerical invariants that classify fibered solid tori up
to fiber preserving homeomorphism, and to obtain geometric
invariants that allow us to unambiguously specify sewings of fibered
solid tori up to fiber preserving deformations.  We will see that
all of the information that we need is contained in the boundaries
of the fibered solid tori, and that in the boundaries, it is the
homology classes of curves that are important.  Understanding how
classes of curves in a torus behave under certain homeomorphisms of
the torus needs an understanding of how subgroups of
\(\ints\oplus\ints\) behave under certain automorphisms of
\(\ints\oplus\ints\).

\subsection{Automorphisms of \protect\(\ints\oplus\ints\protect\)}

Let \(T\) be a fibered solid torus.  Let \(H\) be a fiber of \(T\)
on \(\bd T\), and let \(m\) be a meridian for \(T\).  Fiber
preserving self homeomorphisms of \(T\) will take \(H\) to \(H\) or
its reverse and \(m\) to \(m\) or its reverse.  The induced
automorphism on \(H_1(\bd T)\) will preserve the cyclic subgroups
generated by the classes of \(H\) and of \(m\).  These cyclic
subgroups are maximal cyclic.  They are also different since
\(\hint(H,m)\) is never 0.  This motivates what we look at in this
section.

Let \(G=\ints\oplus\ints\).  We identify \(G\) with column vectors
\(\cv{i}{j}\).  We use the notation \((i,j)\) for the greatest
common divisor of \(i\) and \(j\).  A cyclic subgroup of \(G\) is
generated by some \(\cv{i}{j}\) in \(G\) other than the zero element
\(\cv{0}{0}\).  This subgroup will be maximal cyclic in \(G\) if and
only if \((i,j)=1\).  Every non-zero element of \(G\) is contained
in a unique maximal cyclic subgroup of \(G\).  The element
\(\cv{i}{j}\) is contained in the maximal cyclic subgroup generated
by
\[\cv{\displaystyle\frac{i}{(i,j)}}{\displaystyle\frac{j}{(i,j)}}.\]
We use the convention that \((i,0)=i\) if \(i\ne0\).  It follows
that \(\cv{i}{j}\) and \(\cv{i'}{j'}\) are contained in the same
maximal cyclic subgroup if and only if \(i/j=i'/j'\) where all
fractions with denominator 0 are taken to be the same.  We can thus
use fractional notation to identify maximal cyclic subgroups of
\(G\), and we use \([i/j]\) to refer to the unique maximal cyclic
subgroup of \(G\) containing \(\cv{i}{j}\).

Let \(\phi\) be an automorphism of \(G\).  We identify \(\phi\) with
the \(2\x2\) integer matrix that represents \(\phi\).  We can thus
talk about \(\det(\phi)\) and we know that \(\det(\phi)=\pm1\).  We
say that \(\phi\) is \defit{orientation preserving} if
\(\det(\phi)=1\) and \defit{orientation reversing} if
\(\det(\phi)=-1\).

In the following, we consider automorphisms that fix \([1/0]\)
because homeomorphisms of the boundary of a solid torus must
preserve the meridian if they are to extend to the entire solid
torus.

\begin{lemma}\mylabel{OneGroup} Let \([\nu/\mu]\) and
\([\nu'/\mu']\) be maximal cyclic subgroups of \(G\) with
\(\mu\ne0\) and \(\mu'\ne0\).  There is an orientation preserving
automorphism \(\phi\) of \(G\) fixing \([1/0]\) and taking
\([\nu/\mu]\) to \([\nu'/\mu']\) if and only if
\(\nu/\mu\equiv\nu'/\mu'\mod1\).  There is an orientation reversing
automorphism \(\phi\) of \(G\) fixing \([1/0]\) and taking
\([\nu/\mu]\) to \([\nu'/\mu']\) if and only if
\(\nu/\mu\equiv-\nu'/\mu'\mod1\).  \end{lemma}

\begin{demo}{Proof} Assume a \(\phi\) exists carrying \([\nu/\mu]\)
to \([\nu'/\mu']\).  Because \(\phi\) preserves \([1/0]\), it must
have the form: \[\epsilon\mtbt{\omega}{b}{0}{1}\] where
\(\epsilon=\pm1\) and \(\omega=\pm1\).  Here \(\omega=\det(\phi)\)
and determines whether or not \(\phi\) preserves orientation.  We
assume that \(\nu/\mu\) and \(\nu'/\mu'\) have been reduced to
lowest terms so that \(\cv{\nu}{\mu}\) and \(\cv{\nu'}{\mu'}\) are
generators of the respective maximal cyclic subgroups.  We have
\[\pm\cv {\nu'}{\mu'}=
\epsilon\mtbt{\omega}{b}{0}{1}\cv{\nu}{\mu}=\epsilon\cv{\omega\nu+b\mu}
{\mu}.\] This immediately gives that
\(\nu'/\mu'\equiv\omega\nu/\mu\mod1\).  The converse is a
straightforward reversal of the argument.  \end{demo}

From the above lemma, the orbit of \([\nu/\mu]\) under orientation
preserving automorphisms of \(G\) that preserve \([1/0]\) can be
represented by the fraction \(\nu/\mu\) with \((\nu,\mu)=1\),
\(\mu>0\) and \(0\le\nu<\mu\).  This is the reason for the
restriction in the next lemma.

\begin{lemma}\mylabel{TwoGroups} Let \([\nu/\mu]\) be a maximal
cyclic subgroup of \(G\) with \((\nu,\mu)=1\), \(\mu>0\) and
\(0\le\nu<\mu\).  The only automorphisms of \(G\) that preserve each
of \([\nu/\mu]\) and \([1/0]\) are among
\[\phi_1=\pm\mtbt{1}{0}{0}{1}, \phi_2=\pm\mtbt{-1}{1}{0}{1},
\hbox{\,and}\,\,\, \phi_3=\pm\mtbt{-1}{0}{0}{1}.\] Only \(\phi_1\)
works for all \([\nu/\mu]\), \(\phi_2\) works only for \([1/2]\),
and \(\phi_3\) works only for \([0/1]\).  \end{lemma}

\begin{demo}{Proof} One can check that the given automorphisms work
as stated.  We show there are no more.  From the proof of the lemma
above, we have \[\pm\cv{\nu}{\mu}=\epsilon\cv{\omega\nu+b\mu}{\mu}\]
so \(\nu=\omega\nu+b\mu\).  If \(\omega=1\), then \(b=0\).  If
\(\omega=-1\), then \(2\nu=b\mu\) and either \(b=\nu=0\) or
\(\nu\mid\mu\) since \(0\le\nu<\mu\).  So if \(b\ne0\), then
\(b=1\), \(\nu=1\) and \(\mu=2\).  \end{demo}

\begin{cor} Let \([\nu/\mu]\) be a maximal cyclic subgroup of \(G\)
with \(\mu>0\) and \((\nu,\mu)=1\).  The only automorphisms of \(G\)
that preserve each of \([\nu/\mu]\) and \([1/0]\) are among
\[\phi_1=\pm\mtbt{1}{0}{0}{1},
\phi_2=\pm\mtbt{-1}{\nu}{0}{1},\hbox{\,and}\,\,\,
\phi_3=\pm\mtbt{-1}{0}{0}{1}.\] Only \(\phi_1\) works for all
\([\nu/\mu]\), \(\phi_2\) works only if \(\mu=2\), and \(\phi_3\)
works only if \(\nu=0\).  \end{cor}

\begin{demo}{Proof} This is obtained by conjugating the
automorphisms of the previous lemma by an automorphism that takes
\([\nu/\mu]\) to \([\nu'/\mu]\) with \(0\le\nu'<\mu\).  \end{demo}

[I am indebted to Yuncherl Choi for pointing out an omission of the
automorphism \(\phi_3\) in the last two lemmas.]

The lemmas above apply to isomorphisms between two copies of
\(\ints\oplus\ints\) once a fixed ordered pair of generators is
chosen for each copy.

\subsection{Homeomorphisms of fibered solid tori}

We start by showing that the automorphisms of \(\ints\oplus\ints\)
identified above can be realized by fiber preserving homeomorphisms
of a fibered solid torus \(T\).  We need certain conventions first.

Let \(m\) be a meridian of \(T\) and let \(\ell\) be a longitude for
\(T\).  We identify \(H_1(\bd T)\) with \(\ints\oplus\ints\) and use
\(m\) to represent the generator of the first factor and use \(l\)
to represent the generator of the second factor.  Elements of
\(H_1(\bd T)\) are column vectors with the first (top) component
giving the coefficient of \(m\).  A self homeomorphism of \(\bd T\),
induces an automorphism of \(H_1(\bd T)\) represented by a \(2\x2\)
integer matrix with determinant \(\pm1\) that acts on elements of
\(H_1(\bd T)\) on the left.  The choice of \((m,l)\) as the
generating pair for \(H_1(\bd T)\) determines the orientation of
\(\bd T\) and the intersection pairing on \(H_1(\bd T)\).  Recall
that this pairing is anti-symmetric and has \(\hint(m,l)=1\).

Longitudes of a solid torus are not unique.  There are also two
choices for the orientation of the meridian.  We now make choices
for our convenience by referring to \(D\x I\).

Let \(T\) be a fibered solid torus obtained from \(D\x I\) using the
homeomorphism \(\rho\) determined by \(\nu/\mu\) with \(\nu\)
relatively prime to \(\mu\).  Let \(l\) be the image in \(T\) of the
curve \[L(t)=(e^{2\pi i(\nu t/\mu)},t)\] in \(D\x I\).  The ends of
\(L\) are \((1,0)\) and \((\rho(1),1)\) so that \(l\) is a simple
closed curve in \(T\) and in fact a longitude of \(T\) transverse to
the meridian \(m\) which is the image of \(\bd D\x\{0\}\).  We now
pick orientations.  We orient \(\bd D\) counterclockwise, and we
orient \(L\) from \((\rho(1),1)\) to \((1,0)\).  By declaring that
\(\bd D\) comes ``before'' \(L\), we determine an orientation of
\(\bd(D\x I)\) and thus on \(D\x I\).  Orientations are now
inherited by \(m\), \(l\), \(\bd T\) and \(T\).  We make these
choices because they cooperate well with the fibers of \(T\).  Let
\(H\) be a fiber on \(\bd T\).  If \(H\) is oriented by orienting
the fibers \(\{x\}\x I\) in \(D\x I\) from the \(1\) level towards
the \(0\) level, then \(\hint(H,m)=-\mu\) and \(\hint(H,l)=\nu\).
Thus \(H\) represents \(\nu m+\mu l\) in \(H_1(\bd T)\).  This
recovers \(\nu/\mu\) as the ratio of the coefficients of \(m\) and
\(l\).  The choice of orientation for \(H\) is actually not
important since reversing the orientation of \(H\) has \(H=-\nu
m-\mu l\) and \(\nu/\mu\) is still recovered in the same way.

Note that \(\hint(H,m)\) is never 0.  Thus \(H\) and \(m\) are never
homologous in \(\bd T\) and \(H\) is never null homologous in \(T\).

We now describe three self homeomorphisms of \(D\x I\).  The first
\(h_1\) is the identity.  The second is defined by
\(h_2(x,t)=(\overline{x},1-t)\) where \(\overline{x}\) is the
complex conjugate of \(x\) in \(\comps\).  The orientation of \(D\x
I\) is preserved under \(h_2\).  The third is defined by
\(h_3(x,t)=(\overline{x},t)\).  The orientation of \(D\x I\) is
reversed under \(h_3\).

Let \(\rho\) be the rotation of \(D\) determined by \(\nu/\mu\) and
let \(\rho'\) be the rotation determined by \(-\nu/\mu\).  The two
rotations are inverses of each other.  Let \(T\) and \(T'\) be the
fibered solid tori constructed from \(\rho\) and \(\rho'\)
respectively, oriented consistently with \(D\x I\).  We will show
that \(h_1\) and \(h_2\) induce fiber preserving self homeomorphisms
of \(T\) and \(h_3\) induces a fiber preserving homeomorphism from
\(T\) to \(T'\).  That \(h_1\) induces the identity on \(T\) is
trivial.  A typical two element set \(\{(x,0),\,(\rho(x),1)\}\) in
\(D\x I\) that maps to one point in \(T\) is taken by \(h_2\) to
\(\{(\overline{x},1),\,(\overline{\rho(x)},0)\}\).  Since
\[\overline{\rho(x)}=\overline{x}e^{-2\pi
i(\nu/\mu)}=\rho'(\overline{x}),\] the pair
\(\{(x,0),\,(\rho(x),1)\}\) is taken to
\(\{(\rho'(\overline{x}),0),(\overline{x},1)\}
=\{(\rho'(\overline{x}),0),(\rho\rho'(\overline{x}),1)\}\) and
identified pairs are carried to identified pairs.  Since the fibers
\(\{x\}\x I\) of \(D\x I\) are preserved under \(h_2\), we have an
orientation preserving, fiber preserving self homeomorphism of
\(T\).  Under \(h_3\), the two element set
\(\{(x,0),\,(\rho(x),1)\}\) maps to
\(\{(\overline{x},0),\,(\overline{\rho(x)},1)\}=
\{(\overline{x},0),\,(\rho'(\overline{x}),1)\}\) and pairs
identified in \(T\) are carried to pairs that are identified in
\(T'\).  This gives an orientation reversing, fiber preserving
homeomorphism from \(T\) to \(T'\).  We get a self homeomorphism of
\(T\) from \(h_3\) if \(T=T'\).  This happens when \(\rho=\rho'\)
which is true when \(\nu/\mu=1/2\).

The action on \(H_1(\bd T)\) induced by \(h_1\) is the identity.
The action of \(h_2\) clearly reverses \(m\).  It also reverses
\(l\) since \[h_2(L(t)) = (e^{-2\pi i(\nu t/\mu)},1-t) = (e^{2\pi
i(\nu(1-t)/\mu)}, 1-t).\] This gives two of the automorphisms of
Lemma \ref{TwoGroups}.  Two more are realized when
\(\nu/\mu=1/2\). We have \[h_3(L(t))=(e^{-2\pi i(t/2)},t).\]
Recalling that \(L(t)\) is to be oriented from \(1\) to \(0\), we
have that \(h_3(L(t))-L(t)\) is homologous to one full circuit
around \(\bd D\) in a counterclockwise direction.  Since \(h_3\)
reverses the direction of \(m\), we have that \(h_3\) realizes one
of the remaining automorphisms of Lemma \ref{TwoGroups} and a fourth
is realized by \(h_2h_3\).  It is trivial that the remaining two are
realized by \(h_3\) and \(h_2h_3\) on the ordinary fibered solid
torus.

We are now in a position to prove:

\begin{thm} Let \(T\) and \(T'\) be fibered solid tori determined by
\(\nu/\mu\) and \(\nu'/\mu'\) respectively.  Then there is a fiber
preserving homeomorphism from \(T\) to \(T'\) if and only if
\(\nu'/\mu'\equiv\pm\nu/\mu\mod1\).  The homeomorphism can be taken
to be orientation preserving if and only if
\(\nu'/\mu'\equiv\nu/\mu\mod1\), and can be taken to be orientation
reversing if and only if \(\nu'/\mu'\equiv-\nu/\mu\mod1\). \end{thm}

\begin{demo}{Proof} The if direction follows from the fact that the
defining rotations for \(T\) and \(T'\) are determined by the
fractions \(\nu/\mu\mod1\) and \(\nu'/\mu'\mod1\) and and from the
homeomorphisms \(h_2\) and \(h_3\) that we have described above.  To
prove the only if direction, let \(h\) be the promised fiber
preserving homeomorphism.  There is a (non-fiber preserving)
homeomorphism \(g\) from \(T\) to \(T'\) taking meridian to meridian
and longitude to longitude (as chosen above).  This preserves
orientation by definition.  Thus \(g^{-1}h\) preserves or reverses
orientation as \(h\) does.  A fiber \(H\) of \(T\) in \(\bd T\)
generates the maximal cyclic subgroup \([\nu/\mu]\) in \(H_1(\bd
T)\).  Since \(g^{-1}h(H)\) generates the maximal cyclic subgroup
\([\nu'/\mu']\) and \(g^{-1}h\) preserves the maximal cyclic
subgroup generated by the meridian, the only if direction follows
from Lemma \ref{OneGroup}.  \end{demo}

Some consequences are as follows.  Canonical numerical invariants of
an oriented fibered solid torus can be taken to be pairs of integers
\(\nu\) and \(\mu\) with \((\nu,\mu)=1\), \(\mu>0\), and
\(0\le\nu<\mu\).  (Equivalently, a single numerical invariant can be
taken to be a fraction \(\nu/\mu\) with \(0\le\nu/\mu<1\).  If the
fraction is in reduced terms with positive denominator, then the
numerator and denominator give the pair of integers just mentioned.)
Canonical numerical invariants of unoriented fibered solid tori can
be taken to be pairs of integers \(\nu\) and \(\mu\) with
\((\nu,\mu)=1\), \(\mu>0\), and \(0\le\nu\le\mu/2\).  In the
unoriented case, \(\nu=\mu/2\) implies that \(\nu\) is a divisor of
\(\nu\) and \(\mu\) so \(\nu=1\) and \(\mu=2\) and we have
\(0\le\nu<\mu/2\) whenever \(\mu>2\).  (A single numerical invariant
can be taken to be a fraction \(\nu/\mu\) in the inverval
\([0,1/2]\).)  The numerical invariants for the (unique) ordinary
fibered solid torus are \(\nu=0\) and \(\mu=1\).  The only fibered
solid tori that admit an orientation reversing, fiber preserving
self homeomorphism are the ones determined by \(0/1\) (the ordinary
fibered solid torus) and \(1/2\).

The unique behavior demonstrated by the non-ordinary fibered solid
torus determined by \(1/2\) and the special case in Lemma
\ref{TwoGroups} associated with \([1/2]\) will persist throughout
these notes and be the cause of various special cases.

\subsection{Homeomorphisms of fibered tori}

We will often have occasion to attach a fibered solid torus to a
Seifert fibered space by a homeomorphism defined on the boundary of
the fibered solid torus.  Here we study what is needed to specify
the attaching map, and to what extent the resulting quotient space
is determined.

We need the notion of a fiber preserving isotopy.  Let \(X\) and
\(Y\) be spaces fibered by circles.  A map \(\Phi:X\x I\into Y\) is
a \defit{fiber preserving isotopy} if each \(\Phi_t\) defined by
\(\Phi_t(x)=\Phi(x,t)\) is a fiber preserving homeomorphism from
\(X\) to \(Y\).  For a given fiber \(H\) of \(X\), this requires
that for each \(t\) the image of \(H\) is a fiber of \(Y\), but it
does not require that the image of \(H\) be the same fiber of \(Y\)
for all \(t\).  Two maps that are connected by a fiber preserving
isotopy will be said to be \defit{fiber isotopic}.

The first result applies to the uniqueness of the quotient space.
It says that the fiber structure of a fibered solid torus is
determined by two curves, a fiber and a meridian.

\begin{thm}\mylabel{FiberMeridian} Let \(T\) and \(T'\) be fibered
solid tori.  Let \(h:\bd T\into \bd T'\) be a fiber preserving
homeomorphism.  Then \(h\) extends to a fiber preserving
homeomorphism from \(T\) to \(T'\) if and only if \(h\) takes a
meridian of \(T\) to a curve homologous to a meridian of \(T'\) or
its reverse.  \end{thm}

\begin{demo}{Proof} If \(h\) extends to any homeomorphism from \(T\)
to \(T'\), then \(h\) takes a meridian of \(T\) to a meridian of
\(T'\) up to reversal.  We consider the other direction.  If \(T''\)
is another fibered solid torus, and if \(g\) is a fiber preserving
homeomorphism from \(\bd T'\) to \(\bd T''\) that extends to a fiber
preserving homeomorphism from \(T'\) to \(T''\), then \(h\) extends
to a fiber preserving homeomorphism from \(T\) to \(T'\) if and only
if \(gh\) extends to a fiber preserving homeomorphism from \(T\) to
\(T''\).  We will use this to replace \(h\) by another homeomorphism
that is more convenient.  Let meridians and longitudes be chosen on
\(T\) and \(T'\) so that when these are used as generators for
\(H_1(\bd T)\) and \(H_1(\bd T')\), then the fibers generate maximal
cyclic subgroups \([\nu/\mu]\) and \([\nu'/\mu']\) respectively
where each fraction is in the interval \([0,1)\), is in reduced
terms, and has positive denominator.  The isomorphism induced by
\(h\) on \(H_1\) must be given by one of the matrices of Lemma
\ref{TwoGroups}.  As shown in the previous section, there is a fiber
preserving homeomorphism from \(T'\) to some \(T''\) that also
realizes this matrix.  (We can chose \(T''=T'\) if the determinant
of the matrix is \(+1\).)  Let \(g\) be the restriction of this
homeomorphism to \(\bd T'\).  The action of \(gh\) on \(H_1\) is
given by the identity matrix since each of the matrices of Lemma
\ref{TwoGroups} has order two.

We can now assume that the action of \(h\) on \(H_1\) is given by
the identity matrix and concentrate once again on \(T\) and \(T'\).
We know that \(T\) and \(T'\) are quotients of \(D\x I\) using a
rotation determined by \(\nu/\mu\) for both quotients.  Thus we can
treat \(T\) and \(T'\) as the same space and we are working with a
fiber preserving self homeomorphism \(h\) of \(\bd T\) that induces
the identity on \(H_1(\bd T)\).

If \(h\) is the identity map, then \(h\) extends as desired.  If
\(h\) is fiber isotopic to the identity map, then a homeomorphism
can be constructed that reflects the isotopy on a neighborhood of
the boundary and a restriction of the identity on the rest.  The
result now follows from the next lemma.  \end{demo}

We need the notion of a crossing curve.  Let \(F\) be a fibered
torus.  That is, a torus represented as \(S^1\x S^1\) fibered by the
circles \(\{x\}\x S^1\).  A \defit{crossing curve} for \(F\) is a
simple closed curve on \(F\) that intersects each fiber in exactly
one point.  A crossing curve is necessarily transverse to all the
fibers.  If \(Q\) is a crossing curve for \(F\), then \(Q\) and any
one fiber represent generators for \(H_1(F)\).  All other crossing
curves for \(F\) represent some \(\pm Q+bH\), \(b\in\ints\) and each
such homology class is represented by some crossing curve.

\begin{lemma}\mylabel{StraightLevels} Let \(F\) be a fibered torus.
Let \(h\) be a fiber preserving self homeomorphism that induces the
identity on \(H_1(F)\).  Then \(h\) is fiber isotopic to the
identity.  \end{lemma}

\begin{demo}{Proof} Let \(Q\) be a crossing curve for \(F\).  We let
\(\tilde F\) be the cover of \(F\) determined by the subgroup of
\(\pi_1(F)\) generated by \([Q]\).  Thus \(\tilde F\) is an open
annulus and the curves \(Q\) and \(h(Q)\) lift to \(\tilde F\).  The
cover \(\tilde F\) is ``fibered'' by lines. Since \(h\) is fiber
preserving, both \(Q\) and \(h(Q)\) are crossing curves for \(F\).
Their lifts are crossing curves for \(\tilde F\).  A slide along a
neighborhood of a fiber in \(F\) results in an ``equivariant'' slide
along a neighborhood of a fiber in \(\tilde F\).  An ``equivariant''
slide along a neighborhood of a fiber in \(\tilde F\) induces a
slide along a neighborhood of a fiber in \(F\).  It is ``obvious''
that a lift of \(h(Q)\) can be taken to a lift of \(Q\) by
``equivariant'' slides of neighborhoods of fibers in \(\tilde F\).
Thus a fiber isotopy that keeps each fiber invariant setwise carries
\(h(Q)\) to \(Q\).  We assume that \(h\) has been altered by this
isotopy and now carries \(Q\) to \(Q\).  A similar exercise is now
done to a curve parallel to \(Q\) and disjoint from \(Q\).  The
slides are done so as to keep the function on \(Q\) unchanged
throughout.  We can think of \(Q\) as the zero level and the new
curve as the half level.  This can be repeated for the dyadic
rationals and in the limit, we get that \(h\) is level preserving.
The homeomorphism on \(Q\) is of degree one and is isotopic to the
identity.  This extends to all of \(F\) by doing the same on each
level.  This preserves fibers and deforms \(h\) to the identity.
\end{demo}

Note that if \(T\) is a fibered solid torus determined by a fiber
and meridian where the meridian happens to be a crossing curve for
boundary, then \(T\) is an ordinary solid torus.

The next results apply to the existence of the quotient spaces.

\begin{lemma} Let \(F\) and \(F'\) be fibered tori.  Let \(h:F\into
F'\) be a homeomorphism that takes the homology class of a fiber of
\(F\) to the homology class of a fiber of \(F'\) or its negative.
Then \(h\) is isotopic to a fiber preserving map. \end{lemma}

\begin{demo}{Proof} If \(Q\) and \(Q'\) are crossing curves for
\(F\) and \(F'\) respectively, then \(h(Q)\) is homologous to \(\pm
Q'+bH'\), \(b\in\ints\).  We can thus replace \(Q\) by a new curve
by reversing direction if necessary and adding an integral multiple
of \(H\) so that after the replacement \(h(Q)\) is homologous to
\(Q'\).  Reversing the direction of \(H\) if necessary, we get that
the induced isomorphism on \(H_1\) is represented by the identity
when the bases used are \(H\) and \(Q\) for \(H_1(T)\) and \(H'\)
and \(Q'\) for \(H_1(F')\).  Thus \(h\) is isotopic to a
homeomorphism taking \(H\) to \(H'\) and \(Q\) to \(Q'\).  By
cutting along \(H\cup Q\) and \(H'\cup Q'\) we see that the result
reduces to the Alexander lemma.  \end{demo}

\begin{lemma}\mylabel{FillExist} Let \(F\) be a fibered torus and
let \(J\) be a simple closed curve on \(F\) that is not null
homologous and is not homologous to a fiber.  Then there is a
fibered solid torus \(T\) and a fiber preserving homeomorphism
\(h:\bd T\into F\) that takes a meridian to the homology class of
\(J\).  If \(F\) is oriented, then \(T\) and \(h\) can be chosen so
that \(h\) is orientation preserving.  \end{lemma}

\begin{demo}{Proof} Choose a simple closed curve \(K\) on \(F\) so
that \(\hint(J,K)=1\).  (This is with respect to a given orientation
or to an arbitrary one if none is given.)  The class of a fiber
\(H\) is given as \(\nu J+\mu K\).  Since \(H\) is a simple closed
curve, \((\nu,\mu)=1\), and since \(J\) is not homologous to \(H\),
\(\mu\ne0\).  By reversing \(H\), we can assume that \(\mu>0\).
Since \(K\) can be altered by adding an integral multiple of \(J\)
to \(K\) and still have \(\hint(J,K)=1\), we can get
\(0\le\nu<\mu\).  We know how to build a fibered solid torus \(T\)
with meridian \(m\) and longitude \(l\) so that a fiber represents
\(\nu m+\mu l\).  A homeomorphism exists from \(\bd T\) to \(F\)
taking \(m\) to \(J\) and \(l\) to \(K\).  This preserves the (given
or chosen) orientation and takes a fiber of \(\bd T\) to a curve
homologous to \(H\).  The result follows from the previous lemma.
\end{demo}

\subsection{Attaching fibered solid tori}

Combining the above results, we get:

\begin{thm}\mylabel{UniqueFill} Let \(M\) be a Seifert fibered space
and let \(F\) be a torus boundary component of \(M\).  Let \(J\) be
a simple closed curve on \(F\) that is not homologous on \(F\) to a
fiber.  Then there is a Seifert fibered space \(M'\) that is
obtained by sewing a fibered solid torus \(T\) to \(F\) by a fiber
preserving homeomorphism that takes a meridian of \(T\) to the
homology class of \(J\) or its reverse.  If \(M\) is oriented, we
can require that the orientation of \(M'\) extends the orientation
of \(M\).  If \(M''\) is formed in a similar fashion, then there is
a homeomorphism from \(M'\) to \(M''\) that extends the identity on
\(M\).  \end{thm}

\begin{demo}{Proof} That \(M'\) can be formed as required follows
from Lemma \ref{FillExist}.  To get the orientations to cooperate,
the sewing map should take the orientation of \(\bd T\) as inherited
from \(T\) to the reverse of the orientation of \(F\) as inherited
from \(M\).  If \(M''\) is also formed to fit the requirements, then
let \(h_1:\bd T_1\into F\) and \(h_2:\bd T_2\into F\) be the
respective sewing maps.  The last claim will follow if
\(h_2^{-1}h_1\) extends to a fiber preserving homeomorphism from
\(T_1\) to \(T_2\).  However, this follows from Theorem
\ref{FiberMeridian} since a meridian of \(T_1\) is taken to a
meridian of \(T_2\) or its reverse.  \end{demo}

We do not have to attach a fibered solid torus \(T\) to a Seifert
fibered space along all of the boundary of \(T\).  In the next
result, the spaces \(M_i\) will turn out to be fiber homeomorphic to
\(M\) however we are not yet ready to prove this.  We can at least
prove its independence of the attaching map.

\begin{lemma}\mylabel{AnnularSew} Let \(M\) be a Seifert fibered
space and let \(C\) be a component of \(\bd M\).  Let \(T_i\),
\(i=1,2\), be ordinary solid tori with saturated annuli \(A_i\) in
\(\bd T_i\).  If \(M_i\) are formed by attaching \(T_i\) to \(M\) by
fiber preserving homeomorphisms from \(A_i\) to to saturated annuli
\(A'_i\) in \(C\), then there is a fiber preserving homeomorphism
from \(M_1\) to \(M_2\) that extends the identity on \(M\) minus a
collar on \(C\).  \end{lemma}

\begin{demo}{Proof} The annuli \(A'_i\) are regular neighborhoods of
fibers in \(C\).  Since the fibers are parallel, there is a fiber
preserving isotopy of \(C\) carrying \(A'_1\) to \(A'_2\).  This
implies the existence of a fiber preserving self homeomorphism of
\(M\) fixed off a collar on \(C\) that carries \(A'_1\) to \(A'_2\).
The conclusion now follows in the standard way from the next lemma.
\end{demo}

\begin{lemma} Let \(T_i\), \(i=1,2\), be ordinary solid tori, and
let \(A_i\) be saturated annuli in \(\bd T_i\).  Then any fiber
preserving homeomorphism \(f\) from \(A_1\) to \(A_2\) extends to a
fiber preserving homeomorphism \(\overline{f}\) from \(T_1\) to
\(T_2\).  \end{lemma}

\begin{demo}{Proof} Each \(T_i\) is obtained from \(D\x I\) by
sewing each \((x,0)\) to \((x,1)\).  A fiber can be used as a
longitude and \(\bd D\x\{0\}\) can be used as a meridian.  The map
\(f\) may preserve the orientation of the longitude or reverse it.
The orbit space is \(D\).  The annuli intersect the meridians in a
subinterval and the map \(f\) determines a map on \(\bd D\) that
commutes with \(f\) and the projections to \(D\).  This extends to
\(D\) and to all of \(D\x I\) and thus all of \(T_1\) by crossing
with the identity on \(I\) if \(f\) preserves the orientation of the
longitude, or by crossing with the map taking \(t\) to \(1-t\) if
the longitude is reversed.  The map created \(f'\) takes \(A_1\) to
\(A_2\) and may not agree with \(f\).  However \(f\) and \(f'\)
carry each boundary curve of \(A_1\) to the same boundary curve of
\(A_2\) and with the same effect on orientation.  By the next lemma,
\(f\) and \(f'\) are fiber isotopic and the result follows.
\end{demo}

\begin{lemma}\mylabel{AnnulusClasses} Let \(A_i\), \(i=1,2\), be
fibered annuli and let \(f\) and \(f'\) be fiber preserving
homeomprhisms from \(A_1\) to \(A_2\).  If \(f\) and \(f'\) carry
each boundary curve of \(A_1\) to the same curve of \(A_2\) and with
the same effect on orientation, then \(f\) and \(f'\) are fiber
isotopic.  \end{lemma}

\begin{demo}{Proof} The proof uses the same arguments as those used
in Lemma \ref{StraightLevels} and is simpler.  \end{demo}

\section{The orbit space of a Seifert fibered space}

Let \(M\) be a Seifert fibered space.  If we give the set of fibers
of \(M\) the quotient topology, then we have formed the \defit{orbit
surface} of \(M\).  We have not yet shown that the quotient space is
a surface so the terminology is a little premature.  The facts will
catch up shortly.

\subsection{\mylabel{FSTOrbit}The orbit space of a fibered solid
torus}

Let \(T\) be a fibered solid torus defined by \(\nu/\mu\).  Let the
set of fibers \(\tau\) of \(T\) have the quotient topology.
Regarding \(T\) as a quotient of \(D\x I\), we identify \(D\) with
the disk in \(T\) that is the image of \(D\x\{0\}\).  Each fiber of
\(T\) intersects \(D\) at least once, so \(\tau\) is a quotient of
\(D\).  If \(\rho\) is the rotation of \(D\) that is used to define
\(T\), then the points of \(D\) that lie in the same fiber are the
points of an orbit of \(\rho\).  Thus \(\tau\) is the quotient of
\(D\) under the action of \(\rho\).  This quotient is a disk and the
quotient map is a branched cover, with branch set the center point
of \(D\).  The center point of \(D\) is carried to an interior point
of \(\tau\) that we think of as the center point of \(\tau\).  In
the complement of the center, the quotient map is a \(\mu\)-fold
covering projection.  Thus the projection from \(T\) to \(\tau\)
takes points in the boundary of \(T\) to the boundary of \(\tau\),
points in \(\inter T\) to \(\inter \tau\), and the centerline of
\(T\) to the center point of \(\tau\).  From this it follows that
the orbit surface of a Seifert fibered space is a 2-manifold with
boundary and that the projection map carries boundary to boundary
and interior to interior.  The topological type of the orbit surface
is an invariant of the fiber type of the Seifert fibered space.

If \(T\) is a fibered solid torus, and \(\tau\) is its orbit space,
then we have realized \(\tau\) as the result of a composition of two
qotient maps.  The first map is from \(T\) to \(D\).  This is not
``fiber preserving'' unless \(T\) is an ordinary solid torus.  The
second map is the quotient map under the action of \(\rho\).

There is a second way to realize \(\tau\) as the result of a
composition of two quotient maps.  In the second way, both maps
``preserve fibers.''  Define \(\rho\x1:D\x I\into D\x I\) by
\((\rho\x1)(x,t)=(\rho(x),t)\).  This induces a fiber preserving map
of order \(\mu\) on \(T\).  The quotient space \(\overline{T}\) of
this action is naturally fibered as a product \(\tau\x S^1\) (an
ordinary solid torus).  The quotient map \(q\) is a branched cover
with brach set the centerline of \(T\).  Off the branch set, \(q\)
is a \(\mu\)-fold covering projection.  We now obtain \(\tau\) by
projecting the product \(\overline{T}=\tau\x S^1\) onto the first
factor.

Let \(h_t(\tau)\), \(0\le t\le1\), be an isotopy of \(\tau\) in that
each \(h_t\) is an embedding from \(\tau\) into \(\tau\).  Further
require that each \(h_t\) fix the center point of \(\tau\).  We get
an isotopy \(h_t\x1:\overline{T}\into \overline{T}\) by setting
\((h_t\x1)(x,s)=(h_t(x),s)\) which fixes the centerline of
\(\overline{T}\) pointwise.  This is a fiber preserving isotopy in
that fibers are preserved for each \(t\).  For each \((x,s)\), we
have a path \((h_t\x1)(x,s)\) as \(t\) varies from \(0\) to
\(1\). If \(x\) is not the center of \(\tau\), then the path avoids
the center of \(\tau\x\{s\}\).  If \(x\) is the center of \(\tau\),
then the path is constant.  If \((y,s)\) represents a point in \(T\)
with \(q(y)=x\), then there is a lift \((h_{y,t}\x1)(y,s)\) to
\(D\x\{s\}\) starting at \((y,s)\).  This defines an isotopy of
\(T\) that preserves levels, preserves fibers and fixes the
centerline pointwise.  By its construction, it is constant on any
set \(p^{-1}(S)\) where \(h_t\) is constant on \(S\) and the map
\(p:T\into\tau\) is projection.

\begin{lemma}\mylabel{FSTdeformation} Let \(T\) be a fibered solid
torus and let \(p:T\into\tau\) be the projection of \(T\) to its
orbit surface.  Let \(\tau'\) be a disk in \(\tau\) that contains
the center point in its interior.  Let \(T'=p^{-1}(\tau')\).  Then
there is a fiber preserving homeomorphism from \(T\) to \(T'\) that
is fixed on a neighborhood of the centerline of \(T\).  Further, the
homeomorphism is isotopic to the identity on \(T\).  \end{lemma}

\begin{demo}{Proof} There is an isotopy connecting the identity of
\(\tau\) to a homeomorphism from \(\tau\) to \(\tau'\).  We can
require that the isotopy fix pointwise a neighborhood of center
point of \(\tau\).  Lifting this to \(T\) gives a fiber preserving
isotopy that fixes pointwise a neighborhood of the centerline, and
that connects the identity on \(T\) to a homeomorphism to \(T'\).
\end{demo}

\subsection{Invariance, invariants and transitivity of fibers}

We are ready to show the uniqueness of fiber solid torus
neighborhoods of fibers.

Let \(M\) be a Seifert fibered space and let \(H\) be a fiber in the
interior of \(M\).  Let \(N\) and \(N'\) be fibered solid tori
neighborhoods of \(H\).  There is a fibered solid torus neighborhood
\(N''\) of \(H\) small enough to lie in the interior of \(N\cap
N'\).  If \(\tau\) is the orbit surface of \(N\), then the image of
\(N''\) is a disk \(\tau''\) in \(\tau\) that contains the center
point of \(\tau\) in its interior.  By Lemma \ref{FSTdeformation},
we obtain a fiber preserving homeomorphism from \(N\) to \(N''\)
that is the identity on a neighborhood of the centerline.
Similarly, there is such a homeomorphism from \(N'\) to \(N''\).
Composing one with the inverse of the other gives a fiber preserving
homeomorphism from \(N\) to \(N'\) preserving \(H\) that is the
identity on a neighborhood of \(H\).  We have:

\begin{thm}\mylabel{UniqueFSTNbds} Let \(M\) be a Seifert fibered
space and let \(H\) be a fiber in the interior of \(M\).  Then any
two fibered solid torus neighborhoods of \(H\) are connected by a
fiber preserving homeomorphism that is the identity on a
neighborhood of \(H\).  \BlackBox\end{thm}

For fibers in the boundary of a Seifert fibered space, no such
theorem is necessary since fibered solid torus neighborhoods of such
fibers are ordinary solid tori by definition.

If \(M\) is a Seifert fibered space, we will want to gain
information about preimages of disks in its orbit surface.  The next
lemma is a building block.

\begin{lemma} Let \(T\) be fibered solid torus and let \(T'\) be an
ordinary solid torus.  If \(T\) and \(T'\) are sewn along saturated
annuli in their boundaries by a fiber preserving homeomorphism, then
the result is fiber homeomorphic to \(T\).  \end{lemma}

\begin{demo}{Proof} By Lemma \ref{AnnularSew} we only have to show
this for one particular sewing.  Let \(p:T\into\tau\) be the
projection of \(T\) to its orbit surface.  Let \(E\) be a disk in
\(\tau\) that misses the image of the centerpoint and that
intersects \(\bd\tau\) in an arc.  As argued in \S\ref{LocalModel},
\(p^{-1}(E)\) is an ordinary solid torus that is sewn to
\(p^{-1}(\tau-\inter E)\) along a saturated annulus.  But
\(p^{-1}(\tau-\inter E)\) is fiber homeomorphic to \(T\) and the
result follows.  \end{demo}

A similar agrument gives:

\begin{lemma} Let \(M\) be Seifert fibered space and let \(T'\) be
an ordinary solid torus.  If \(M\) and \(T'\) are sewn along
saturated annuli in their boundaries by a fiber preserving
homeomorphism, then the result is fiber homeomorphic to \(M\).
\BlackBox\end{lemma}

If \(M\) is a Seifert fibered space, and \(H\) is a fiber of \(M\),
then we say that \(H\) is an \defit{ordinary fiber} of \(M\) if it
has a fibered solid torus neighborhood that is an ordinary solid
torus.  If not, we say that \(H\) is an \defit{exceptional} fiber of
\(M\).  Since every fiber of a fibered solid torus other than the
centerline is ordinary, it follows that the exceptional fibers are
isolated in \(M\).  We also know that the exceptional fibers are
located in \(\inter M\).  From Theorem \ref{UniqueFSTNbds}, we know
that the numerical invariant \(\nu/\mu\) of a fibered solid torus
neighborhood of a fiber \(H\) of \(M\) is also an invariant of
\(H\).

The denominator of \(\nu/\mu\), if the fraction is in reduced terms,
gives the intersection number of a meridian and a fiber in the
boundary of a fibered solid torus neighborhood of \(H\).  Thus a
fiber in the boundary of a fibered solid torus neighborhod
represents \(\mu\) times a generator of \(\pi_1\) of the
neighborhood, while \(H\) represents a generator.  Thus fibers close
to \(H\) are ``\(\mu\) times as long'' as \(H\).  We say that
\(\mu\) is the \defit{index} (or \defit{order}) of \(H\).

Let \(M\) be a Seifert fibered space, and let \(G\) be its orbit
surface.  A point in \(G\) is called \defit{exceptional} if it is
the image of an exceptional fiber, and it is called \defit{ordinary}
if not.  An exceptional point has a neighborhood in which it is the
only exceptional point.  Thus the exceptional points are isolated in
\(G\).  Associated with an exceptional point is the invariant
\(\nu/\mu\) of its corresponding fiber.  The topological type of
\(G\), the number of exceptional points of \(G\) and the collection
of invariants (with multiplicities) \(\nu/\mu\) associated with the
exceptional points are all invariants of the fiber type of \(M\).
These are not enough however to classify \(M\), and we will have to
further investigate the structure of \(M\).

We return to the task of extracting information from the orbit
surface.

\begin{thm} Let \(M\) be a Seifert fibered space and let \(G\) be
the orbit surface of \(M\).  Let \(H\) be a fiber in \(M\), let
\(h\) be its image in \(G\), and let \(E\) be a disk in \(G\)
containing \(h\) in its interior.  If \(E\) contains no exceptional
point except possibly for \(h\), then the preimage of \(E\) in \(M\)
is a fibered solid torus neighborhood of \(H\).  \end{thm}

\begin{demo}{Proof} By Theorem \ref{UniqueFSTNbds}, the theorem is
true if there is a fibered solid torus neighborhood of \(H\) which
contains the preimage of \(E\).  If not, then \(E\) can be
triangulated so that \(h\) is in the interior of some 2-cell \(E'\)
of the triangulation, and so that each 2-cell of the triangulation
is in the orbit disk of some fibered solid torus in \(M\).  Then the
preimage of each 2-cell other than \(E'\) is an ordinary solid
torus, and the preimage of \(E'\) is a fibered solid torus
neighborhood of \(H\).  If there is an ordering of the 2-cells of
the triangulation with \(E'\) as the first so that each 2-cell
intersects the union of its predecessors in an arc, then the result
follows from repetitions of Lemma \ref{AnnularSew}.  However, such
triangulations of arbitrarily small mesh are easy to find.
\end{demo}

The next result shows that ordinary fibers in a Seifert fibered
space are all pretty much the same.  We say that an isotopy of a
space \(X\) carries one subset \(S_0\) to another \(S_1\) if there
is an isotopy \(\Phi:X\x I\into X\) for which \(\Phi_0\) is the
identity on \(X\) and for which \(\Phi_1(S_0)=S_1\).  Note that if
\(M\) is a Seifert fibered space, \(p:M\into G\) is the projection
to the orbit surface and \(\Phi:M\x I\into M\) is a fiber preserving
isotopy of \(M\), then \(p\Phi:G\x I\into G\) defined by
\(p\Phi(x,t)=p\Phi(p^{-1}(x),t)\) is a well defined isotopy of
\(G\).

\begin{thm}\mylabel{IsotopyExist} Let \(M\) be a connected Seifert
fibered space and let \(H\) and \(H'\) be ordinary fibers in \(M\).
Then there is a fiber preserving isotopy \(\Phi\) of \(M\) that
carries \(H\) to \(H'\).  If \(p:M\into G\) is the projection to the
orbit surface and \(\alpha\) is a path from \(p(H)\) to \(p(H')\)
that avoids all exceptional points, then \(p\Phi\) can be required
to carry \(p(H)\) to \(p(H')\) along \(\alpha\).  If \(T\) and
\(T'\) are fibered solid torus neighborhoods of \(H\) and \(H'\)
respectively, then \(\Phi\) can be required to carry \(T\) to
\(T'\).  \end{thm}

\begin{demo}{Proof} An isotopy of the surface carrying \(p(T)\) to
\(p(T')\) while carrying \(p(H)\) to \(p(H')\) along \(\alpha\) can
be contructed in a finite number of stages so that each stage move
points only in a disk that contains no exceptional points.  (For
example, the first and last stages can be based on isotopies that
sqeeze \(T\) and \(T'\) to very thin fibered solid torus
neighborhoods of \(H\) and \(H'\) repsectively and the other stages
can be restricted to solid tori over a chain of small disks that
cover the path \(\alpha\).)  The preimage of such a disk is an
ordinary solid torus and the result follows by lifting these
isotopies to \(M\) by the techniques of \S\ref{FSTOrbit}.
\end{demo}

\begin{cor}\mylabel{UniqueDrill} Let \(M\) be a connected Seifert
fibered space and let \(H\) and \(H'\) be ordinary fibers in \(M\)
with fibered solid torus neighborhoods \(T\) and \(T'\)
respectively.  Then there is a fiber preserving homeomorphism from
\(M-\inter T\) to \(M-\inter T'\).  If \(M\) is orientable, then the
homeomorphism can be chosen to be orientation preserving.
\BlackBox\end{cor}

Uniqueness of fibered solid torus neighborhoods can be strengthened.

\begin{thm} Let \(M\) be a Seifert fibered space and let \(H\) be a
fiber of \(M\).  If \(N\) and \(N'\) are fibered solid torus
neighborhoods of \(H\), then there is a fiber preserving isotopy of
\(M\) fixed on a neighborhood of \(H\) that carries \(N\) to \(N'\).
\end{thm}

\begin{demo}{Proof} This follows from a similar statement that
assumes that \(N'\) is contained in \(N\) so we make that
assumption.  The image of \(N\) in the orbit surface \(G\) of \(M\)
is a disk \(E\) containing no exceptional point except possible the
image \(h\) of \(H\).  There is a disk neighborhood \(E''\) of \(E\)
with the same property.  The preimage of \(E''\) is another fibered
solid torus neighborhood of \(H\).  An isotopy of \(G\) fixed on a
neighborhood of the image of \(H\) and off \(E''\) carries \(E\) to
\(E'\) the image of \(N'\).  A lift of this isotopy is the desired
isotopy.  \end{demo}

\begin{cor} Let \(M\) be a Seifert fibered space and let \(H\) be a
fiber of \(M\).  If \(N\) and \(N'\) are fibered solid torus
neighborhoods of \(H\), then \(M-\inter N\) and \(M-\inter N'\) are
fiber homeomorphic.  \BlackBox\end{cor}

\subsection{Properties of the orbit surface}

We give some remarks connecting the topology of \(M\) and its orbit
surface \(G\).  If \(M\) is connected, we know that \(G\) is
connected.  The converse is also true.  Paths in \(G\) obviously
lift locally to \(M\), so the quotient map has (non-unique) path
lifting.  Now if \(G\) is connected, then any two fibers of \(M\)
are connected by a path, and \(M\) is connected.  If \(M\) is
compact, then so is \(G\).  If \(G\) is compact, then the images of
a finite number of neighborhoods of fibers of \(M\) cover \(G\).
Thus \(M\) is the union of a finite number of fibered solid tori and
is compact.  If \(M\) is compact, then it has only a finite number
of exceptional fibers and \(G\) has only a finite number of
exceptional points.  From previous remarks, we have that \(G\) has
no boundary if and only if \(M\) has no boundary.  Combining all
three observations, \(M\) is a closed, connected 3-manifold if and
ony if \(G\) is a closed, connected surface.  We add an observation
about \(\pi_1\).

\begin{lemma} Let \(M\) be a connected Seifert fibered space, let
\(G\) be the orbit surface of \(M\), and let \(p\) be the projection
map between them.  Then \(p\) induces a surjection on \(\pi_1\).
\end{lemma}

\begin{demo}{Proof} Paths in \(G\) lift to paths in \(M\), so loops
in \(G\) lift to paths in \(M\) that start and end in the same
fiber.  Such a path can be completed in the fiber to cover the
original loop in \(G\).  \end{demo}

A consequence of this lemma is that any Seifert fibering of \(S^3\)
must have \(S^2\) as the orbit space.  A Seifert fibering of a
closed, connected 3-manifold with finite fundamental group must have
\(S^2\) or \(P^2\) as the orbit space.

\section{Examples}

We first consider lens spaces.  In dimension 3, a lens space can be
viewed as a union of two solid tori sewn together along their
boundaries by a homeomorphism.  We first look at lens spaces without
considering fiber structures.

Let \(T_1\) and \(T_2\) be the two solid tori that make up a lens
space \(L\), and think of the sewing of their boundaries as
determined by a homeomorphism \(h\) from \(\bd T_1\) to \(\bd T_2\).
Once generators for the homologies of the boundaries are chosen,
\(h\) can be given as a \(2\x2\) integer matrix of determinant
\(\pm1\) since the homeomorphism type of \(L\) is determined if
\(h\) is determined up to isotopy.  We chose the generators to be
meridian-longitude pairs \((m_i,l_i)\) for each torus \(T_i\).  The
matrix \(A=\mtbt{q}{r}{p}{s}\) of \(h\) states that the image of
\(m_1\) is \(qm_2+pl_2\) and the image of \(l_1\) is \(rm_2+sl_2\).
The choice of letters is made to agree with standard conventions.

If \(h_1\) is the restriction to \(\bd T_1\) of a self homeomorphism
of \(T_1\) and \(h_2\) is a restriction to \(\bd T_2\) of a self
homeomorphism of \(T_2\), then \(h_2hh_1\) determines the same space
\(L\) as \(h\).  Self homeomorphisms of \(T_i\) exist that reverse
the meridians and fix the longitudes, reverse the longitudes and fix
the meridians, or reverse both.  Thus we may multiply any row or
column of \(A\) by \(-1\).  This allows us to negate any two entries
of \(A\) simultaneously.  (A diagonal is negated by negating one row
and one column.)  We may thus assume that \(p\ge0\) and \(|A|=1\).
There are self homeomorphisms of \(T_1\) (twists) that fix the
meridian and add integer multiples of the meridian to the longitude.
A matrix for this is \(\mtbt{1}{b}{0}{1}\) which combined with \(A\)
gives \[\mtbt{q}{r}{p}{s}\mtbt{1}{b}{0}{1}=
\mtbt{q}{bq+r}{p}{bp+s}.\] However, the right column of the result
gives all possible right columns which combine with \(\cv{q}{p}\) to
give determinant 1.  Thus all possible matrices of determinant 1
with left column \(\cv{q}{p}\) determine the same space which we
refer to as \(L_{p,q}\).  Similar self homeomorphisms of \(T_2\)
combined with \(A\) give \[\mtbt{1}{b}{0}{1}\mtbt{q}{r}{p}{s}=
\mtbt{q+bp}{r+bs}{p}{s}.\] Thus \(q\) need only be specified mod
\(p\).  Since the main diagonal can be negated without violating any
of our conventions, \(q\) need only be specified up to sign.
Lastly, we can think of the sewing as being accomplished by
attaching \(T_2\) to \(T_1\) instead of the other way around.  The
matrix of this is \(A^{-1}=\mtbt{s}{-r}{-p}{q}\), which can be
replaced by \(\mtbt{s}{r}{p}{q}\) by negating the secondary
diagonal.  Thus \(q\) can be replaced by \(s\) which is determined
up to multiples of \(p\) by the equivalence \(qs\equiv1\mod p\).
Since \(s\) can be negated, it can be specified by
\(qs\equiv\pm1\mod p\).

Since \(\cv{q}{p}\) is the image of \(m_1\) which is null-homotopic
in \(T_1\), the curve \(qm_2+pl_2\) in \(\bd T_2\) becomes trivial
in \(L_{p,q}\).  However this curve represents \(p\) times a
generator of \(\pi_1(T_2)\).  Adding \(T_1\) to \(T_2\) is
accomplished by adding a 2-handle which is a neighborhood of a
meridinal disk, and then adding a 3-handle.  Thus
\(\pi_1(L_{p,q})=\ints_p\) where we use \(\ints_0=\ints\) and
\(\ints_1=\{1\}\).  Thus \(p\) cannot be altered without changing
\(L_{p,q}\).  We have established all of one direction and part of
the converse of the following.  The rest of the converse can be done
by calculations of Whitehead torsion.

\begin{thm} The lens spaces \(L_{p,q}\) and \(L_{p',q'}\) are
homeomorphic if and only if \(p'=p\) and \(q'=\pm q^{\pm1}\mod
p\). \BlackBox\end{thm}

We now consider fiber structures.  We can get a Seifert fibered
structure on \(L_{p,q}\) by putting fiber structures on \(T_1\) and
\(T_2\) that are compatible with the sewing homeomorphism.  If
\(T_1\) is a fibered solid torus determined by \(\nu/\mu\), then a
fiber on \(\bd T_1\) represents the element \(\cv{\nu}{\mu}\), and
under the action of \(A\), will represent the element
\[\cv{\nu'}{\mu'}= \mtbt{q}{r}{p}{s}\cv{\nu}{\mu}=
\cv{q\nu+r\mu}{p\nu+s\mu}.\] This fibers \(L_{p,q}\) with (possibly)
two exceptional fibers of indices \(\mu\) and \(\mu'=p\nu+s\mu\).
Since \(\cv{r}{s}\) is determined only up to added multiples of
\(\cv{q}{p}\), there is considerable freedom in choosing \(\nu/\mu\)
and \(\nu'/\mu'\).  We thus get infinitely many examples of Seifert
fibered spaces that are homeomorphic but not under fiber preserving
homeomorphisms.  If \(\mu\) or \(\mu'\) is 1, then there will be
fewer than two exceptional fibers.  If \(\mu=1\), then \(\nu=0\) and
\(\mu'=s\).  We know that the possibilities for \(s\) are given by
\(qs\equiv\pm1\mod p\), so that any lens space can be fibered in
infinitely many ways with one exceptional fiber.  The lens spaces
\(L_{p,1}\) can be fibered (among other ways) with no exceptional
fibers.

The space \(S^3\) is obtained from the sewing
\(\mtbt{0}{-1}{1}{0}\).  By our remarks above, \(S^3\) is also
obtained from the sewings \(\mtbt{q}{sq-1}{1}{s}\) for any \(q\) and
\(s\).  Here \(\cv{\nu'}{\mu'}=\cv{q\nu+(sq-1)\mu}{\nu+s\mu}\) which
can be replaced by \(\cv{-\mu}{\nu+s\mu}\) since \(\nu'/\mu'\) is
only determined mod 1.  Since \(s\) is arbitrary, and \(\nu\) is
only required to be relatively prime to \(\mu\), it is seen that
\(S^3\) can be fibered with two exceptional fibers of any two
indices as long as they are relatively prime, with one exceptional
fiber of any index, or with no exceptional fibers.  We will see
later that there are no other possibilities.  (Actually, it seems
that we won't.)

The space \(S^2\x S^1\) is obtained from the sewing
\(\mtbt{1}{0}{0}{1}\) and in fact from any of the sewings
\(\mtbt{1}{r}{0}{1}\).  This gives
\(\cv{\nu'}{\mu'}=\cv{\nu+r\mu}{\mu}\) which can be replaced by
\(\cv{\nu}{\mu}\).  Thus \(S^2\x S^1\) can be fibered with two
exceptional fibers having identical invariants or with no
exceptional fibers.

Note that all lens spaces have an orbit surface that is a union of
two disks sewn along the boundaries.  Thus the orbit surface is
\(S^2\).  If \(G\) is any closed surface, then \(G\x S^1\) is a
Seifert fibered space with no exceptional fibers and with orbit
surface \(G\).  Note that \(G\x S^1\) is non-orientable if \(G\) is
non-oritentable.  A Seifert fibered space need not be orientable
even though each fiber has an orientable neighborhood and is thus an
orientation preserving curve.  The Seifert fibered space might be
orientable even if the orbit surface is not.  If \(G\) is a
non-orientable suface, then there is an orientable circle bundle
over \(G\).  It can be realized as the double of the orientable
interval bundle over \(G\).  It contains a copy of the orientable
double cover of \(G\) (the boundary of the orientable interval
bundle over \(G\)) as a separating surface.  If \(G\) is \(P^2\),
then this separating surface is \(S^2\) and the complementary
domains are orientable line bundles over \(P^2\) which are copies of
\(P^3\) minus a point.  Thus the orientable circle bundle over
\(P^2\) is the connected sum of two copies of \(P^3\).  We will see
that this is the only Seifert fibered space that is not prime.  The
only other Seifert fibered space that is not irreducible is \(S^2\x
S^1\).

If \(M\) is the orientable circle bundle over the Klein bottle, then
an alternate fibering exists for \(M\).  The Klein bottle is
realized as a quotient of \(S^1\x I\) in which reflection of \(S^1\)
across an axis is used to sew the 0 level to the 1 level.  The
orientable circle bundle is realized as a quotient of \(S^1\x S^1\x
I\) in which reflection of each \(S^1\) factor across an axis is
used to sew the 0 level to the 1 level.  Reflection one \(S^1\)
reverses orientation, while reflecting two preserves orientation.
If we think of \(S^1\) as the unit circle in \(\comps\), then
reflection can be realized by complex conjugation.  Thus each factor
has \(1\) and \(-1\) as two fixed points, and the action on \(S^1\x
S^1\) has four fixed points.  All other points have order two and
\(M\) is seen to have a Seifert fibration with four exceptional
fibers each with index 2.  The orbit surface has \(S^1\x S^1\) as a
2-fold branched cover with four branch points.  If the orbit surface
is decomposed as a cell complex using the four exceptional points as
vertices and using \(e\) edges and \(f\) faces, then the branched
cover will have four vertices, \(2e\) edges and \(2f\) faces.  Since
the Euler characteristic of a torus is 0, we have \(4-2e+2f=0\) and
\(-e+f=-2\).  Thus the orbit surface has Euler characteristic
\(4-e+f=2\) and is \(S^2\).  Thus \(M\) fibers over \(S^2\) with
four exceptional fibers of index 2, and it also fibers over the
Klein bottle with no exceptional fibers.

\section{Fiber structures of compact Seifert fibered spaces}

We now restrict attention to Seifert fibered spaces that are
connected, compact manifolds.  Such spaces have a finite number of
exceptional fibers, and a finite number of boundary components, each
of which is a compact fibered surface and is thus a torus or a Klein
bottle.  The orbit surfaces of such spaces are connected, compact
surfaces.

Recall that a compact, connected surface is bounded by a finite
number (possibly zero) of circles.  Such a surface is associated
with a unique closed surface that is obtained by sewing disks to
each boundary component.  Two compact, connected surfaces are said
to be \defit{of the same genus} if they yield homeomorphic closed
surface when disks are added to all boundary components.  If the
associated closed surface is orientable, then the genus is usually
identified by the number of handles \([1-(\chi/2)]\) where \(\chi\)
is the Euler characteristic of the closed surface.  If the closed
surface is non-orientable, then the genus is usually identified with
the number of projective plane summands, called the number of
crosscaps (which are actually \Mob{}s), in the surface \([2-\chi]\).
A compact, connected surface is characterized by its orientability,
its genus (expressed as either a number of handles or as a number of
crosscaps), and the number of its boundary components.

Let \(M\) be a connected, compact Seifert fibered space.  There are
only a finite number of exceptional fibers in \(M\) and these are in
\(\inter M\).  Let \(H\) be an exceptional fiber in \(M\), let \(N\)
be a fibered solid torus neighborhood of \(H\), and let \(J\) be a
meridian of \(N\).  If we remove the interior of \(N\), then a
connected, compact Seifert fibered space \(M_0\) with one more torus
boundary component than \(M\) results.  We can create another
Seifert fibered space \(M'\) by sewing in another fibered solid
torus \(N'\) using a fiber preserving homeomorphism from \(\bd N'\)
to the new torus boundary component \(C\) of \(M_0\).  The space
\(M'\) is completely determined if a curve on \(C\) is picked out to
be a meridian of \(N'\).  We are only required to choose the curve
to be non trivial in the homology of \(C\) and not homologous in
\(C\) to a fiber.  If the curve is chosen to be a crossing curve for
\(C\), then \(N'\) will be an ordinary solid torus, and \(M'\) will
have one fewer exceptional fiber than \(M\) and the same number of
boundary components.  (Obviously, a Seifert fibered space with no
exceptional fibers can be obtained if this step is repeated a finite
number of times.)  We can recover \(M\) from \(M'\) by removing the
interior of \(N'\) and replacing \(N\).  All that we have to do to
correctly replace \(N\) is remember the location of \(J\).  If we do
not remember where \(J\) was, then we need a technique to recover
it.

This establishes an outline.  We will build Seifert fibered spaces
by adding exceptional fibers one at a time, starting with a space
with no exceptional fibers.  We will do this by removing the
interiors of ordinary solid tori, and then replacing them by fibered
solid tori.  We will need to study the process of removing and
replacing fibered solid tori, we will need to study Seifert fibered
spaces having no exceptional fibers, and we will need to establish a
mechanism for locating meridians of removed fibered solid tori.

\subsection{Drilling and filling}

Let \(M\) be a Seifert fibered space, let \(H\) be a fiber in the
interior of \(M\), and let \(N\) be a fibered solid torus
neighborhood of \(H\).  We know that the fiber type of \(M-\inter
N\) depends only on \(H\).  We say that \(M-\inter N\) is obtained
from \(M\) by \defit{drilling out} the fiber \(H\).  Note that when
\(M_0\) is obtained from \(M\) by drilling out an ordinary fiber,
then the fiber type of \(M_0\) is independent of the specific fiber.

Let \(M\) be a Seifert fibered space, let \(C\) be a torus boundary
component of \(M\), and let \(J\) be a simple closed curve in \(C\)
that is not null homologous and is not homologous in \(C\) to a
fiber.  We know that a fibered solid torus can be sewn to \(M\)
along \(C\) so that \(J\) becomes a meridian for the fibered solid
torus and we know that the fiber type of the result \(M'\) depends
only on \(J\).  We say that \(M'\) is obtained from \(M\) by
\defit{filling} the torus boundary component \(C\).

Note that drilling removes a disk from the orbit surface of a
Seifert fibered space and filling adds a disk.  Thus if a Seifert
fibered space is altered by an equal numbers of drillings and
fillings, then the orbit surface of the result is homeomorphic to
the orbit surface of the original.  If an unequal number of
drillings and fillings is done, then the genus of the orbit surface
remains the same.  Also, drilling and filling will not change the
number of boundary components that are Klein bottles.  There is more
that is invariant.  To see this, we need to define another
structure.

\subsection{Classes preserved by drilling and filling}

The next lemma is the basis of this section.  Recall that if
\(\Phi\) is a fiber preserving isotopy of a Seifert fibered space
\(M\) and \(p:M\into G\) is the projection to the orbit surface,
then \(p\Phi\) is the induced isotopy of \(G\).  The action of
\(\Phi\) on a fiber \(H\) of \(M\) gives rise to a path
\((p\Phi_t)|_{p(H)}\) in \(G\).

\begin{lemma} Let \(M\) be a Seifert fibered space, let \(H\) be an
ordinary fiber of \(M\) in \(\inter M\), let \(p:M\into G\) be the
projection to the orbit surface, and let \(\Phi\) be a fiber
preserving isotopy of \(M\) that carries \(H\) to a fiber \(H'\) in
\(M\).  Then the unpointed homotopy class of the map
\(\Phi_1|_H:H\into H'\) is determined by the homotopy class of the
path \((p\Phi_t)|_{p(H)}\) rel its endpoints.  \end{lemma}

\begin{demo}{Proof} Since the maps on the fibers are homeomorphisms,
there are only two unpointed homotopy classes and these are
determined by the effect of the maps on the orientations of the
fibers.  If counterexamples exist (two isototpies that induce the
same path starting at some fiber \(H\), but whose 1 levels give maps
on \(H\) that are not homotopic), then following one isotopy with
the reverse (in \(t\)) of the second isotopy would give an example
of a fiber preserving isotopy that carries \(H\) to itself in an
orientation reversing way along a loop in \(G\) that is null
homotopic in \(G\).  In particular this would demonstrate that the
element of \(\pi_1(M)\) represented by the fiber would be conjugate
to its inverse.  This cannot happen if \(M\) is a fibered solid
torus since \(\pi_1\) of a fibered solid torus is \(\ints\) in which
a fiber is a non-trivial element and is not conjugate to its
inverse.  The easiest way to enlarge this observation to an
arbitrary \(M\) without heavy details is to view an isotopy as a
path in a function space.

We are only concerned with the restriction of the isotopy to the
fiber \(H\).  We thus want to look at the space \(A\) of all
embeddings of \(H\) into fibers of \(M\).  Since we are only
concerned with the homotopy class of a given embedding, in which any
element in a class can be obtained from any other by composing on
the right with an orientation preserving self homeomorphism of
\(H\), we can let \(B\) be the space of all orientation preserving
self homeomorphisms of \(H\) and let it act on \(A\) on the right.
The isotopies that we are concerned with induce paths in the space
\(A/B\).  The space \(A/B\) is a double cover of \(G\) since it has
two points for each fiber in \(M\).  We want to know if paths in
\(G\) lift uniquely to \(A/B\).  However, our observation that this
is the case if \(M\) is a fibered solid torus says that lifting is
locally unique, therefore it is globally unique.  \end{demo}

Using the lemma above, we can construct a well defined homomorphism
from \(\pi_1(G,x)\) to \(\ints_2\) (regarded as \(\{-1,1\}\) under
multiplication) where \(x\) is an ordinary point of \(\inter G\).
If \(\alpha\) is a loop (which can be taken to avoid exceptional
points), then we take \(\alpha\) to \(+1\) if an isotopy carrying
\(H\), the fiber above \(x\), to itself around \(\alpha\) maps \(H\)
to itself in an orientation preserving way, and to \(-1\) if it is
carried in an orientation reversing way.  The lemma shows that this
is well defined, and it is clearly a homomorphism.  Its kernel is
the image of \(\pi_1\) of the space \(A/B\) of the proof of the
lemma since a loop lifts to a loop in \(A/B\) if and only the loop
is taken to \(+1\).  We denote this homomorphism by \(\phi\) and
call it the \defit{classifying homomorphism} of \(M\).

We say that two such homomorphisms \(\phi_1:\pi_1(G_1)\into\ints_2\)
and \(\phi_2:\pi_1(G_2)\into\ints_2\) are \defit{equivalent} if
there is a homeomorphism \(h:G_1\into G_2\) so that
\(h_{\#}:\pi_1(G_1)\into\pi_1(G_2)\) gives \(\phi_1=\phi_2h_{\#}\).
Fiber homeomorphic Seifert fibered spaces have equivalent
classifying homomorphisms.

The classifiying homomorphism distinguishes the torus boundary
components from the Klein bottle boundary components.  If \(C\) is a
compact boundary component of a Seifert fibered space \(M\), then
its image on the orbit surface is a boundary circle.  If a loop is
conjugate to this circle then its image under the classifying map is
\(+1\) if \(C\) is a torus, and \(-1\) if \(C\) is a Klein bottle.

\begin{lemma} Let \(M\) be a Seifert fibered space and let \(M'\) be
obtained from \(M\) by an equal number of drillings and fillings.
Then \(M\) and \(M'\) have equivalent classifying homomorphisms.
\end{lemma}

\begin{demo}{Proof} If the number of drillings and fillings is zero,
then there is nothing to prove.  Otherwise there is a consecutive
pair of operations one of which is a drilling and the other of which
is a filling.  We first consider the case where the filling comes
first.

There are three spaces involved.  Space 1 exists just before the
filling, space 2 after the filling and space 3 after the drilling.
Since exceptional fibers may be involved, there may not be an easy
way to relate the three Seifert fibered spaces, but we can relate
the three orbit surfaces.  Let these be denoted \(G_1\), \(G_2\) and
\(G_3\) in order.  Let \(E_1\) be the disk in \(G_2\) that is the
disk added during the filling, and let \(E_2\) be the disk that is
removed during the drilling.  There is an isotopy of \(G_2\)
carrying \(E_1\) to \(E_2\) so we get a homeomorphism \(f\) from
\(G_1=G_2-\inter E_1\) to \(G_3=G_2-\inter E_2\) that extends to a
homeomorphism of \(G_2\) that is isotopic to the indentity.

If \(l\) is a loop in \(G_1\), then \(fl\) is a loop in \(G_3\) and
we wish to show that the classifying homomeophisms on \(G_1\) and
\(G_3\) take the same values on these two loops.  But the
classifying homomorphism on \(G_2\) does give these two loops the
same value.  The value of \(l\) on \(G_1\) agrees with that on
\(G_2\) since the only change is filling a hole whose boundary had
value \(+1\).  Similarly the value of \(fl\) on \(G_2\) agrees with
that on \(G_3\) and the classifying homomorphisms on \(G_1\) and
\(G_3\) are equivalent.

It the drilling comes first, then the analysis is similar except
that now \(G_2\) has two boundary components \(J\) and \(K\)
(perhaps the same), one filled with a disk to create \(G_1\) and the
other filled with a disk to create \(G_3\).  The classifying
homomorphisms must take the value \(+1\) on both \(J\) and \(K\)
since the boundary components over \(J\) and \(K\) are involved in
drilling and filling and must be tori.  We can form \(G'_2\) by
attaching disks to both \(J\) and \(K\) (if \(J=K\) only one disk is
attached) and extending the classifying homomorphism.  There is an
ambient isotopy of \(G'_2\) that carries one attached disk to the
other and the remaining details are left to the reader.

The general induction step now involves surfaces \(G'_1\), \(G_1\),
\(G_2\) and \(G'_2\) with classifying homomorphisms where the
homomorphisms on \(G_1\) and \(G_2\) are known to be equivalent and
each \(G'_i\) is obtained from \(G_i\) by removing a disk or each
\(G'_i\) is obtained from \(G_i\) by adding a disk.  The details are
left to the reader.  \end{demo}

There are restrictions on the values of a classifying homomorphism
on boundary components.  Compact surfaces are built from punctured
disks, once punctured tori and \Mob{}s.  The boundary of a once
punctured torus is a commutator and the boundary of a \Mob is a
square.  Thus these curves must be taken to \(+1\) by a classifying
homomorphism.  The curves on the boundary of a punctured disk are
related by having their product equal 1.  Thus it cannot be that an
odd number of them is taken to \(-1\) by a classifying homomorphism.
A compact surface can be regarded as a punctured disk which has some
of its boundary components sewn to once punctured tori and \Mob{}s.
The curves that are used for attaching have value \(+1\) under a
classifying homomorphism and an even number of the remaining curves
(possibly zero) have value \(-1\).  Thus the number of Klein bottle
boundary components of a compact Seifert fibered space is even.

Let \(M\) be a connected, compact Seifert fibered space.  We define
the \defit{class} of \(M\) to be the class of all closed Seifert
fibered spaces that can be obtained from \(M\) by a finite sequence
of drillings, fillings and fiber preserving homeomorphisms.  We
would like to pick out a well defined representative of a class.
Since all exceptional fibers can be drilled, we can always get a
space with no exceptional fibers.  Since all torus boundary
components can be filled, we can always get a space with no torus
boundary components.  Unfortunately, a class can have more than one
fiber type of space with no exceptional fibers and no torus boundary
components.  However, we will see that this is a problem only of
closed manifolds.  As soon as boundary is introduced, ambiguity goes
away.  We therefore define (the) \defit{classifying space} of a
class to be (the) representative with no exceptional fibers and
exactly one torus boundary component.  We could do without the torus
boundary component if there are Klein bottle boundary components
present, but there seems to be no gain by doing that.  Uniqueness
must be established.  We need some techniques.

Let \(G\) be a compact, conected surface with non-empty boundary and
let \(\phi\) be a homomorphism from \(\pi_1(G)\) to \(\ints_2\).
The basepoint is no problem because \(\ints_2\) is abelian and
\(\phi\) is impervious to change by conjugation.  There are a finite
number of pairwise disjoint arcs \(\{\alpha_i\}\) in \(G\) with
boundaries in the boundary of \(G\) so that cutting \(G\) along
these arcs yields a disk \(E\).  The surface \(G\) is recovered by
sewing pairs of arcs together in the boundary of \(E\).  These arcs
\(\{\alpha'_i, \alpha''_i\}\) in \(\bd E\) are pairwise disjoint and
have two copies \(\alpha'_i\), \(\alpha''_i\) corresponding to each
\(\alpha_i\).  Let \(M_0\) be a Seifert fibered space with no
exceptional fibers and with orbit surface \(E\).  We know that
\(M_0\) must be an ordinary solid torus.  We fix one fiber in \(\bd
M_0\) as a longitude with a given orientation for future reference.
This orients all of the fibers in \(M_0\) by requiring that they all
represent the same homology class as the longitude.  The preimages
of the \(\alpha'_i\) and \(\alpha''_i\) in \(M_0\) are a set of
pairwise disjoint saturated annuli in \(\bd M_0\).  We recover a
Seifert fibered space whose orbit surface is \(G\) if we sew these
annuli together in pairs with some care.  With some extra care, we
also get a space whose classifying homomorphism is \(\phi\).  The
sewings have to be with fiber preserving homeomorphisms.  The care
needed to recover \(G\) is to see that the sewing of the annuli
reflects the sewing of the \(\alpha'_i\) to the \(\alpha''_i\).
This is accomplished by seeing that the right pairs of annuli are
sewn together, and seeing that orientations of the \(I\) factors in
the annuli are handled to reflect the orientations of the sewings of
the \(\alpha'_i\) to the \(\alpha''_i\).  This requires that the
boundaries of the annuli be matched correctly.  To recover \(\phi\),
we need to handle the orientations of the \(S^1\) factors correctly.
This is done by finding a simple closed curve \(J\) in \(G\) that
pierces a given \(\alpha_i\) once and noting the value of \(\phi\)
on \(J\).  If the value is \(+1\), then the orientations of the
\(S^1\) fibers in the appropriate annuli are to be sewn so as to
match the orientations as inherited from the fixed longitude of
\(M_0\).  If the value is \(-1\), then the orientations of the
fibers in one annulus are to be matched to the negatives of the
orientations of the fibers in the other.  This will create a space
with the right value of the classifying homomorphism on \(J\).
Since a set of curves chosen one per \(\alpha_i\) that pierces
\(\alpha_i\) once is a set of generators for \(\pi_1(G)\), we have
reconstructed \(\phi\).

We call this construction a realization of \(\phi\) along the arcs
\(\{\alpha_i\}\).  Note that the space constructed has no
exceptional fibers.  Since all realizations start with an ordiary
solid torus, the only freedom in contructing the realization is in
the sewing maps.  However, the pairing of the annuli to be sewn is
determined by \(G\) and the set \(\{\alpha_i\}\), the treatment of
the orientations of the \(I\) factors is determined, and the
treatment of the \(S^1\) factors is determined.  By Lemma
\ref{AnnulusClasses}, the sewings used in one realization are fiber
isotopic to the sewings used any other realization.  Thus any two
realizations along \(\{\alpha_i\}\) are fiber homeomorphic.  It is
also clear that every compact, connected Seifert fibered space \(M\)
with non-empty boundary and an absence of exceptional fibers is a
realization because it is obtained from an ordinary solid torus by
sewings of certain annuli that are obtained as saturated sets over
an appropriate set of cutting arcs of the orbit surface of \(M\).

Let \(\{\beta_i\}\) be another set of arcs that cut \(G\) to a disk
\(E'\).  Let \(M\) be one realization along \(\{\alpha_i\}\).  The
preimages of the \(\beta_i\) are saturated annuli in \(M\) that cut
\(M\) to a Seifert fibered space \(M'_0\) with no exceptional fibers
with orbit surface \(E'\).  Thus \(M'_0\) is seen to be an ordinary
solid torus and \(M\) is seen to be a realization along
\(\{\beta_i\}\).  Thus realizations along \(\{\beta_i\}\) are fiber
homeomorphic to realizations along \(\{\alpha_i\}\).

Now let \(G\) and \(G'\) be surfaces possessing homomorphisms
\(\phi\) and \(\phi'\) respectively, and let \(h:G\into G'\) be a
homeomorphism that demonstrates that \(\phi\) and \(\phi'\) are
equivalent.  Let \(M\) realize \(\phi\) along arcs \(\{\alpha_i\}\)
in \(G\).  If \(p:M\into G\) is projection, then \(hp:M\into G'\) is
also a projection and \(M\) is seen also to be a realization of
\(\phi'\) along \(\{h(\alpha_i)\}\) in \(G'\).  Thus realizations of
\(\phi\) and \(\phi'\) are fiber homeomorphic.  We have shown the
following.

\begin{thm} (a) Let \(G\) be a compact, connected surface with
non-empty boundary and let \(\phi\) be a homomorphism from
\(\pi_1(G)\) to \(\ints_2\).  Then there is a Seifet fibered space
with orbit surface \(G\), with no exceptional fibers and with
classifying homomorphism \(\phi\).  Any two such are fiber
homeomorphic.

(b) Compact, connected Seifert fibered spaces with no exceptional
fibers, with non-empty boundary, and with equivalent classifying
homomorphisms are fiber homeomorphic.  \BlackBox\end{thm}

\begin{cor} Let \(M\) be a class of Seifert fibered spaces.  Then
any two classifying spaces for the class are fiber homeomorphic.
\end{cor}

\begin{demo}{Proof} Two such spaces are obtained from each other by
drillings and fillings.  Since the number of torus boundary
components of the two spaces is the same, the number of drillings
and fillings is equal.  Thus the two spaces have equivalent
classifying homomorphisms and are fiber homeomorphic.  \end{demo}

\begin{cor} Two compact, connected Seifert fibered spaces with the
same number of torus boundary components are in the same class if
and only if they have equivalent classifying homomorphisms.
\end{cor}

\begin{demo}{Proof} On each space, drillings and fillings replace
exceptional fibers by ordinary fibers.  For either space, the number
of drillings and fillings is the same and the classifying
homomorphisms stay in their equivalence classes.  If the resulting
spaces have boundary, we are done.  If not, then one drilling each
gives spaces with boundary and the classifying homomorphisms become
equivalent if and only if the originals were equivalent.  However,
the classifying homomorphisms become equivalent if and only if the
spaces become fiber homeomorphic.  \end{demo}

The spaces of the corollary above need not be fiber equivalent nor
even homeomorphic.  The examples of fiberings of lens spaces given
above are all in the same class.  Since the orbit surfaces are all
\(S^2\), which is simply connected, the classifying homomorphisms
are all trivial.

The next two corollaries are really corollaries of the realization
construction.

\begin{cor}\mylabel{ReverseFibers} Every compact, connected Seifert
fibered space with non-empty boundary and no exceptional fibers
admits a fiber preserving self homeomorphism that carries each fiber
to itself and reverses the orientations of all the fibers.
\end{cor}

\begin{demo}{Proof} The realization construction starts with a disk
\(E\) and an ordinary solid torus that has the structure of \(E\x
S^1\).  We think of \(S^1\) as the unit circle in \(\comps\).
Annuli \(\alpha'\x S^1\) are sewn to annuli \(\alpha''\x S^1\) and
the sewings are only specified by orientation requirements on the
arcs \(\alpha'\) and \(\alpha''\) and on the \(S^1\) factors.  We
have enough freedom to cover the possibilities if we insist that the
sewings in the \(S^1\) direction use the either identity on \(S^1\)
or complex conjugation.  This commutes with the self map of \(E\x
S^1\) that carries each \((e,x)\) to \((e,\overline{x})\).  This
induces the desired homeomorphism on the realization.  \end{demo}

If \(M\) is a Seifert fibered space, \(G\) its orbit surface and
\(p:M\into G\) the projection, then a \defit{section} for \(G\) is
an embedding \(\sigma:G\into M\) for which \(p\sigma\) is the
identity on \(G\).

\begin{cor}\mylabel{SectionExists} Every compact, connected Seifert
fibered space with non-empty boundary and no exceptional fibers has
a section for its orbit surface.  \end{cor}

\begin{demo}{Proof} As in the proof of the previous corollary, we
start with the ordinary solid torus \(E\x S^1\) and its section
embedding \(E\) into \(E\x\{1\}\).  Since fibers are sewn using the
identity or complex conjugation on the \(S^1\) factor, the sewing of
\(E\) to recreate \(G\) is reproduced by the image of \(E\x\{1\}\).
\end{demo}

We will concentrate on understanding closed, connected Seifert
fibered spaces in the next section.  The obvious question that
remains for this section is how many equivalence classes of
homomorphisms from \(\pi_1\) of a surface to \(\ints_2\) are there.
Before we answer that we will show that if we answer this question
for closed surfaces, then we will not only have enough information
to understand this question for arbitrary surfaces but we will also
have enough information for an understanding of compact, connected
Seifert fibered spaces modulo an understanding of closed, connected
Seifert fibered spaces.

Let \(G\) be a compact, connected surface with boundary and let
\(\phi\) be a homomorphism from \(\pi_1(G)\) to \(\ints_2\).  We
know that the number of boundary components on which \(\phi\) takes
the value \(-1\) is even.  Thus if \(G\) has one boundary component,
then \(\phi\) takes the value \(+1\) on this component and \(\phi\)
has a unique extention \(\hat{\phi}\) to the surface \(\hat{G}\)
which is obtained from \(G\) by sewing a disk onto the boundary
curve.  We know that restricting \(\hat{\phi}\) to any surface
obtained from \(\hat{G}\) by removing the interior of a disk results
in a homomorphism that is equivalent to \(\phi\).  Thus if \(\bd G\)
is connected, then the equivalence classes of homomorphisms from
\(\pi_1(G)\) to \(\ints_2\) are in one-to-one correspondence with
the equivalence classes of homomorphisms from \(\pi_1(\hat{G})\) to
\(\ints_2\).

If \(G\) has several boundary components, then let \(J\) be a
separating, simple closed curve on \(G\) that bounds a disk with
holes \(E\) containing all the boundary components of \(G\).  Let
\(\hat{G}\) be the result of sewing a disk to \(G-\tinter E\) along
\(J\).  A homomorphism from \(\pi_1(G)\) to \(\ints_2\) determines a
unique \(\hat{\phi}\) from \(\pi_1(\hat{G})\) to \(\ints_2\), and a
restriction to \(E\).  The homomorphism \(\phi\) can be
reconstructed from the restrictions of \(\hat{\phi}\) to \(G-\tinter
E\) and \(\phi\) to \(E\).  As mentioned above, the class of the
restriction of \(\hat{\phi}\) to \(G-\tinter E\) is determined by
\(\hat{\phi}\) and the fact that \(\hat{G}\) and \(G-\tinter E\)
differ by a disk.  The restriction of \(\phi\) to \(E\) is
determined by the number of boundary components of \(G\) that
\(\phi\) evaluates to \(-1\).  This is because homeomorphisms of
\(E\) can realize any permutation of the boundary components of
\(E\) keeping \(J\) fixed, and because \(\pi_1(E)\) is generated by
the boundary components of \(E\) with \(J\) omitted.  Thus the
equivalence classes of all homomorphisms from fundamental groups of
compact, connected surfaces to \(\ints_2\) are determined by those
from closed, connected surfaces, and from those on disks with holes.
The classes on closed surfaces will be discussed shortly, and the
classes on disks with holes are easily understood.

Now let \(M\) be a compact, connected Seifert fibered space with
non-empty boundary.  Assume that all exceptional fibers have been
drilled out and replaced with ordinary solid tori.  Let \(G\) be the
orbit surface, let \(p:M\into G\) be the projection, and let
\(\phi\) be the classifying homomorphism.  As before, let \(J\) be a
separating simple closed curve on \(G\) that bounds a disk with
holes \(E\) that contains all the boundary components of \(G\).  The
saturated torus \(T=p^{-1}(J)\) splits \(M\) into two Sefiert
fibered spaces \(M'\) with \(\bd M'=T\), and \(M''\) with \(\bd
M''=T\cup\bd M\).  We can sew an ordinary solid torus to \(M'\)
along \(T\) to create a closed, connected Seifert fibered space
\(\hat{M}\) with no exceptional fibers.  The classifying
homomorphism of \(\hat{M}\) is determined by \(\phi\).  The fiber
type of \(M'\) is determined by the classifying homomorphism of
\(\hat{M}\) and the fact that \(M'\) is obtained from \(\hat{M}\) by
drilling out a single fiber.

The space \(M''\) has \(E\) as its orbit surface and is determined
by the restriction of \(\phi\) to \(\pi_1(E)\) which in turn is
determined by the number of boundary curves of \(E\) that \(\phi\)
takes to \(-1\) which equals the number of boundary components of
\(M\) that are Klein bottles.  We recover \(M\) as \(M'\cup M''\).
If spaces fiber homeomorphic to \(M'\) and \(M''\) are supplied
instead (which we still refer to as \(M'\) and \(M''\)), then we can
try to recover a space fiber homeomorphic to \(M\) by sewing \(M'\)
to \(M''\) by a fiber preserving homeomorphism from one torus
boundary component of \(M''\) to the unique boundary component of
\(M'\).  No matter what fiber preserving homeomorphism is used, the
result of the sewing will have orbit surface homeomorphic to \(G\)
and the homeomorphism can be arranged to demonstrate that the
classifying homomorphism is equivalent to \(\phi\).  Thus we always
get a space fiber homeomorphic to \(M\).

We now consider the question of determining the equivalence classes
of homomorphisms from fundamental groups of closed, connected
surfaces to \(\ints_2\).  We will see that for orientable surfaces
there are 1 or 2 classes, and for non-orientable surfaces, there are
as many as four classes.  The fact that there are at least this many
classes is shown by constructing the homomorphisms and showing by
various arguments that they are not equivalent.  The fact that they
exhaust all the classes is shown by knowing enough homeomorphisms of
closed surfaces to reduce any homomorphism to one of the classes.
This last argument is done mostly by pictures.

Let \(G\) be a closed, connected, orientable surface with \(g\)
handles.  The fundamental group is generated by \(2g\) curves in
\(G\) with one point in common whose complement in \(G\) is an open
disk.  The only relation that applies to these generators is that a
certain product of commutators is the identity.  Thus any assignment
of \(\pm1\) to each curve results in a valid homomorphism to
\(\ints_2\).  There are two obviously different homomorphisms ---
one \(\phi_1\) that has image \(\{+1\}\) and another \(\phi_2\) that
has image \(\{+1,-1\}\).  The homomorphism \(\phi_1\) is trivial and
needs no more details.  A canonical example for \(\phi_2\) is to
define \(\phi_2\) to be \(-1\) on each of the \(2g\) generators
described above.  We will argue later that these represent all
classes.  Note that \(\phi_2\) does not apply to \(S^2\).

Let \(G\) now be a closed, connected, non-orientable surface with
\(k\) crosscaps.  Not all curves on \(G\) are alike.  Some preserve
orientation and some reverse orientation.  The subgroup \(O\) of
orientation preserving curves is of index 2 in \(\pi_1(G)\).  Two
homomorphisms that are immediate are the trivial \(\phi_1\) which
takes all of \(\pi_1(G)\) to \(+1\), and another \(\phi_3\) that has
image all of \(\ints_2\) with kernel \(O\).  That is, it takes every
orientation reversing curve to \(-1\).  We now find two more
examples.

Since \(\ints_2\) is abelian, a homomorphism from \(\pi_1(G)\)
factors through \(H_1(G)\).  We will find that there is a unique
involution in \(H_1(G)\) and will use that to detect differences
between homomorphisms.  Generators for \(\pi_1(G)\) and \(H_1(G)\)
are a set of \(k\) orientation reversing curves in \(G\), one that
crosses each crosscap.  The only relation in \(\pi_1(G)\) is that
the product of their squares is the identity.  In \(H_1(G)\), this
becomes the sole relation that the sum of their doubles is zero.  If
\(0\ne x\in H_1(G)\) has \(2x=0\), then \(2x\) is an even multiple
of the sum of the the generators, and \(x\ne0\) is an odd multiple
of the sum of the the generators.  Since twice the sum of the
generators is zero, \(x\) is reduced to the sum of the generators
and is seen to be unique.  We now build two more homomorphisms by
defining their values on the generators.  If there is more than one
generator (\(G\) is not \(P^2\)), then \(\phi_4\) takes one
generator to \(-1\) and all the rest to \(+1\).  If there are more
than 2 generators (\(G\) is neither \(P^2\) nor the Klein bottle),
then \(\phi_5\) takes exaclty two of the generators to \(-1\) and
all the rest to \(+1\).  Neither \(\phi_4\) nor \(\phi_5\) is equal
to the trivial \(\phi_1\).  Neither is equal to \(\phi_3\) since
each is \(+1\) on some orientation reversing curve.  They are not
equivalent since on \(H_1(G)\) they disagree on the unique
involution \(x\).

If \(G\) is a closed, connected surface, we have constructed one
homomorphism if \(G\) is \(S^2\), two if \(G\) is orientable and not
\(S^2\), two if \(G\) is \(P^2\), three if \(G\) is the Klein
bottle, and four if \(G\) is non-orientable and is neither \(P^2\)
nor the Klein bottle.

\begin{thm} If \(G\) is a closed, connected surface and
\(\phi:\pi_1(G)\into\ints_2\) is a homomorphism, then \(\phi\) is
equivalent to one of \(\phi_i\), \(i\in\{1,2,3,4,5\}\).  \end{thm}

\begin{demo}{Proof} If \(G\) is \(S^2\), \(P^2\), the Klein bottle,
or the connected sum of three copies of \(P^2\), then the statement
follows because there are not enough generators to create more
examples.  In the other cases, assume that \(\pi_1(G)\) has been
endowed with the standard generating sets described above.  If \(G\)
is orientable and not \(S^2\), then what must be shown is that if
\(\phi\) takes at least one generator to \(-1\), then it is
equivalent to \(\phi_2\).  Inductively we will show that if \(\phi\)
does not take all generators to \(-1\), then it is equivalent to a
homomorphism that takes more generators to \(-1\).  If \(G\) is
non-orientable, then what must be shown is that if \(\phi\) takes
more than two but not all generators to \(-1\), then it is
equivalent to a homomorphism that takes fewer generators to \(-1\).
These demonstrations will suffice because the homeomorphisms of
\(G\) are rich enough to take any one set of standard generating
curves for \(\pi_1(G)\) to any other set of standard generating
curves.

\FIGURE{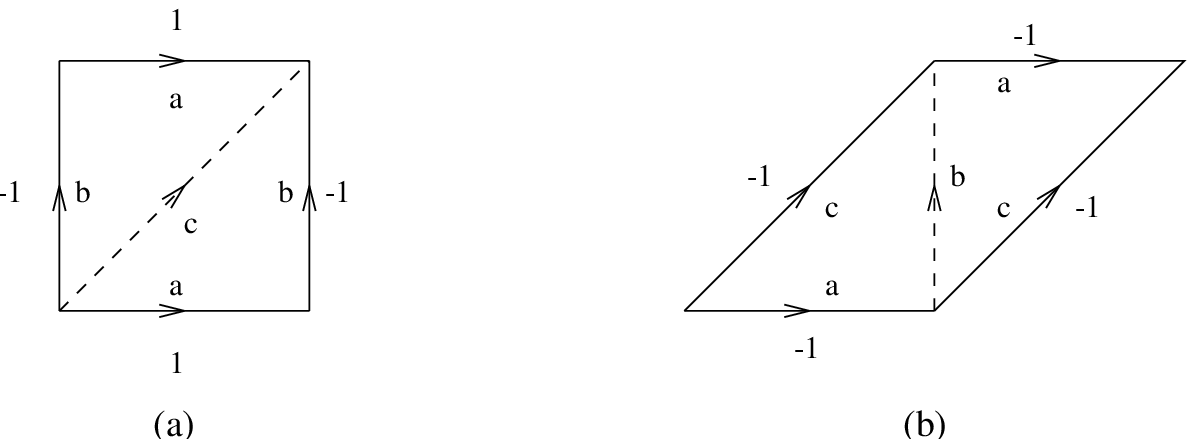}{One orientable handle}

If \(G\) is orientable and has a handle with one generating curve
taken to \(+1\) and one to \(-1\), then we can find two generating
curves that are taken to \(-1\).  In Figure \ref{1handle}(a) a curve
piercing the ``side'' labeled \(a\) is taken to \(+1\) and a curve
piercing the ``side'' lableled \(b\) is taken to \(-1\).  The handle
is homeomorphic to the diagram in Figure \ref{1handle}(b) where the
generating curves that pierce the two shown sides pass across the
line \(b\) and both take value \(-1\).

\FIGURE{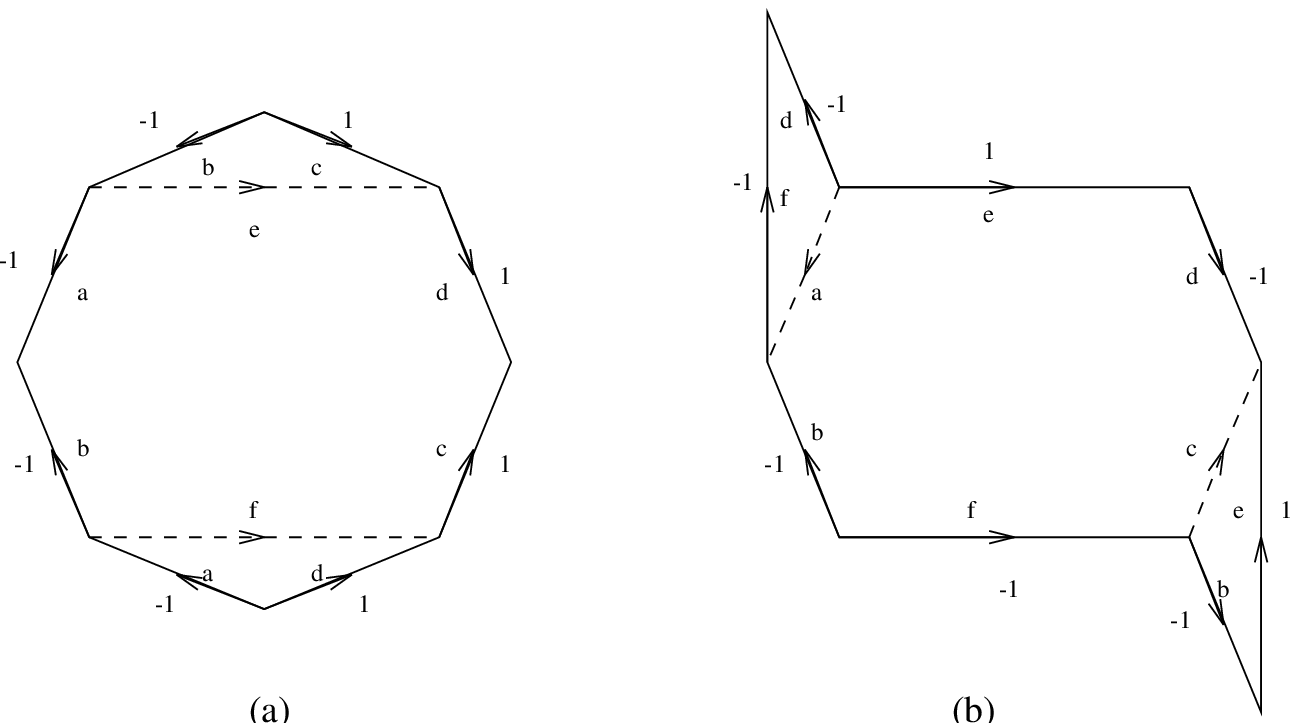}{Two orientable handles}

If \(G\) has two handles, one with generators taken to \(-1\) and
one with generators taken to \(+1\), then we wish to raise the
number of generators in the two handles that are taken to \(-1\) to
at least three.  In Figure \ref{2handles}(a), the ``sides''
\((a,b)\) form a handle whose generators are taken to \(-1\), and
the ``sides'' \((c,d)\) form a handle whose generators are taken to
\(+1\).  The surface shown is homeomorphic to the diagram in Figure
\ref{2handles}(b) where \((b,f)\) and \((d,e)\) form handles and the
generators piercing \(b\), \(d\) and \(f\) are taken to \(-1\),
while the generator piercing \(e\) is taken to \(+1\).

The two changes are enough to complete the induction.

\FIGURE{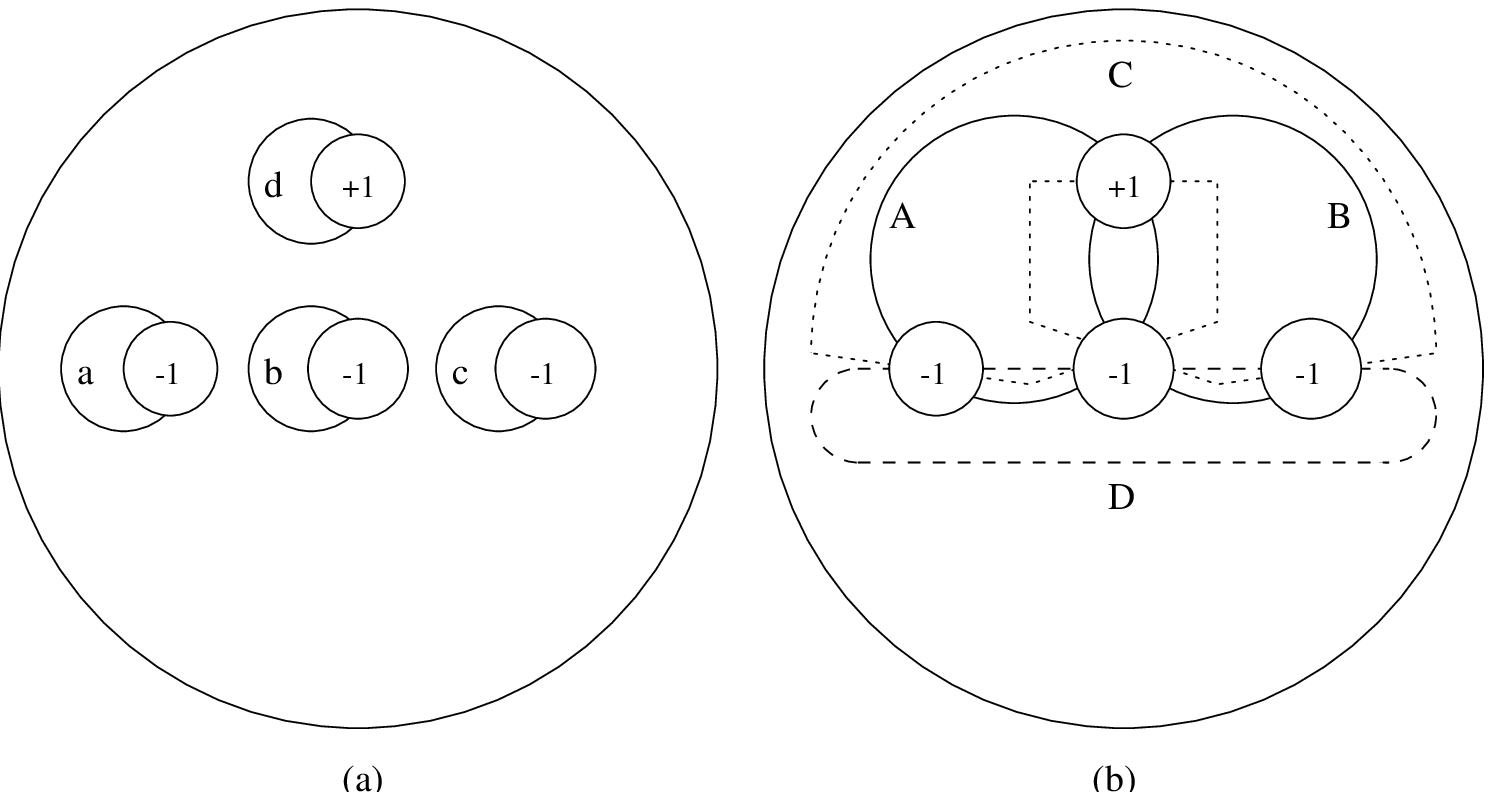}{Four crosscaps}

If \(G\) is non-orientable, has at least four crosscaps, and has
three of the corresponding generators taken to \(-1\) and one taken
to \(+1\), then we wish to find new generators that have fewer taken
to \(-1\).  In Figure \ref{4crosscaps}(a), the four small circles
represent four crosscaps, and the four labeled curves piercing them
are four disjoint, orientation reversing curves that are generators.
The values given the various generators are shown in the circles.
In Figure \ref{4crosscaps}(b), the four labeled curves are also four
disjoint, orientation reversing curves (they each pierce an odd
number of crosscaps), and there is a homeomorphism taking the
diagram of Figure \ref{4crosscaps}(a) to Figure \ref{4crosscaps}(b)
taking the curves \((a,b,c,d)\) to the curves \((A,B,C,D)\).  The
values given the new generators can be read from Figure
\ref{4crosscaps}(b) by counting the number of times they cross the
crosscaps with value \(-1\).  Thus \(A\), \(B\) and \(C\) are taken
to \(+1\) and \(D\) is taken to \(-1\).  That exactly two fewer
generators are taken to \(-1\) is not surprising since the value on
the unique involution in \(H_1(G)\) has to stay the same.
\end{demo}

We consider what this result tells us about classes of Seifert
fibered spaces.  We look first at the orientability of the spaces in
various classes.

Let \(M\) be a compact, connected Seifert fibered space, \(G\) its
orbit surface, \(p:M\into G\) the projection, and
\(\phi:\pi_1(G)\into\ints_2\) the classifying homomorphism.  Let
\(\rho:\pi_1(G)\into\ints_2\) be the homomorphism taking all
orientation reversing curves in \(G\) to \(-1\).  If \(G\) is
orientable, then \(\rho\) is trivial.  Let \(l\) be a loop in \(M\).
We wish to see if \(l\) reverses orientation.  Since we can homotop
\(l\) without changing its behavior on orientation, we ask that
\(l\) avoid exceptional fibers.  We can establish an ``orientation
frame'' along each point of \(l\) where two of the directions lie in
a meridinal disk of a fibered solid torus neighborhood and the third
lies in the direction of a fiber.  The trace of the two meridinal
directions can be followed in the orbit surface \(G\) as the loop
\(pl\) runs below \(l\).  We get an inversion of orientation of the
two meridinal directions if \(pl\) is an orientation reversing loop
in \(G\).  That is, if \(\rho(pl)=-1\).  We get an inversion in the
fiber direction if \(\phi(pl)=-1\).  The orientation of the
3-dimensional space \(M\) is reversed if one or the other but not
both of these orientations is reversed.  Thus \(l\) reverses
orientation in \(M\) if and only if \(\rho(pl)\phi(pl)=-1\).  This
says that \(M\) is orientable if and only if the pointwise product
\(\rho\phi\) takes on only the value \(+1\).  That is, if and only
if \(\phi=\rho\).  We can now look at cases.

If \(G\) is orientable, then \(\rho\) is trivial and \(M\) is
non-orientable if and only if \(\phi\) takes on negative values.
Thus if \(G\) is \(S^2\), there is only one class, and it has only
orientable spaces.  If \(G\) is orientable and not \(S^2\), then
there are two possible classes, one with all orientable spaces and
the other with all non-orientable spaces.  If \(G\) is
non-oreintable, then there are always at least two classes.  One,
corresponding to trivial \(\phi\), with non-orientable spaces, and
one, with \(\phi=\rho\) with orientable spaces.  If \(G\) has more
than one crosscaps, there is a third class, and if \(G\) has more
than two crosscaps, there is a fourth class.  The spaces in these
extra classes are all non-orientable since \(\phi\) and \(\rho\)
disagree on at least some curves.

Interesting things can happen in the non-orientable spaces.  Let
\(M\), \(G\), \(\phi\) and \(\rho\) be as above with \(M\)
non-orientable.  Let \(\alpha\) be a loop in \(G\).  Let \(\alpha\)
be based at an ordinary point \(x\) and let \(H\) be the fiber in
\(M\) over \(x\).  By Theorem \ref{IsotopyExist}, there is a fiber
preserving isotopy \(\Phi\) of \(M\) starting at the identity that
drags \(H\) along a path of fibers over \(\alpha\).  We can require
that the isotopy take a fiber solid torus neighborhood \(T\) of
\(H\) back to itself.  We know that \(\Phi_1:T\into T\) reverses
orientation if and only if \(\rho(\alpha)\phi(\alpha)=-1\).  Since
\(M\) is non-orientable, such an \(\alpha\) must exist.  We can now
prove the following.

\begin{lemma}\mylabel{NonOFlexible} Let \(M\) be a non-orientable,
Seifert fibered space, let \(C\) be a torus boundary component of
\(M\), let \(H\) be a fiber in \(C\), and let \(Q\) be a crossing
curve for \(C\).  Let \(n\) be in \(\ints\).  Then there is a fiber
preserving, self homeomorphism of \(M\) that is fixed on all
boundary components of \(M\) except \(C\) so that the image of \(Q\)
is homologous to one of \(\pm(Q+2nH)\).  The homeomorphism can be
chosen to be orientation preserving on \(C\), and it can be chosen
to be orientation reversing on \(C\).  It may not be possible to
chose which of \(Q+2nH\) or \(-Q+2nH\) is the image of \(Q\).
\end{lemma}

\begin{demo}{Proof} We can sew an ordinary solid torus \(T\) to
\(C\) and use the isotopy of the discussion above.  This will give
us a homeomorphism from \(M\) to \(M\) that has certain effects on
\(C\).  The ability to make things happen comes from the fact that
the way we sew \(T\) to \(C\) is up to us.  We will first get a
homeomorphism that reverses the orientation on \(C\), and then we
alter it so that it has the same effect on \(Q\) but preserves
orientation.

Let \(T\) be sewn to \(C\) so that the meridian of \(T\) goes to the
crossing curve \(Q+nH\).  Denote the resulting Seifert fibered space
by \(\hat{M}\).  We can drag the centerline of \(T\) around an
orientation reversing curve of \(\hat{M}\) so that \(T\) is returned
to itself with its orientation reversed.  This can be done by an
isotopy \(\Phi\) that moves no points on any of the boundary
components of \(\hat{M}\).  The map \(\Phi_1|M:M\into M\) is our
first candidiate.  Since the isotopy is fiber preserving and the
meridian of \(T\) must be returned to itself, the action of
\(\Phi_1\) on \(C=\bd T\) is mostly determined.  We use column
vectors to denote \(H_1(T)\) and we let \((H,Q)\) be a generating
pair for \(H_1(T)\).  In the simple case where \(n=0\), we are
returning \(Q\) to \(\pm Q\) and \(H\) to \(\pm H\).  Since
\(\Phi_1\) is orientation reversing on \(T\), the matrix involved
must be \(\pm\mtbt{1}{0}{0}{-1}\).

In the case where \(n\ne0\), then \(H\) and \(Q+nH\) are each fixed
up to sign, and we get the matrix by conjugating by the change of
basis matrix that takes \(H\) to \(H\) and \(Q\) to \(Q+nH\).  The
resulting matrix is
\[\pm\mtbt{1}{n}{0}{1}\mtbt{1}{0}{0}{-1}\mtbt{1}{-n}{0}{1}=
\pm\mtbt{1}{-2n}{0}{-1}.\] This gives the result with a map
\(\Phi_1|_M\) that reverses orientation on \(C\).  We get a map that
preserves the orientation on \(C\) by composing \(\Phi_1|_M\) with a
second orientation reversing map constructed by using the curve
\(Q+2nH\) as the image of the meridian of \(T\) so that \(Q+2nH\) is
preserved up to sign.  \end{demo}

\subsection{Closed Seifert fibered spaces without exceptional
fibers}

If \(M\) is a closed, connected Seifert fibered space with no
exceptional fibers, then a single fiber can be drilled out by
removing the interior of a fibered solid torus \(N\) and we are left
with a classifying space \(M_0\) for the class of \(M\).  That is,
\(M_0\) has no exceptional fibers and has a single torus boundary
component.  The classifying homomorphism for \(M_0\) is determined
by that for \(M\) and in turn this completely determines \(M_0\).
Thus \(M\) is determined by its classifying homomorphism and the way
that \(N\) was attached to \(M_0\).  The way that \(N\) is attached
is important because it turns out that different spaces can be
obtained from \(M_0\) by attaching \(N\) in different ways.  The
lens spaces \(L_{p,1}\) that can be realized with no exceptional
fibers are examples.

Let \(N\) be an ordinary solid torus, and let \(h\) be a fiber
preserving homeomorphism from \(\bd N\) to \(\bd M_0\).  Let \(M_h\)
be the Seifert fiber space obtained by sewing \(N\) to \(M_0\) using
\(h\).  Let \(m\) be a meridian for \(N\) and let \(H\) be a fiber
on \(\bd N\).  We can also use \(H\) as a longitude.  By Theorem
\ref{UniqueFill}, we know that the fiber type of \(M_h\) is
determined by the homology class (up to sign) of \(h(m)\).  Since
\(N\) is an ordinary solid torus, \(h(m)\) must be a crossing curve
on \(\bd M_0\).  Given one fixed crossing curve \(Q\) on \(\bd
M_0\), all other crossing curves \(Q'\) have the form \(\pm Q+nH\),
\(n\in\ints\).  Since the sign is irrelevant in determining \(M_h\),
we have a \(\ints\) indexed set of possibilities \(Q+nH\) to use for
the image of \(m\) under \(h\).  We will show that for orientable
\(M_0\), these each lead to different fiber types for \(M_h\), and
that for non-orientable \(M_0\) these collapse to two types.

We consider the structure of \(M_0\).  By Corollary
\ref{SectionExists}, we know that there is at least one curve \(Q\)
on \(\bd M_0\) that is a crossing curve and that is the boundary of
a section \(\sigma\) of the orbit surface.  We will show that if
\(M_0\) is orientable, then \(Q\) is essentially unique, and if
\(M_0\) is non-orientable, then \(Q\) is determined only up to an
even multiple of a fiber.

Let \(G_0\) be the orbit surface for \(M_0\) and let \(\sigma'\) be
another section for \(G_0\).  Let \(Q'\) be the curve bounding the
image of \(\sigma'\).  We can adjust \(\sigma\) and \(\sigma'\) by
sliding along fibers so that \(Q\) and \(Q'\) are transverse.  We
are interested in \(\hint(Q,Q')\).

Let \(\{\alpha_i\}\) be a set of arcs in \(G_0\) that cut \(G_0\)
into a disk \(E\).  The preimages of these arcs under the projection
are saturated annuli \(A_i\) in \(M_0\) that cut \(M_0\) into an
ordinary solid torus \(T_0\).  After the cut we get two copies
\(\alpha'_i\) and \(\alpha''_i\) for each arc and two copies
\(A'_i\) and \(A''_i\) for each annulus.  We can ``cut'' \(\sigma\)
and \(\sigma'\) by restriction into sections for \(E\) into \(T_0\).
Let \(\bd E\) be given an orientation and let \(R\) and \(R'\) be
the images under the ``cut'' \(\sigma\) and \(\sigma'\) of \(\bd
E\).  Let \(R\) and \(R'\) inherit their orientations from \(\bd
E\).  Further adjustment of \(\sigma\) can assure that \(R\) and
\(R'\) are transverse and that they do not intersect at the
boundaries of the \(A'_i\) and \(A''_i\).

The curves \(R\) and \(R'\) contain all of \(Q\) and \(Q'\) broken
into various arcs.  These arcs are separated by arcs in \(R\) and
\(R'\) that map into the \(A'_i\) and \(A''_i\).  Thus the
intersections of \(R\) and \(R'\) include all intersections of \(Q\)
and \(Q'\) plus extra intersections that occur inside the \(A'_i\)
and \(A''_i\).  Symbolically we can say \begin{equation}\hint(R,R')=
\hint(Q,Q')+\sum\limits_i\hint(R\cap (A'_i\cup A''_i),
R'\cap(A'_i\cup A''_i)).\mylabel{ScSum}\end{equation}

Each intersection of \(R\) and \(R'\) inside some \(A'_i\) has a
corresponding intersection in \(A''_i\).  Since \(M_0\) is
orientable, \(A'_i\) is sewn to \(A''_i\) by a homeomorphism that
reverses the orientations that \(A'_i\) and \(A''_i\) inherit from
\(T_0\).  We need to know how this homeomorphism deals with the
orientations of \(R\) and \(R'\) at these corresponding
intersections.  Since \(\sigma\) and \(\sigma'\) are sections, each
of \(R\) and \(R'\) intersects an \(A'_i\) in a single arc that is
transverse to the fibers of \(A'_i\).  These arcs are the images of
\(\alpha'_i\).  These arcs cross \(A'_i\) in the same direction in
that they originate at the same boundary component of \(A'_i\).  The
same statements apply to the intersection of \(R\) and \(R'\) with
\(A''_i\).  Thus the sewing of \(A'_i\) to \(A''_i\) sews \(R\cap
A'_i\) to \(R\cap A''_i\) and \(R'\cap A'_i\) to \(R'\cap A''_i\) in
a way that preserves the orientations of both arcs or reverses the
orientations of both arcs.  Thus intersections of \(R\) and \(R'\)
inside some \(A'_i\) are mapped to corresponding intersections in
\(A''_i\) by a map that behaves the same on the orientations of both
\(R\) and \(R'\) and that reverses the orientations of \(A'_i\) and
\(A''_i\).  Thus the corresponding intersections have opposite sign
and cancel in the sum \ref{ScSum}.  This gives
\(\hint(R,R')=\hint(Q,Q')\).  But \(R\) and \(R'\) are meridians on
the solid torus \(T_0\) and have intersection zero.  Thus for
\(M_0\) orientable, boundaries of sections are homologous up to sign
on \(\bd M_0\).

If \(M_0\) is not orientable, then we do not know how \(A'_i\) is
sewn to \(A''_i\).  Thus intersections of \(R\) and \(R'\) in
\(A'_i\) may cancel with intersections in \(A''_i\) or double with
them.  The only conclusion that we can reach is that \(\hint(Q,Q')\)
is an even number.  Thus for \(M_0\) non-orientable, if \(Q\) is a
boundary for a section, then any other boundary of a section must be
homologous up to sign to \(Q+2nH\), \(n\in\ints\).  However the
converse is also true.  Each homology class of the form \(Q+2nH\) is
represented by the boundary of a section because for \(M_0\)
non-orientable, Lemma \ref{NonOFlexible} gives a fiber preserving,
self homeomorphism of \(M_0\) that takes \(Q\) to the class
\(Q+2nH\).

The information that we need to analyze sewings is summarized in the
following.

\begin{lemma} Let \(M_0\) be a compact, connected Seifert fibered
space with no exceptional fibers and whose boundary consists of a
single torus \(C\).  If \(M_0\) is orientable, then there is only
one simple closed crossing curve on \(C\) up to isotopy and reversal
that is the boundary of a section for the orbit surface of \(M_0\).
If \(M_0\) is non-orientable, then there are infinitely many simple
closed crossing curves on \(C\) that are boundaries of sections of
the orbit surface, and given any one such curve \(J_1\), then any
crossing curve \(J_2\) on \(C\) is the boundary of a section if and
only if it differs from \(J_1\) by an even multiple of the fiber in
which case there is a fiber preserving self homeomorphism of \(M_0\)
that carries \(J_1\) to \(J_2\) or its reverse.  \end{lemma}

\begin{demo}{Proof} The last provision can be deduced from Lemma
\ref{StraightLevels}.  \end{demo}

We return to \(M\), a closed, connected Seifert fibered space with
no exceptional fibers.  What we have to say about \(M\) will depend
on the orientation of \(M\) if \(M\) is orientable.  Thus we assume
that \(M\) is oriented or that \(M\) is non-orientable.  Let \(N\)
be a fibered solid torus neighborhood of some fiber and let
\(M_0=M-\inter N\).  Let \(Q\) be a boundary of a section of the
orbit surface of \(M_0\) and let \(m\) be a meridian of \(N\).  Both
\(Q\) and \(m\) lie on \(\bd M_0=\bd N\) and we can look at
\(\hint(m,Q)\).

First assume that \(M\) is oriented.  Let \(N\) inherit the
orientation, and let \(\bd N\) be oriented consistently with \(N\).
Both \(Q\) and \(m\) are crossing curves (since \(N\) is an ordinary
solid torus), so for a fiber \(H\) on \(\bd N\), both \((H,Q)\) and
\((H,m)\) are valid generating pairs for \(H_1(\bd N)\).  Pick an
orientation for \(H\).  The pairing \(\hint(\,\,,\,)\) is determined
by the orientation of \(\bd N\).  There are unique orientations on
\(Q\) and \(m\) so that \(\hint(H,Q)=\hint(H,m)=+1\).  Let
\(b=\hint(m,Q)\).  Using the generating pair \((H,Q)\) for \(H_1(\bd
N)\), this expresses \(m\) as
\begin{equation}m=\hint(m,Q)H+\hint(H,m)Q=bH+Q.
\mylabel{ObstructionMeridian}\end{equation}

If \(M\) is non-orientable, then chose an orientation for \(N\)
arbitrarily and use the rest of the steps above to calculate
\(b=\hint(m,Q)\mod2\) given as either \(+1\) or \(0\).

Various choices were made in the calculation of \(b\).  If the
opposite choice of orientation is made for \(H\), then the
requirements on \(Q\) and \(m\) will force simultaneous reversal of
orientations on them and we will preserve \(\hint(m,Q)\).  If a
different fiber \(H'\) and fibered solid torus neighborhood \(N'\)
are chosen for the drilling, then Theorem \ref{IsotopyExist} gives
an (orientation preserving), fiber preserving self homeomorphism of
\(M\) carrying \(N\) to \(N'\).  This will carry a section of
\(M-\inter N\) to a section of \(M-\inter N'\) and a meridian of
\(N\) to meridian of \(N'\) (in an orientation preserving way if an
orientation is present).  This will also preserve \(\hint(m,Q)\).
If \(M\) is not orientable, then the opposite choice of orientation
for \(N\) simply reverses the sign of \(b\) which makes no
difference mod 2.  From this it is clear that we get the same
calculation for \(b\) from any Seifert fibered space that has the
same (oriented) fiber type as \(M\).

Consider the triple \((G,\phi,b)\) associated with \(M\) where \(G\)
is the orbit surface, \(\phi:\pi_1(G)\into\ints_2\) is the
classifying homomorphism and \(b\) is as calculated above.  This
triple is a well defined invariant of the (oriented if applicable)
fiber type of \(M\).  We now wish to know whether the triple
distinguishes (oriented if applicable) fiber types.

Let \(M\) and \(M'\) be closed, connected Seifert fibered spaces
without exceptional points that are oriented or non-orientable and
let them have equivalent triples \((G,\phi,b)\) and
\((G',\phi',b')\) in that \(\phi\) and \(\phi'\) are equivalent and
\(b=b'\).  Since the orientability of a Seifert fibered space is
determined by its classifying homomorphism, we know that \(M\) and
\(M'\) are either both oriented or both non-orientable.  Drill out a
fiber from each of \(M\) and \(M'\) to create \(M_0\) and \(M'_0\)
respectively.  Since the classifying homomorphisms are equivalent,
there is a fiber preserving homeomorphism \(h\) from \(M_0\) to
\(M'_0\).

First assume that the spaces are orientable.  We can arrange that
\(h\) be orientation preserving since by Corollary
\ref{ReverseFibers}, there is an orientation reversing, fiber
preserving homeomorphism from \(M_0\) to itself.  Choose an
orientation for a fiber \(H\) in \(\bd M_0\) and orient \(h(H)\)
consistently.  Orient \(\bd M_0\) and \(\bd M'_0\) consistently with
the orientations of the ordinary solid tori that were removed from
\(M\) and \(M'\) to form \(M_0\) and \(M'_0\).  This forces
orientations on \(Q\), \(m\), \(Q'\) and \(m'\) with notation
continuing the pattern above so that \(h\) carries \(Q\) to a curve
homologous to \(Q'\) and \(m\) to a curve homologous to \(m'\).  Now
\(M\) is obtained from \(M_0\) by sewing an ordinary solid torus to
\(\bd M_0\) with meridian going to \(Q+bH\) and \(M'\) is obtained
from \(M'_0\) by sewing an ordinary solid torus to \(\bd M'_0\) with
meridian going to \(Q'+bH'\).  But \(h\) carries \(Q+bH\) to a curve
homologous to \(Q'+bH'\) and from Theorem \ref{UniqueFill} we know
that there is a homeomorphism from \(M\) to \(M'\) that extends
\(h\).

Assume that the spaces are non-orientable.  Now the curves \(Q\) and
\(Q'\) used in computing \(b\) and \(b'\) may not be related by
having \(h(Q)\) homologous to \(Q'\).  However, they are known to
differ by an even multiple of the fiber and Lemma \ref{NonOFlexible}
gives a fiber preserving, self homeomorphism of \(M_0\) which can be
used to alter \(h\) so that \(h(Q)\) is homologous to \(Q'\).  We
finish the argument as in the orientable case.  We have shown.

\begin{thm} Let \(M\) be a closed, connected Seifert fibered space
with no exceptional fibers, either oriented or non-orientable.  Then
the triple \((G,\phi,b)\) is a well defined invariant of the
(oriented if applicable) fiber type of \(M\) that completely
determines the (oriented if applicable) fiber type of \(M\).
\BlackBox\end{thm}

We refer to the integer \(b\) as the obstruction to a section for
the orbit surface.

We can now enumerate the closed, connected Seifert fibered spaces
with no exceptional fibers.  There are four possible combinations of
the orientability of the Seifert fibered space and its orbit
surface.  We indicate these by \((O,o)\), \((O,n)\), \((N,o)\) and
\((N,n)\).  The first (upper case) letter refers to the
orientability \((O)\) or non-orientability \((N)\) of the Seifert
fibered space.  The second (lower case) letter refers to orbit
surface.  For a given surface \(G\) there is only one class of
Seifert fibered spaces over \(G\) that gives the combinations
\((O,o)\), \((O,n)\) and \((N,o)\).  This is because there are only
two classes of classifying homomorphisms for an orientable orbit
surface, one resulting in orientable Seifert fibered spaces, and the
other in non-orientable spaces, and there is only one class of
homomorphisms, the one containing the orientation homomorphism
\(\rho\), that give orientable Seifert fibered spaces over a
non-orientable orbit surface.  For a given non-orientable surface,
there are up to three classes that give \((N,n)\).  We distinguish
these by labeling them \((N,n,\hbox{I})\), \(N,n,\hbox{II})\) and
\((N,n,\hbox{III})\).  The first uses the class of the trivial
homomorphism from \(\pi_1(G)\) to \(\ints_2\), the second and third
use the non-trivial homomorphisms that disagree with \(\rho\) where
the second takes the unique involution of \(H_1(G)\) to \(-1\), and
the third takes takes it to \(+1\).  Since an orientable surface is
determined by its genus \(g\) (number of handles) and a
non-orientable surface is determined by the number \(k\) of
crosscaps, we have a complete list of classes in the symbols
\((O,o,g)\), \((O,n,k)\), \((N,o,g)\), \((N,n,\hbox{I},k)\),
\((N,n,\hbox{II},k)\) and \((N,n,\hbox{III},k)\).  There are some
restrictions on \(g\) and \(k\).  For example there is no class
\((N,n,\hbox{III},2)\) since the Klein bottle does not have enough
crosscaps to admit the required homomorphism.  Within each class
there are certain closed, connected Seifert fibered spaces with no
exceptional fibers.  These are identified up to (oriented if
applicable) fiber type by \((O,o,g\mid b)\), \((O,n,k\mid b)\),
\((N,o,g\mid b)\), \((N,n,\hbox{I},k\mid b)\),
\((N,n,\hbox{II},k\mid b)\) and \((N,n,\hbox{III},k\mid b)\) where
\(b\in\ints\) if the Seifert fibered space is oriented and
\(b\in\{0,1\}\) if it is non-orientable.  Note that for a given
closed, connected orbit surface, there are only finitely many (up to
a maximum of 6) fiber types of non-orientable Seifert fibered spaces
with no exceptional fibers.

\subsection{Crossing curves of fibered solid tori}

We will analyze Seifert fibered spaces with exceptional fibers by
relating them to spaces in the same class with no exceptional
fibers.  From a space \(M\) with an exceptional fiber, we can create
a space \(M'\) with one fewer exceptional fiber by drilling out the
exceptional fiber and replacing it with an ordinary fiber.  That is,
a fibered solid torus neighborhood \(N\) of the exceptional fiber is
replaced by an ordinary solid torus \(T\).  What \(M\) and \(M'\)
have in common is the space \(M_0=M-\inter N=M'-\inter N'\).  To go
in a well defined manner from \(M\) to \(M'\), we need to be able to
tell where the meridian of \(N'\) should go on \(\bd M_0\).  This is
a crossing curve since \(N'\) is an ordinary solid torus.  To go in
a well defined manner from \(M'\) to \(M\), we need to be able to
tell where the meridian of \(N\) should go on \(\bd M_0\).  Thus we
need a way of starting with a fibered solid torus \(N\) and picking
out a well defined crossing curve and a technique for recovering the
original meridian after the fibered solid torus has been forgotten.
This is the subject of the current section.  We will need to apply
what we say to both orientable and non-orientable manifolds, so we
consider both settings here.

Let \(N\) be a fibered solid torus.  We assume first that \(N\) is
oriented.  Let \(m\) be a meridian arbitrarily oriented.  Let \(H\)
be a fiber in \(\bd N\), oriented so that \(\mu=\hint(H,m)\) is
positive.  (Recall that \(\hint(H,m)\) can never be zero.)  In what
follows we do not need to refer to a longitude, so we avoid it.  A
crossing curve \(Q\) for \(\bd N\) will intersect \(H\) once, so
\(H\) and \(Q\) can be used to generate \(H_1(\bd N)\).  Since \(N\)
is oriented, we get an orientation on \(\bd N\).  We can use this
orientation to declare that certain orientations on \(H\) and \(Q\)
are consisitent with the orientation.  This condition can be
captured by the equality \(\hint(H,Q)=+1\).  If one crossing curve
\(Q'\) is found so that \(\hint(H,Q')=+1\), then all other crossing
curves \(Q''\) with \(\hint(H,Q'')=+1\) take the form \(Q''=Q'+nH\),
\(n\in\ints\).  Thus all such crossing curves \(Q''\) have
\(\hint(m,Q'')\) defined up to integral multiples of \(\mu\).  There
is then a unique crossing curve \(Q\) and a unique integer \(\beta\)
with \(\hint(H,Q)=+1\) and \(\beta=\hint(m,Q)\) with
\(0\le\beta<\mu\).  If \(m\) were chosen with the opposite
orientation, then simultaneously reversing the orientation of \(H\)
and \(Q\) would then satisfy all the conditions above with the same
integer \(\beta\).  The only choice made in the above discussion was
the orientation of \(m\), so \(N\) and its orientation determine
\(Q\) up to reversal and determine \(\beta\) completely.

Now assume that the orientation on \(\bd N\) is known, that \(H\)
and \(Q\) are known without fixed orientations, and that \(\mu\) and
\(\beta\) are known.  If we demand that \(\hint(H,Q)=+1\), then
there are two possible choices for the orientations of \(H\) and
\(Q\) as a pair, and these differ from each other by simultaneous
reversal.  These two choices are identical to the two orientation
arrangements discussed in the previous paragraph.  Thus setting
\(m=\beta H+\mu Q\) recovers \(m\) up to reversal of orientation.

We thus have for oriented \(N\) that \(H\) and \(Q\) are determined
up to orientation and the pair \((\mu,\beta)\) is uniquely
determined.  Further, the orientation on \(\bd N\), the unoriented
curves \(H\) and \(Q\) and the pair \((\mu,\beta)\) determine \(m\)
up to orientation.

What has actually taken place in the previous paragraphs is that
\(H\) and \(Q\) can be expressed in terms of a meridian-longitude
generating pair \((m,l)\) for \(H_1(\bd N)\) as the columns of
\(\mtbt{\nu}{\alpha}{-\mu}{\beta}\), and \(m\) and \(l\) can be
expressed in terms of the generating pair \((H,Q)\) as the columns
of the inverse matrix \(\mtbt{\beta}{-\alpha}{\mu}{\nu}\) which
expresses \(m\) as \(m=\beta H+\mu Q\).

If \(N\) is not oriented, then \(H\) and the set of possible
crossing curves \(Q''\) are still the same, but there is no way to
chose among the possible orientations.  The set of intersection
numbers \(\hint(m,Q'')\) are ambiguous as to sign.  Thus \(m\),
\(H\) and the crossing curves are oriented arbitrarily and we define
\(\mu\) as \(|\hint(H,m)|\) and \(\beta\) as the minimum of
\(|\hint(m,Q'')|\) over all crossing curves \(Q''\).  This gives
\(0\le\beta\le\mu/2\).  For \(\mu>2\), there will only be one
crossing curve \(Q\), up to change in orientation, that achieves ths
minimum.  For \(\mu=2\), \(Q\) is not determined.  In terms of a
generating meridian-longitude pair, we can have \(H=\cv{1}{-2}\)
with both \(Q=\cv{0}{1}\) and \(Q=\cv{1}{-1}\) giving
\(\hint(H,Q)=1\) and \(|\hint(m,Q)|=1\).  This ambiguity is
unavoidable and will be shown later not to be important.

Now assume that \(\bd N\) has no fixed orientation, that \(H\) and
\(Q\) are known, that the pair \((\mu,\beta)\) is known, and that
\(\mu>2\).  We can express \(m\) as \(m=\beta H+\mu Q\).  There are
four combinations of orientations of \(H\) and \(Q\) that break into
two groups of two combinations.  In a given group, one combination
is obtained from the other by simultaneous reversal of \(H\) and
\(Q\).  Within each group \(m\) is determined up to reversal.
However, the two groups give different curves for \(m\).  Up to
sign, the two possibilities for \(m\) are represented by \(\pm\beta
H+\mu Q\).

Thus for unioriented \(N\) with centerline of index at least 3, we
get uniquely determined \(Q\), up to reversal, and uniquely
determined \((\mu,\beta)\).  These determine two possible curves, up
to orientation, for \(m\).  We will also see later that this
ambiguity is not important.

For a fibered solid torus \(N\), oriented or not, we call
\((\mu,\beta)\) the \defit{crossing invariants determined by} \(N\).
We call \(Q\) the \defit{crossing curve determined by} \(N\) except
when \(N\) is unoriented and its centerline is a fiber of index 2.

\subsection{Closed, oriented Seifert fibered spaces}

If \(M\) is a closed, connected, oriented Seifert fibered space,
then there is a finite set \(\{H_1,H_2,\dots,H_n\}\) of exceptional
fibers.  Let \(\{N_1,N_2,\dots,N_n\}\) be pairwise disjoint fibered
solid torus neighborhoods of the exceptional fibers.  Let
\(\{Q_1,Q_2,\dots,Q_n\}\) be the crossing curves determined by the
neighborhoods, and let \((\mu_1,\beta_1,\dots,\mu_n,\beta_n)\) be
the crossing invariants.  We can assume that the \(H_i\) have been
numbered so that the \(\mu_i\) are in non-decreasing order, and the
\(\beta_i\) in a group of constant \(\mu_i\) are in non-decreasing
order.  The \(N_i\) can be removed and replaced by ordinary solid
tori \(\{N'_1,N'_2,\dots,N'_n\}\) where the sewing is determined by
requiring that a meridian of each \(N'_i\) be sewn to \(Q_i\).  This
gives a space \(M'\) with no exceptional fibers.  Note that the
invariance of fibered solid torus neighborhoods of fibers and the
invariance of the crossing curves shows that \(M'\) is completely
determined up to oriented fiber type by \(M\).  The space \(M'\) is
determined up to oriented fiber type by one of the symbols
\((O,o,g\mid b)\) or \((O,n,k\mid b)\).  We can thus associate one
of the following symbols with \(M\): \[(O,o,g\mid
b,\mu_1,\beta_1,\dots,\mu_n,\beta_n)\] or \[(O,n,k\mid
b,\mu_1,\beta_1,\dots,\mu_n,\beta_n).\] The symbol is completely
determined by \(M\).

If one of the symbols above is given, then we build a realization.
We start with an oriented space \(M'\) with no exceptional fibers
determined by \((O,o,g\mid b)\) or \((O,n,k\mid b)\).  We chose
\(n\) fibers and pairwise disjoint ordinary solid torus
neighborhoods \(N'_i\) of these fibers.  Each \(N'_i\) has a
meridian \(Q_i\).  The orientation that \(N'_i\) inherits from
\(M'\) determines an orientation on \(\bd N'_i\).  The orientation
on \(\bd N'_i\), a fiber \(H_i\) in \(\bd N'_i\), and the meridian
\(Q_i\) determine a curve \(M_i=\beta_i H_i+\mu_i Q_i\) up to
orientation.  There is a fibered solid torus \(N_i\) that can be
sewn to \(\bd N'_i\) in a way that preserves fibers and takes a
meridian to \(m_i\).  The resulting space is a realization of the
symbol.  We argue that any two realizations have the same oriented
fiber type.  The oriented fiber type of \(M'\) is determined.  If
different fibers and neighborhoods are chosen for the \(N'_i\) then
the transitivity of fibers will give an orientation preserving,
fiber preserving homeomorphism that carries one set of neighborhoods
to the other.  It will also carry one set of meridians to the other.
The curves for attaching the new meridians are determined up to
reversal and Theorem \ref{UniqueFill} guarantees that the resulting
space is determined up to oriented fiber type.  If the given symbol
came from an oriented space \(M\), then \(M\) is a realization
because of the way that the symbol was derived from \(M\).  Thus all
realizations are of the same oriented fiber type of \(M\).  We have
the following.

\begin{thm} Let \(M\) be a closed, connected, oriented Seifert fiber
space.  Then the appropriate symbol \[(O,o,g\mid
b,\mu_1,\beta_1,\dots,\mu_n,\beta_n)\] or \[(O,n,k\mid
b,\mu_1,\beta_1,\dots,\mu_n,\beta_n)\] is a well defined invariant
of the oriented fiber type of \(M\) that completely determines the
oriented fiber type of \(M\).  \BlackBox\end{thm}

\subsection{Closed, non-orientable Seifert fibered spaces}

Let \(M\) be a closed, connected, non-orientable Siefer fibered
space.  As in the previous section there is a finite set
\(\{H_1,\dots,H_n\}\) of exceptional fibers with corresponding
fibered torus neighborhoods \(\{N_1,\dots,N_n\}\).  The possibility
that some of the fibers are of index 2 changes the analysis, so we
start by assuming that there are none.  Thus we get a corresponding
set \(\{Q_1,\dots,Q_n\}\) of crossing curves determined by the
\(N_i\) and crossing invariants
\((\mu_1,\beta_1,\dots,\mu_n,\beta_n)\) ordered as in the previous
section.  As before, we can remove each \(N_i\) and replace it with
an ordinary solid torus \(N'_i\) with meridian sewn to \(Q_i\).
This creates a space \(M'\) with no exceptional fibers.  This space
is determined by one of the symbols \((N,o,k\mid b)\),
\((N,n,\hbox{I},k\mid b)\), \((N,n,\hbox{II},k\mid b)\), or
\((N,n,\hbox{II},k\mid b)\).  We then associate one of \[(N,o,k\mid
b,\mu_1,\beta_1,\dots,\mu_n,\beta_n),\] \[(N,n,\hbox{I},k\mid
b,\mu_1,\beta_1,\dots,\mu_n,\beta_n),\] \[(N,n,\hbox{II},k\mid
b,\mu_1,\beta_1,\dots,\mu_n,\beta_n),\] or \[(N,n,\hbox{III},k\mid
b,\mu_1,\beta_1,\dots,\mu_n,\beta_n)\] with \(M\).  As before, the
symbol is determined by \(M\).

Given a symbol as above, we build a realization in a well defined
manner as before up to the point of adding a new solid torus where
an ordinary solid torus \(N'_i\) was before.  We are faced with two
possible curves \(m=\pm\beta_i H_i+\mu_i Q_i\) to chose from.
However, it follows from Lemma \ref{NonOFlexible} that there is a
fiber preserving, orientation reversing map on \(M'-N'_i\) that
takes \(Q_i\) to \(\pm Q_i\) and \(H_i\) to \(\mp H_i\).  This is
seen to take one possibility for \(m\) to the other or its reverse.
Thus the two possible ways of attaching \(N_i\) are seen to yield
spaces of the same fiber type.

If there are fibers of index 2 present, then our ordering of the
crossing invariants has \(\mu_i=2\) for all \(i\) less than \(s\)
where \(s\) is the number of exceptional fibers of index 2.  All
\(\alpha_i\), \(i\le s\), equal 1.  The crossing curves \(Q_i\) for
\(i>s\) are determined.  There is a problem in chosing crossing
curves for \(i\le s\).  This will create no problem in recovering
\(M\) as is shown by the next lemma.

\begin{lemma} Let \(M_0\) be a non-orientable Seifert fibered space
with at least one torus boundary component \(C\).  Let \(N\) be a
fibered solid torus of type \(1/2\) and assume that it is possible
to obtain \(M\) from \(M_0\) by attaching \(N\) to \(M_0\) with a
fiber preserving homeomorphism \(h\) from \(\bd N\) to \(C\).  Then
the fiber type of \(M\) is independent of which fiber preserving
homeomorphism \(h\) is used.  \end{lemma}

\begin{demo}{Proof} To give a frame of reference for the duration of
the proof, we orient \(C\) and \(\bd N\) arbitrarily.  We know there
is a crossing curve \(Q\) on \(\bd N\) for which \(m=2Q+H\).  Since
the sewings of \(\bd N\) to \(C\) must be fiber preserving, the
sewings are distinguished by what the sewing maps do on \(Q\) and
what the maps do to the orientation of the fibers.  The possible
images of \(Q\) are all crossing curves on \(C\), and these are all
of the form \(\pm Q'+nH\), \(n\in\ints\), where \(Q'\) is one fixed
crossing curve on \(C\).  We have already seen that reversing the
effect of the orientations on either \(Q\) or \(H\) has no effect on
the result of the sewing.  Thus we can consider only curves of the
form \(Q'+nH\), \(n\in\ints\), and assume that \(h\) preserves the
orientations of the fiber.  We also know from Lemma
\ref{NonOFlexible} that curves that differ by even multiples of the
fiber give equivalent results.  We can thus consider only \(Q'\) and
\(Q'-H\) as the possible images of \(Q\).  The images of \(m\) under
these two sewings are \(2Q'+H\) and \(2Q'-H\).  But again Lemma
\ref{NonOFlexible} gives a fiber preserving, self homeomorphism of
\(M_0\) that fixes \(Q'\) and \(H\) up to sign and reverses the
orientation on \(C\).  This carries one image of \(m\) to the other
or its reverse and the resulting spaces are fiber homeomorphic.
\end{demo}

There is thus a lack of ambiguity in determining \(M\) from another
space.  However, fibers of index 2 do not determine crossing curves
and we have an ambiguity in determining another space from \(M\).
To get around this we will not drill out the fibers of index 2.
This leads to the following general procedure.

Let \(M\) be a closed, non-orientable Seifert fibered space.  Let
\(s\) be the number of excpetional fibers of index 2, let
\(\{H_1,\dots,H_n\}\) be the exceptional fibers of index greater
than 2, let \(\{N_1,\dots,N_n\}\) be pairwise disjoint fibered solid
torus neighborhoods of the \(H_i\), let \(\{Q_1,\dots,Q_n\}\) be the
crossing curves determined by the \(H_i\), and let
\((\mu_1,\beta_1,\dots,\mu_n,\beta_n)\) be the crossing invariants
of the \(H_i\) ordered in the usual way.  We create a space by
removing all the \(N_i\) and sewing in ordinary fibered solid tori
\(N'_i\) with meridians going to the curves \(Q_i\).  This creates a
space with \(s\) exceptional fibers of index 2.  We denote this
space by \(M_s\).  We can drill out the fibers of index 2 from
\(M_s\) to create a space \(M_0\).  Since \(M_0\) has boundary, it
is determined by its class and thus by \(M\).  Since it is
non-orientable, it is determined by one of the symbols \((N,o,g\mid
-)\), \((N,n,\hbox{I},k\mid -)\), \((N,n,\hbox{II},k\mid -)\) or
\((N,n,\hbox{III},k\mid -)\).  Since sewing in the \(s\) exceptional
fibers of index 2 gives a space that is independent of the sewings,
we have that \(M_s\) is dependent only on \(M\) and not on any
choices.  We now associate one of the following with \(M\):
\[(N,o,k\mid(b,s),\mu_1,\beta_1,\dots,\mu_n,\beta_n),\]
\[(N,n,\hbox{I},k\mid(b,s),\mu_1,\beta_1,\dots,\mu_n,\beta_n),\]
\[(N,n,\hbox{II},k\mid(b,s),\mu_1,\beta_1,\dots,\mu_n,\beta_n),\] or
\[(N,n,\hbox{III},k\mid(b,s),\mu_1,\beta_1,\dots,\mu_n,\beta_n)\]
where \(s\) is the number of exceptional fibers of index 2 and \(b\)
is irrelevant if \(s>0\).  The determination of the symbol by \(M\)
has been discussed above.  The determination of \(M\) by the symbol
follows from the fact that the symbol determines uniquely the space
\(M_s\) and the numbers \((\mu_1,\beta_1,\dots,\mu_n,\beta_n)\)
determine \(M\) from \(M_s\).  We have the result summarized as the
following.

\begin{thm} Let \(M\) be a closed, connected, non-orientable Seifert
fiber space.  Then the appropriate symbol
\[(N,o,k\mid(b,s),\mu_1,\beta_1,\dots,\mu_n,\beta_n),\]
\[(N,n,\hbox{I},k\mid(b,s),\mu_1,\beta_1,\dots,\mu_n,\beta_n),\]
\[(N,n,\hbox{II},k\mid(b,s),\mu_1,\beta_1,\dots,\mu_n,\beta_n),\] or
\[(N,n,\hbox{III},k\mid(b,s),\mu_1,\beta_1,\dots,\mu_n,\beta_n)\] is
a well defined invariant of the fiber type of \(M\) that completely
determines the fiber type of \(M\).  \BlackBox\end{thm}

\subsection{Compact Seifert fibered spaces with boundary}

Compact Seifert fibered spaces with boundary have been discussed to
some extent above.  We do not include a complete discussion of
exceptional fibers here, but we inlcude one lemma that shows that
even in the orientable case, there needs to be less information kept
about the location of crossing curves.  The use of the following is
when \(C\) is to be used as a place to attach a fibered solid torus.

\begin{lemma} Let \(M\) be a Seifert fibered space with at least two
boundary components \(C\) and \(C'\) with \(C\) a torus.  Let \(Q\)
and \(Q'\) be two crossing curves for \(C\).  Then there is a fiber
preserving, self homeomorphism of \(M\) that is fixed on all
boundary components of \(M\) except \(C\) and \(C'\) that carries
\(Q\) to \(Q'\) or its reverse.  \end{lemma}

\begin{demo}{Proof} Let \(G\) be the orbit surface and let \(J\) and
\(J'\) be the images of \(C\) and \(C'\).  Let \(\alpha\) be an arc
from \(J\) to \(J'\) that avoids exceptional points.  We can split
\(G\) along \(\alpha\) and we can split \(M\) to get \(M_-\) along
the annulus \(A\) that is the preimage of \(\alpha\).  We get two
copies \(A'\) and \(A''\) of \(A\) in the splitting.  There is a
fiber preserving isotopy of \(M_-\) (in the strong sense in that
each fiber is mapped to itself throughout the isotopy) that rotates
the fibers of \(A'\) through one full circle and that only moves
points in a small neighborhood of \(A'\).  The end of this isotopy
gives a fiber preserving self homeomorphism of \(M\) that carries a
given crossing curve to itself plus one copy of the fiber.
Repetitions of this suffice to carry any crossing curve to any other
or its reverse.  \end{demo}

\subsection{The diagram of a closed Seifert fibered space}

The ``symbol'' associated with a closed, connected Seifert fibered
space has the shape \((\,\,\mid\,\,)\) where the information to the
left of the vertical bar determines the class of the space, and the
information to the right of the bar determines the fiber type of the
space within that class.  The determination of the class of the
space can also be made with a geomeric object --- the classifying
space of the class --- the unique space in the class with no
exceptional fibers and exactly one torus boundary component.  We can
also associate a geometric object with the information to the right
of the vertical bar.  We will call this object a \defit{diagram}
associated with a closed, connected Seifert fibered space.

Let \(M\) be a closed, connected Seifert fibered space that is
either oriented or non-orientable.  Let \(\{H_1,\dots,H_n\}\) be the
set of exceptional fibers and let \(\{N_1,\dots,N_n\}\) be pairwise
disjoint fibered solid torus neighborhoods of the \(H_i\).  Let
\((\mu_1,\beta_1,\dots,\mu_n,\beta_n)\) be the crossing invariants
ordered in the usual way.  Let \(s\) be the number of \(\mu_i\)
equal to 2.  If \(s\ne0\) and \(M\) is not orientable, then crossing
curves for the corresponding \(\bd N_i\) are not determined.  Let
\(\{Q_1,\dots,Q_n\}\) be crossing curves determined by the \(N_i\)
where possible and chosen arbitrarily where not.  Let \(M'\) be
obtained from \(M\) by replacing the \(N_i\) by ordinary solid tori
\(\{N'_1,\dots,N'_n\}\) so that a meridian of each \(N'_i\) is sewn
to \(Q_i\).  We will remember where the meridian \(m_i\) of \(N_i\)
is on \(\bd N'_i\).  Let \(N_0\) be a fibered solid torus in \(M'\)
that contains all the \(N'_i\) in its interior.  Note that
\(M_0=M'-\inter N_0\) is the classifying space for the class of
\(M\).  Our diagram for \(M\) will consist of the fibered space
\(V=N_0-\bigcup\inter N'_i\) together with certain curves on the
boundary components of \(V\) and (if \(M\) is oriented) an
orientation of \(V\).  The fiber structure of \(V\) is simply the
product of an \(n\) punctured disk with \(S^1\).  On each \(\bd
N'_i=\bd N_i\) we will preserve \(m_i\) the meridian of \(N_i\).  On
\(\bd N_0=\bd M_0\) we will preserve a curve \(Q_0\) that is the
boundary of a section for \(M_0\).  This wil be ambiguous if \(M\)
is non-orientable.  Note that \(Q_0\) is a crossing curve on \(\bd
N_0\) for the fiber structure.  If \(M\) is oriented, then \(M\) and
\(M'\) share the space \(M_0\) (they actually share more) so \(M_0\)
and thus \(M'\) inherits an orientation from \(M\).  Then \(N_0\)
and thus \(V\) inherit an orientation from \(M'\).  The diagram is
the fibered space \(V\) with its orientation, if available, together
with the curves \(Q_0\) and the \(m_i\).

We now need to argue that \(V\) recovers all the information in the
symbol for \(M\) that is to the right of the vertical bar.  That is,
the numbers \(b\), \(s\), the \(\mu_i\) and the \(\beta_i\).  We
will see that the ambiguity in \(Q_0\) will not interfere with these
determinations.

By Theorem \ref{FiberMeridian}, there is only one fiber type of
fibered solid torus that can be sewn by a fiber preserving
homeomorphism of the boundary to \(\bd N'_i\) so that a meridian is
sewn to \(m_i\).  This fibered solid torus (which we may as well
refer to as \(N_i\)) can be oriented consistently with \(V\) (if an
orientation of \(V\) is available).  This determines crossing
invariants \((\mu_i,\beta_i)\) and sometimes a crossing curve
\(Q_i\) in the usual way.  This gives us the number \(s\) of
exceptional fibers of index 2, the \(\mu_i\) and the \(\beta_i\).
We only need \(b\), the obstruction to the section.

The crossing curve \(Q_i\) is determined up to orientation unless
\(s\ne0\) and \(V\) is non-orientable in which case \(Q_i\) is
ambiguous.  However, when \(V\) is non-orientable and \(s\ne0\),
there is no need to determine \(b\).  We thus assume that \(V\) is
oriented, or that \(s=0\).

Assume first that \(V\) is oriented.  If an orientation of an
arbitrary fiber \(H\) in \(V\) is chosen, then all fibers can be
oriented consistently with \(H\).  Then the orientation of each
crossing curve \(Q_i\) is determined by requiring that
\(\hint(H_i,Q_i)=+1\) where \(H_i\) is a fiber in \(\bd N_i\) and
\(\bd N_i\) is oriented to agree with the orientation on \(N_i\).
The curves \(Q_i\) are the meridians of the ordinary solid tori
\(N'_i\).  If the \(N'_i\) are restored by adding them to \(V\),
then \(N_0\) is recovered.  We orient \(\bd N_0\) consistently with
\(N_0\).  Since \(m_0\) and \(Q_0\) are crossing curves on \(\bd
N_0\), we can orient \(m_0\) and \(Q_0\) by requiring that
\(\hint(H_0,m_0)=\hint(H_0,Q_0)=+1\) for some fiber \(H_0\) in \(\bd
N_0\).  The orientations used on the \(\bd N_i\) were inherited from
the \(N_i\) and are the reverse of the orientations that they would
inherit from \(V\).  Since \(m_0\) is a meridian for \(N_0\) in
which the \(N'_i\) are solid tori with meridians \(Q_i\), we have
\(m_0=\sum Q_i\) in \(H_1(V)\).  However, \(m_0\) and \(Q_0\) were
oriented as required to determine \(b\) as \(\hint(m_0,Q_0)\) and,
from \ref{ObstructionMeridian}, determine \(m_0\) as \(bH_0+Q_0\).
Thus \(\sum Q_i=bH_0+Q_0\) and we have shown that \(-Q_0+\sum Q_i\)
is homologous in \(V\) to a multiple of a fiber.  The number \(b\)
is that multiple.  If the opposite orientation is chosen for an
arbitrary fiber \(H\), then all orientations of \(Q_0\), \(m_0\) and
the \(Q_i\) reverse and the calculations remain the same.

If \(M\) is non-orientable and \(s=0\), then the curves \(Q_i\) and
\(m_0\) are still determined, but their orientations are not.  Also,
\(Q_0\) is determined only up to an even multiple of the fiber.
However, \(b\) is still determined mod 2 and this is all we need.

We illustrate a use for the diagram by showing the effect of
reversal of orientation on the invariants of a closed, connected,
oriented Seifert fibered space \(M\).  Let \(-M\) denote the space
\(M\) endowed with the opposite orientation of \(M\).  Let the \(M\)
be determined by \[(O,o,g\mid b,\mu_1,\beta_1,\dots,\mu_n,\beta_n)\]
or \[(O,n,k\mid b,\mu_1,\beta_1,\dots,\mu_n,\beta_n).\] Reversing
the orientation of \(M\) changes neither the orbit surface nor the
classifying homomorphism.  We must look at the effects on the
numbers to the right of the vertical bar.  If \(V\), \(m_i\) and
\(Q_0\) make up the diagram for \(M\), then we obtain the diagram
for \(-M\) by reversing the orientation of \(V\) to get \(-V\).  If
a fiber \(H\) was used to calculate the invariants for \(M\), then
it will be convenient to use the fiber \(-H\) to calculate the
invariants \(b'\), \(\mu'_i\) and \(\beta'_i\) for \(-M\).

Using \(-H\) instead of \(H\) means that we can still use the
orientations on the \(m_i\) that were used for \(M\).  We have
\(\mu'_i=\hint(-H_i,m_i)=\mu_i\) since \(\hint\) is the negative of
the intersection pairing used with \(M\).  If \(Q_i\) is a crossing
curve appropriate for \(M\), we must determine how to alter it to
give a curve \(Q'_i\) apropriate for \(-M\).  We have
\(Q'_i=Q_i+xH_i\).  This is correct since in \(-V\) we get the
correct intersection with \(-H\) without reversing the direction of
\(Q_i\).  We want \(\beta'_i=\hint(m_i,Q'_i)\) or the coefficient of
\(-H_i\) in \(m_i\) to be in the interval \([0,\mu_i)\).  But
\begin{align*} m_i &=\beta_i H_i+\mu_i Q_i \\ &=\beta_i
H_i+\mu_i(Q'_i-xH_i)\\ &=(\beta_i-\mu_i x)H_i-\mu_i Q'_i \\ &=(\mu_i
x-\beta_i)(-H_i)-\mu_i Q'_i \end{align*} so that we get the the
right value when \(x=1\) and \(\beta'_i=\mu_i-\beta_i\).

We need to determine \(b'\).  We have \(Q'_0=Q_0\) since we have
replaced \(H\) by \(-H\).  We determine \(b'\) from the equality
\(-Q'_0+\sum Q'_i=b'(-H)\).  The left side is
\(-Q_0+\sum(Q_i+H_i)=-Q_0+nH+\sum Q_i\) since all fibers are
homologous in \(V\).  But this is \(bH+nH=(-n-b)(-H)\) so
\(b'=-n-b\).  Note that the functions \(f_i(\beta_i)=\mu_i-\beta_i\)
and \(f(b)=-n-b\) have period two.  We have the following.

\begin{thm} Let \(M\) be a closed, connected, oriented Seifert
fibered space determined by either \[(O,o,g\mid
b,\mu_1,\beta_1,\dots,\mu_n,\beta_n)\] or \[(O,n,k\mid
b,\mu_1,\beta_1,\dots,\mu_n,\beta_n).\] Then the Seifert fibered
space with the opposite orientation is detmined by \[(O,o,g\mid
-n-b,\mu_1,\mu_1-\beta_1,\dots,\mu_n,\mu_n-\beta_n)\] or
\[(O,n,k\mid -n-b,\mu_1,\mu_1-\beta_1,\dots,\mu_n,\mu_n-\beta_n).\]
\BlackBox\end{thm}

\chapter{Topology of Seifert fibered spaces}
\setcounter{equation}{0}

Most of this chapter is based on pages 83--101 of Chapter VI of W.
Jaco's book ``Lectures on Three-manifold topology,'' Regional
Conference Series in Mathematics, number 43, AMS, Providence, 1980.
Many of the details were taken from Chapter II, Sections 3 and 4 of
W. Jaco and P. Shalen: ``Seifert fibered spaces in 3-manifolds,''
Mem. Amer. Math. Soc., number 220, (1979).  The material on covers
of Seifert fibered spaces is taken from Section 9 of Seifert's
paper, identified in the opening paragraph of Chapter 1.

We wish to classify Seifert fibered spaces up to homeomorphism.  We
need invariants that do not depend on the fibers, but that give
topological information.  The most powerful invariant is the
fundamental group.  In this chapter, we calculate the fundamental
group of Seifert fibered spaces, and use the results of the
calculation to gather information about the spaces.  In particular
we will classify the compact, connected Seifert fibered spaces up to
homeomorphism.  [The course ended before this was done.  Compact
Seifert fibered spaces with non-empty boundary were classified.]
Along the way we will figure out which Seifert fibered spaces are
aspherical, have incompressible boundary and have an incompressible
surface.

We need some standard concepts from the topology of 3-manifolds.  We
know that a 3-manifold is irreducible if every embedded 2-sphere
bounds a 3-cell.  We say that a 3-manifold is
\defit{\(P^2\)-irreducible} if it is irreducible and it contains no
embedded 2-sided projective plane.  The projective plane theorem of
Epstein, a generalization of the sphere theorem, says that a
3-manifold with non-trivial \(\pi_2\) fails to be
\(P^2\)-irreducible.  A \(P^2\)-irreducible 3-manifold has trivial
\(\pi_2\), and if it has non-empty boundary or infinite \(\pi_1\),
then its universal cover has trivial \(\pi_1\), trivial \(\pi_2\)
and trivial \(H_n\) for all \(n\ge3\).  Thus a \(P^2\)-irreducible
3-manifold with infinite \(\pi_1\) or with non-empty boundary is
aspherical.

A surface \(S\) is \defit{properly embedded} in a 3-manifold \(M\)
if it is embedded and \(\bd S=S\cap\bd M\).  A properly embedded
surface \(S\) is \defit{boundary parallel} in \(M\) if there is an
embedding of \(S\x I\) into \(M\) carrying \(S\x\{0\}\) onto \(S\)
and \((\bd S\x I)\cup (S\x\{1\})\) into \(\bd M\).  A surface in
\(M\) is said to be \defit{incompressible} if either it is a
2-sphere bounding no 3-cell, or is properly embedded, two sided and
not a disk and the inclusion induces an injection on \(\pi_1\), or
if it is a properly embedded disk that is not boundary parallel, or
if it is embedded in the boundary and is not a disk and the
includsion into \(M\) induces an injection on \(\pi_1\).  The
technicalities are to avoid certain trivial situations.

We say that a 3-manifold whose boundary components are all
incompressible is \defit{boundary irreducible}.  This is satisfied
vacuously if the boundary of the manifold is empty.  We say that a
properly embedded, two sided surface \(S\) in a 3-manifold \(M\) is
\defit{boundary compressible} if there is a disk \(D\) in \(M\) so
that \(D\cap S\) is an arc \(A\) in \(\bd D\), \(D\cap\bd M\) is an
arc \(B\) in \(\bd D\) with \(\bd D=A\cup B\) and so that \(A\) is
not boundary parallel in \(S\).  (Shift the definition of boundary
parallel down one dimension.)  A properly embedded surface is
\defit{boundary incompressible} if it is not boundary compressible.
Unfortunately, boundary irreducible and boundary incompressible were
named independently and have stuck.

A compact, orientable, irreducible 3-manifold with an incompressible
surface is called a \defit{Haken} manifold.  We will try to have all
our statements apply to both orientable and non-orientable
3-manifolds, but this will not be possible.  Towards the end of this
chapter, we will restrict ourselves to orientable manifolds.

We are now ready to study some of these properties.  Seifert fibered
spaces with boundary are easier to handle than closed Seifert
fibered spaces, and we start with some elementary observations about
bounded spaces.

\section{Bounded, compact Seifert fibered spaces}

\begin{lemma} Let \(T\) be a fibered solid torus.  Then a fiber of
\(T\) represents a non-trivial element of \(\pi_1(T)\).  \end{lemma}

\begin{demo}{Proof} Let \(T\) be determined by \(\nu/\mu\mod1\) with
\(\nu\) and \(\mu\) in reduced terms and with \(\mu>0\).  Then the
centerline of \(T\) represents a generator of \(\pi_1(T)\), and any
other fiber represents \(\mu\) times the generator.  \end{demo}

\begin{cor} Let \(T\) be a fibered solid torus.  Then a saturated
annulus in \(\bd T\) is an incompressible surface in \(T\).
\BlackBox\end{cor}

We need standard results about 3-manifolds which first need standard
results from group theory.  Let \(X\) be a simplicial complex and
let \(A\) and \(B\) be pairwise disjoint subcomplexes that form
connected, closed subsets of \(X\).  Let \(h:A\into B\) be a
homeomorphism.  Let \(X_h\) be formed from \(X\) by using \(h\) to
identify \(A\) with \(B\).  We make no assumptions on the
connectivity of \(X\), but do assume that \(X_h\) is connected.
Thus \(X\) has at most two components \(X_1\) and \(X_2\).  We
assume that \(\pi_1(A)\) injects into the fundamental group of the
component of \(X\) containing \(A\) and make a similar assumption
about \(B\).  If \(X\) is connected, then we say that \(\pi_1(X_h)\)
is an HNN extension of \(\pi_1(X)\) along \(\pi_1(A)\) and
\(\pi_1(B)\).  If \(X\) is not connected, then \(\pi_1(X_h)\) is the
free product with amalgamation of \(\pi_1(X_1)\) and \(\pi_1(X_2)\)
(along the amalgamating subgroups \(\pi_1(A)\) and \(\pi_1(B)\)).
These constructions can be given purely algebraic definitions so
that given isomorphic subgroups \(A\) and \(B\) (with a given
isomorphism bewteen them) of one or two groups \(G\) and \(H\), we
can form an HNN extension along \(A\) and \(B\) or a free product
with amalgamation along \(A\) and \(B\) depending on whether one or
two parent groups are involved.  We give the next lemma in
topological terms although it can be given in algebraic terms.

\begin{lemma} With \(X_h\), \(X\), \(X_1\), \(X_2\), \(A\) and \(B\)
as above, we have the following.  Each of \(\pi_1(X)\),
\(\pi_1(X_1)\), \(\pi_1(X_2)\), (whichever exist) and \(\pi_1(A)\)
and \(\pi_1(B)\) inject into \(\pi_1(X_h)\).  Any torsion element of
\(\pi_1(X_h)\) is conjugate to an element of \(\pi_1(X)\),
\(\pi_1(X_1)\) or \(\pi_1(X_2)\).  \end{lemma}

\begin{demo}{Proof} See pages 178--187 of {\it Combinatorial Group
Theory} by Lyndon and Schupp.  \end{demo}

\begin{cor} With the notation as above, \(\pi_1(X_h)\) has torsion
if and only if one of \(\pi_1(X)\), \(\pi_1(X_1)\) or \(\pi_1(X_2)\)
has torsion.  \BlackBox\end{cor}

\begin{lemma}\label{IrredSew} Let \(M\) be a \(P^2\)-irreducible
3-manifold, not necessarily connected.  Let \(F\) and \(G\) be
disjoint, connected, incompressible surfaces in \(\bd M\) not
necessarily in different components of \(M\) and let \(h\) be a
homeomorphism from \(F\) to \(G\).  Let \(M_h\) be the 3-manifold
obtained by identifying each \(x\in F\) to \(h(x)\) and assume that
\(M_h\) is connected..  Then \(M_h\) is \(P^2\)-irreducible, the
common image of \(F\) and \(G\) in \(M_h\) is incompressible in
\(M_h\), and the fundamental group of each component of \(M\)
injects into \(\pi_1(M_h)\).  \end{lemma}

\begin{demo}{Proof} Because of the previous lemma, only the first
conclusion needs proof.  The proof of the other conclusions could
also be done geometrically with similar techniques.

If \(S\) is a 2-sphere or 2-sided projective plane in \(M_h\), then
it can be put in general position with respect to the image of \(F\)
and \(G\) in \(M_h\) that we continue to refer to as \(F\).  No
circle of \(S\cap F\) can be orientation reversing on \(S\) since it
would be one sided on \(S\) and could not reside on the two sided
surface \(F\).  Thus every circle of \(S\cap F\) bounds at least one
disk on \(S\).  Since \(F\) is incompressible in \(M\), a circle of
\(S\cap F\) that is innermost on \(S\) bounds a disk in \(S\) and a
disk \(D\) in \(F\) that create a sphere whose only intersection
with \(F\) is \(D\).  This sphere can be regarded as a subset of
\(M\) where it bounds a 3-cell.  An isotopy of \(S\) across this
3-cell lowers the number of components of intersection of \(S\) with
\(F\).  Eventually \(S\) and \(F\) are disjoint and \(S\) is seen to
bound a 3-cell in \(M_h\) if \(S\) is a sphere, or not exist if
\(S\) is a 2-sided projective plane.  \end{demo}

\begin{lemma}\label{IrredWBd} Let \(M\) be a compact, connected
Seifert fibered space with non-empty boundary.  Then \(M\) is
\(P^2\)-irreducible, \(\pi_1(M)\) is torsion free and each fiber of
\(M\) represents a non-trivial element of \(\pi_1(M)\).  \end{lemma}

\begin{remark} The trivial observation that \(\pi_1(M)\) is infinite
under the hypotheses of the lemma is worth making at this point.
\end{remark}

\begin{demo}{Proof} There are pairwise disjoint arcs in the orbit
surface so that the result of cutting along the arcs yields a set of
disks with no more than one exceptional point in each.  The spaces
over these disks are fibered solid tori, and \(M\) is recovered from
these tori by sewings along saturated annuli in their boundaries.
However, a saturated annulus in the boundary of a fibered solid
torus is incompressible in the fibered solid torus, and each solid
torus is \(P^2\)-irreducible.  Thus \(M\) is \(P^2\)-irreducible.

The fact that \(\pi_1(M)\) is torsion free follows from the
construction just given, from the fact that the fundamental group of
a solid torus is torsion free and from the corollary above.

Note that each exceptional fiber has a multiple that is homotopic to
an ordinary fiber.  Also, any two ordinary fibers are homotopic.
Since \(\pi_1(M)\) is torsion free, we only have to show that one
ordinary fiber is non-trivial.  But this is true in a fibered solid
torus, and the result follows from the construction just given and
from Lemma \ref{IrredSew}.  \end{demo}

The proof of the next lemma contains arguments that will be used
repeatedly in the rest of the chapter.

\begin{lemma}\label{BdRed} Let \(M\) be a compact, connected Seifert
fibered space with a compressible boundary component.  Then \(M\) is
a fibered solid torus.  In particular, \(M\) has at most one
exceptional fiber.  \end{lemma}

\begin{remark} This lemma says that the only compact, connected
Seifert fibered space that fails to be boundary irreducible is a
fibered solid torus.  \end{remark}

\begin{demo}{Proof} Let \(C\) be the compressible boundary
component.  Either \(C\) is a torus or a Klein bottle and its
compression yields a 2-sphere.  Since \(M\) is irreducible, the
2-sphere bounds a 3-cell and \(M\) is homeomorphic to a solid torus
or a solid Klein bottle.  In either case \(\pi_1(M)\) is \(\ints\).
Let \(p:M\into G\) be projection to the orbit surface.  Since a
fiber of \(M\) is non-trivial in \(\pi_1(M)\) and \(p\) induces a
surjection on \(\pi_1\), we have that \(\pi_1(G)\) is torsion.  But
the only surface with boundary whose fundamental group is torsion is
a disk.  We will be done if we can show that there are fewer than
two exceptional fibers.

We argue that it will suffice to show that if there are exactly two
exceptional fibers then \(\pi_1(M)\) cannot be cyclic.  This will
rule out the possibility that there are exactly two exceptional
fibers.  If there are at least three exceptional fibers, then there
is an arc in \(G\) which splits \(G\) into a disk with exactly two
exceptional points, and another disk with the rest of the
exceptional points.  The preimage of this arc is an incompressible
annulus in \(M\) that separates \(M\) into two Seifert fibered
spaces that are joined along the annulus.  The fundamental groups of
these two spaces inject into the fundamental group of \(M\), and one
space has two exceptional fibers and orbit surface a disk.  Since we
will show that the space with two exceptional fibers has non-cyclic
\(\pi_1\) we are done.

We now assume that there are two exceptional fibers.  Let \(\alpha\)
be an arc in \(G\) that splits \(G\) into two disks with one
exceptional point each.  The annulus \(A\) over \(\alpha\) splits
\(M\) into two fibered solid tori \(T_1\) and \(T_2\) with defining
numbers \(\nu_1/\mu_1\) and \(\nu_2/\mu_2\).  Since the centerlines
of these fibered solid tori are exceptional fibers, we have
\(\mu_1>1\) and \(\mu_2>1\).  The splitting gives two copies \(A_1\)
and \(A_2\) of \(A\).  The fundamental group of \(A_i\) is carried
by an ordinary fiber of \(T_i\) that represents \(\mu_i\) times a
generator of \(\pi_1(T_i)\).  The fundamental group of \(M\) is thus
generated by two elements \(a\) and \(b\) and the only relation is
\(a^{\mu_1}=b^{\mu_2}\).  Since each \(\mu_i\) is at least 2, the
group is a free product with amalgamation of two copies of \(\ints\)
along a proper subgroup of each and is thus not cyclic.  \end{demo}

\begin{cor}\label{IncompTorus} Let \(M\) be a compact, connected
Seifert fibered space, and let \(F\) be a two sided, saturated torus
or Klein bottle in \(M\) conatinaing no exceptional fibers and that
is the preimage of a simple closed curve in the orbit surface that
does not bound a disk with fewer than two exceptional points.  Then
\(F\) is incompressible in \(M\).  \end{cor}

\begin{demo}{Proof} If false, then splitting \(M\) along \(F\) would
result in a fibered solid torus, which contradicts the hypothesis.
\end{demo}

\begin{remark} There can be saturated Klein bottles that contain
exceptional fibers and that are not saturated over simple closed
curves.  Consider the oriented \(I\)-bundle over the Klein bottle
fibered by circles that contain orientation reversing curves of the
Klein bottle.  This has orbit surface a disk with two exceptional
fibers of index 2.  The image in the orbit surface of the saturated
Klein bottle is an arc connecting the two exceptional points.
\end{remark}

\section{Fundamental groups of Seifert fibered spaces}

We would like to extend the results above (where possible) to closed
Seifert fibered spaces.  The questions are a little more delicate
especially since some of the generalizations are false.  The
3-sphere is a Seifert fibered space in which all fibers are trivial
in the fundamental group.  Also, \(S^2\x S^1\) is not irreducible.
We need the power of the fundamental group.

We show first that results such as in the previous section can be
obtained easily for a large class of Seifert fibered spaces.

\begin{lemma}\label{IrredClosedI} Let \(M\) be a closed, connected
Seifert fibered space so that \(M\) either has an orientable orbit
surface of genus at least one, or orbit surface \(S^2\) with at
least 4 exceptional fibers, or non-orientable orbit surface with at
least 2 crosscaps, or orbit surface \(P^2\) with at least 2
exceptional fibers.  Then \(M\) is \(P^2\)-irreducible, boundary
irreducible and has an embedded incompressible torus or Klein
bottle.  It follows that \(M\) has infinite \(\pi_1\) and is
aspherical.  \end{lemma}

\begin{demo}{Proof} If we can find an incompressible, saturated
torus or Klein bottle \(F\) in \(M\), then we can split \(M\) along
\(F\) and use Lemma \ref{IrredWBd} to conclude that the result of
the splitting is \(P^2\)-irreducible.  The incompressibility of
\(F\) would then imply that \(M\) is irreducible.  The other results
would follow because the fundamental group of \(F\) is infinite and
a \(P^2\)-irreducible 3-manifold with infinite fundamental group is
aspherical.

From Corollary \ref{IncompTorus}, we get our surface \(F\) if there
is a simple closed curve in the orbit surface that does not bound a
disk with fewer than two exceptional points in it.  The hypotheses
that we are given have been tailored so as to guarantee the
existence of such a curve.  \end{demo}

Lemma \ref{IrredClosedI} motivates much of what we do in this
section.  We first calculate the fundamental groups of Seifert
fibered spaces, and then apply our calculation to some of the spaces
singled out as problems by Lemma \ref{IrredClosedI}.

\subsection{Calculating the fundamental group}

We can use the realization construction to calculate the fundamental
group of a Seifert fibered space.  Let \(M\) be a compact, connected
Seiert fibered space, let \(G\) be the orbit surface, let
\(\{H_1,\dots,H_n\}\) be the exceptional fibers of \(M\) with
pairwise disjoint fibered solid torus neighborhoods
\(\{N_1,\dots,N_n\}\), crossing curves \(\{Q_1,\dots,Q_n\}\) and
crossing invariants \((\mu_1,\beta_1,\dots,\mu_n,\beta_n)\).  Let
\(b\) be the obstruction to a section for the orbit surface if \(M\)
is closed.  Let \(\phi\) be the classifying homomorphism.  Since
this is enough information to determine \(M\), it should be enough
to determine \(\pi_1(M)\).  Some of the crossing curves \(Q_i\) may
be ambiguous, but since these ambiguities do not make \(M\)
ambiguous, they will not make \(\pi_1(M)\) ambiguous.

Let \(M'\) be obtained from \(M\) by repacing each \(N_i\) be an
ordinary solid torus \(N'_i\) with meridian sewn to \(Q_i\).  Let
\(M_0\) be obtained from \(M'\) by removing the interiors of all the
\(N'_i\) and the interior of an ordinary solid torus \(N_0\) that is
disjoint from the \(N'_i\).  We get \(M_0\) from \(M\) by drilling
out all the exceptional fibers and one ordinary fiber.  The orbit
surface \(G_0\) of \(M_0\) is a compact surface with boundary.  It
is obtained from \(G\) by removing \(n+1\) disks.  We will label the
boundary components of \(G_0\) depending on their origin.  We let
\(\{c_1,\dots,c_n\}\) be the components that arise from the removed
\(N'_i\).  We let \(\{d_1,\dots,d_m\}\) be the boundary components
that arise from the boundary components of \(M\).  We let \(e\) be
the boundary component that arises from the removal of \(N_0\).
There are pairwise disoint arcs \(\{\alpha_1,\dots,\alpha_r\}\) in
\(G_0\) with boundaries in \(e\) so that if \(G_0\) is cut along
these arcs, then a single disk with holes results.  If \(G_0\) is
orientable of genus \(g\), then \(r=2g\).  If \(G_0\) is
non-orientable with \(k\) crosscaps, then \(r=k\).  There are
pairwise disjoint arcs \(\{\alpha_{r+1},\dots,\alpha_{r+m+n}\}\),
disjoint from the arcs just mentioned, so that each \(c_i\) and each
\(d_i\) is connected to \(e\) by one of these arcs.  If \(G_0\) is
cut along all of the \(\alpha_i\), then a single disk \(E\) results.
We recover \(G_0\) by sewing together the copies of the arcs.  We
recover \(M_0\) from the ordinary solid torus \(T\) over \(E\) by
sewing together saturated annuli that are the pre-images of the
copies of the arcs.

The fundamental group of \(G_0\) is free.  Each sewn arc corresponds
to a free generator represented by a loop that pierces the arc once.
If \(G_0\) is orientable, then the generators corresponding to
\(\{\alpha_1,\dots,\alpha_{2g}\}\) are usually grouped in pairs, a
pair for each handle, and are denoted \(\{a_1,b_1,\dots,a_g,b_g\}\).
If \(G_0\) is non-orientable, then the generators corresponding to
\(\{\alpha_1,\dots,\alpha_k\}\) are usually denoted
\(\{x_1,\dots,x_k\}\).  The generators corresponding to the
remaining \(m+n\) arcs are represented by the boundary components
\(c_i\) and \(d_i\) so we abuse notation and use
\(\{c_1,\dots,c_n\}\) and \(\{d_1,\dots,d_m\}\) to denote the
remaining generators.  The fundamental group of \(G_0\) is free on
these generators, but it is useful to record how \(e\) relates to
these generators.  If \(G_0\) is orientable, then
\begin{equation}e=\prod[a_i,b_i]\prod c_i\prod
d_i.\label{OSurfRel}\end{equation} If \(G_0\) is non-orientable,
then \begin{equation}e=\prod x_i^2\prod c_i\prod
d_i.\label{NOSurfRel}\end{equation}

The reconstruction of \(M_0\) starts with the ordinary solid torus
\(T\).  The fundamental group of \(T\) is free on one generator that
we denote \(h\) since it is carried by any fiber.  Each sewing of a
pair of annuli results in another generator and a relation that
comes from the fact that the annuli have non-trivial fundamental
groups.  We again abuse notation and label each generator with the
same letter as the generator of \(\pi_1(G_0)\) that is associated
with the arc under the identified annuli.  The identified annuli
have fundamental groups generated by conjugates of \(h\).  One
annulus can be thought of as being generated by \(h\) and the other
can be thought of a generated by \(h\) conjugated by the new
generator.  If the annuli are sewn to preserve the orientations of
the fibers, then the conjugate is identified with \(h\).  Otherwise
the conjugate is identified with \(h^{-1}\).  The way that the
orientations of the fibers are identified is given by classifying
homomorphism \(\phi\).  Thus if \(G_0\) is orientable, then
\(\pi_1(M_0)\) is generated by \(h\), the \(a_i\), \(b_i\), \(c_i\)
and \(d_i\) with relations \(y_ihy_i^{-1}=h^{\phi(y_i)}\) where
\(y\) is any of the letters \(a\), \(b\), \(c\) or \(d\).  Note that
\(\phi(c_i)=+1\) since each \(c_i\) bounds a disk when \(G_0\) is
viewed as a subset of \(G\).  If \(G_0\) is non-orientable, then
\(\pi_1(M_0)\) is generated by \(h\), the \(x_i\), \(c_i\) and
\(d_i\) with the same relations.

Note that the realization of \(M_0\) has a section of \(G_0\)
naturally embedded in it.  The generators of \(\pi_1(G_0)\)
correspond to the like named generators of \(\pi_1(M_0)\).  In
particular, the generators \(c_i\), \(d_i\) and \(e\) correspond to
curves on the boundary components of \(M_0\) that are the boundaries
of the section of \(G_0\).  Since these curves are crossing curves,
we have on each of the boundary components of \(M_0\), a pair of
curves that generate the fundamental group of that boundary
component.  One is \(h\) and the other is one of \(c_i\), \(d_i\) or
\(e\).  Note that if the boundary components \(c_i\) and \(e\) of
\(G_0\) are filled in with disks, then the orbit surface of \(M'\)
is obtained.  We obtain \(M'\) from \(M_0\) by sewing ordinary solid
tori to the boundary components of \(M_0\) over the \(c_i\) with
meridians sewn to the \(c_i\) and an ordinary solid torus to the
component over \(e\) with meridian sewn to \(e+bh\).  Note that this
trivializes each \(c_i\) in the above presentations and trivializes
the curve \(eh^b\) where \(e\) is as given in either \ref{OSurfRel}
or \ref{NOSurfRel}.

To reconstruct \(M\) instead of \(M'\), we have to sew fibered solid
tori to the boundary components over the \(c_i\), instead of
ordinary solid tori.  A meridian for the fibered solid torus sewn to
the component over \(c_i\) is to be sewn to
\(h^{\beta_i}c_i^{\mu_i}\).  We can now write out the presentation
for \(\pi_1(M)\).

\begin{thm}\label{PiOneSeifert} Let \(M\) be a compact, connected
Seifert fibered space with \linebreak \(m\) boundary components,
with \(n\) exceptional fibers with crossing invariants \linebreak
\((\mu_1,\beta_1,\dots,\mu_n,\beta_n)\), with section obstruction
\(b\), and with classifying homomorphism \(\phi\).  If the orbit
surface of \(M\) is orientable of genus \(g\), then \(\pi_1(M)\) has
the presentation \begin{equation}\begin{split} \langle
h,a_1,b_1,\dots,a_g,b_g,c_1,\dots,c_n,d_1,\dots,d_m\mid &
a_iha_i^{-1}=h^{\phi(a_i)}, \\ & b_ihb_i^{-1}=h^{\phi(b_i)}, \\ &
c_ihc_i^{-1}=h, \\ & d_ihd_i^{-1}=h^{\phi(d_i)}, \\ &
c_i^{\mu_i}h^{\beta_i}=1, \\ & \prod[a_i,b_i]\prod c_i\prod
d_ih^b=1\rangle.\end{split} \label{PiOneOrbl}\end{equation} If the
orbit surface of \(M\) is non-orientable with \(k\) crosscaps, then
\(\pi_1(M)\) has the presentation \begin{equation}
\label{PiOneNOrbl} \begin{split} \langle
h,x_1,\dots,x_k,c_1,\dots,c_n,d_1,\dots,d_m\mid &
x_ihx_i^{-1}=h^{\phi(x_i)}, \\ & c_ihc_i^{-1}=h, \\ &
d_ihd_i^{-1}=h^{\phi(d_i)}, \\ & c_i^{\mu_i}h^{\beta_i}=1, \\ &
\prod x_i^2\prod c_i\prod d_ih^b=1\rangle.  \end{split}
\end{equation} \BlackBox\end{thm}

Note that if \(\phi\) is the trivial homomorphism, then \(h\) is a
central element.  If \(\phi\) takes on the value \(-1\) on some
element, then \(h\) is not preserved under conjugation by that
element, but the cyclic subgroup generated by \(h\) is preserved.
Thus the cyclic subgroup generated by \(h\) is normal in
\(\pi_1(M)\).  The element \(h\) is represented by some ordinary
fiber in \(M\).  Another ordinary fiber in \(M\) is isotopic in
\(M\) to the one representing \(h\) and so represents an element
conjugate to \(h\) or \(h^{-1}\).  Since the cyclic subgroup
generated by \(h\) is normal in \(\pi_1(M)\), any other ordinary
fiber must represent either \(h\) or \(h^{-1}\) and generate the
same cyclic subgroup as that generated by \(h\).  This is summarized
in the next lemma.

\begin{cor} The fundamental group of a Seifert fiber space has a
unique cyclic normal subgroup (possibly trivial) which is generated
by any ordinary fiber.  \BlackBox\end{cor}

Note that this does not claim that there are no other cyclic normal
subgroups in \(\pi_1(M)\).

Using the notation of Theorem \ref{PiOneSeifert}, we let \(\langle
h\rangle\) denote the cyclic, normal subgroup of \(\pi_1(M)\)
generated by an ordinary fiber.  We can form the quotient group
\(\pi_1(M)/\langle h\rangle\).  This has one of the following
presentations depending on the orientability of the orbit surface of
\(M\): \begin{equation}\begin{split} \langle
a_1,b_1,\dots,a_g,b_g,c_1,\dots,c_n,d_1,\dots,d_m\mid &
c_i^{\mu_i}=1, \\ & \prod[a_i,b_i]\prod c_i\prod
d_i=1\rangle,\end{split}\label{FuchsOrbl} \end{equation} or
\begin{equation}\begin{split} \langle
x_1,\dots,x_k,c_1,\dots,c_n,d_1,\dots,d_m\mid & c_i^{\mu_i}=1, \\ &
\prod x_i^2\prod c_i\prod d_i=1\rangle.\label{FuchsNOrbl}
\end{split}\end{equation} Groups with these presentations are known
as Fuchsian groups about which a great deal is known.  They appear
in complex analysis.  Later we will use these groups to show that
some of the spaces excluded from Lemma \ref{IrredClosedI} are
\(P^2\)-irreducible.  We first give a lemma relating the size of the
fundamental group of some Seifert fibered spaces to the size of the
cyclic normal subgroup generated by an ordinary fiber and the size
of the corresponding Fuchsian quotient.

\begin{lemma}\label{InfinitePiOne} Let \(M\) be a compact,
connected, \(P^2\)-irreducible and boundary irreducible Seifert
fibered space with infinite fundamental group.  Let \(h\in\pi_1(M)\)
be represented by an ordinary fiber.  Then the cyclic normal
subgroup \(\langle h\rangle\) generated by \(h\) is infinite and the
Fuchsian group \(\pi_1(M)/\langle h\rangle\) is infinite.
\end{lemma}

\begin{demo}{Proof} Since \(M\) is \(P^2\)-irreducible with infinite
fundamental group it is aspherical.  Since it is finite dimensional,
its fundamental group is torsion free.  Thus \(h\) is trivial or has
inifnite order.  If \(h\) is trivial, then \(\pi_1(M)\) is a
Fuchsian group which must then be torsion free and infinite.  Also,
\(\bd M\) must be empty by Lemma \ref{IrredWBd}.  From the
presentation of a Fuchsian group, we know that a torsion free,
infinite Fuchsian group is the fundamental group of an aspherical
surface \(S\).  Thus \(M\) and \(S\) have the same homotopy type.
Since \(M\) is closed, \(H_3(M;\ints_2)\) is not zero and we have a
contradiction.

If \(\pi_1(M)/\langle h\rangle\) is finite, then a finite cover of
\(M\) has infinite cyclic fundmanetal group.  The cover will also be
aspherical and boundary irreducible.  However, this is impossible.
(This is standard.  There is a homotopy equivalence to a circle.
Make the map tranverse to a point and look at the preimage of the
point.  Make the preimage an incompressible surface.  It can only be
a non separating 2-sphere or disk violating either the
irreducibility or the boundary irreducibility.)  \end{demo}

\begin{cor} Let \(M\) be a compact connected, Seifert fibered space
with non-empty boundary that is not a fibered solid torus.  Let
\(h\in\pi_1(M)\) be represented by an ordinary fiber.  Then
\(\pi_1(M)\), the cyclic normal subgroup \(\langle h\rangle\)
generated by \(h\), and the Fuchsian group \(\pi_1(M)/\langle
h\rangle\) are all infinite.  \end{cor}

\begin{demo}{Proof} This follows from Lemma \ref{IrredWBd}, Lemma
\ref{BdRed} and Lemma \ref{InfinitePiOne}.  \end{demo}

We now consider some very special cases of Fuchsian groups.

\subsection{The dihedral groups and the group
\protect\(\ints_2*\ints_2\protect\)}

We will be interested in knowing which spaces are not
\(P^2\)-irreducible.  A first step will be to determine the spaces
that are not irreducible.  The fundamental groups of such spaces
will be free products.  (If there are homotopy spheres that are not
\(S^3\), then the free product may be trivial.)  We will show that
\(\ints_2*\ints_2\) is the only non-trivial free product that has a
non-trivial cyclic normal subgroup.  This is the ``infinite
dihedral'' group, and it and the finite dihedral groups turn out to
be relevant. [Actually, the finite ones may not be all that
relevant.  Read on.]

The dihedral group \(D_n\) of order \(2n\) has presentation
\[\langle x,t\mid x^2=t^n=1,\, xtx=t^{-1}\rangle\] or \[\langle
x,t\mid x^2=t^n=xtxt=1\rangle\] or \begin{equation}\langle x,t\mid
x^2=t^n=(xt)^2=1\rangle.\label{Dnxt}\end{equation} Letting \(y=xt\)
so that \(t=xy\) we get the alternate presentation
\begin{equation}\langle x,y\mid
x^2=y^2=(xy)^n=1\rangle.\label{Dnxy}\end{equation} The infinite
dihedral group \(D_{\infty}\) eliminates the cylicity in the
generator \(t\) and has presentations \begin{equation}\langle
x,t\mid x^2=(xt)^2=1\rangle\label{Dinfxt}\end{equation} and
\begin{equation}\langle x,y\mid
x^2=y^2=1\rangle\label{Dinfxy}\end{equation} using the same
substitutions.  However the last is just a presentation of
\(\ints_2*\ints_2\).  We stretch notation and use \(D_n\) to refer
to both the finite and infinite dihedral groups.  We will specify
\(n\) finite when needed and use \(D_{\infty}\) when needed.

Note that \(D_n\) for finite \(n\) is isomorphic to the Fuchsian
group \[\langle c_1,c_2,c_3\mid
c_1^2=c_2^2=c_3^n=c_1c_2c_3=1\rangle\] since \(c_3=(c_1c_2)^{-1}\)
and \((c_1c_2)^{-n}=1\) if and only if \((c_1c_2)^n=1\).  Similarly
\(D_{\infty}\) is isomorphic to the Fuchsian group \[\langle c_1,
c_2,d_1\mid c_1^2=c_2^2=c_1c_2d_1=1\rangle.\]

The subgroup \(\langle t\rangle\) of \(D_n\) generated by \(t\) is
cyclic and normal as are all of the subgroups of \(\langle
t\rangle\).  We investigate the other normal subgroups.

It is clear from the presentations \ref{Dnxy} and \ref{Dinfxy} that
elements of \(D_n\) are words that alternate in the letters \(x\)
and \(y\).  A word with an even number of letters is a power
(positive or negative) of \(t=xy\), and a word with an odd number of
letters is a conjugate of the letter in the center of the word,
either \(x\) or \(y\).  Thus words with an odd number of letters are
involutions.  The elements \(x\) and \(y\) are not conjugate since
setting \(x=1\) in either \ref{Dnxy} or \ref{Dinfxy} does not kill
the element \(y\).  Note that setting \(x=1\) gives a presentation
for \(\ints_2\), so the subgroup of \(D_n\) normally generated by
\(x\) is of index 2 in \(D_n\).  Similarly for \(y\).  Since
\(\langle t\rangle\) is normal in \(D_n\) it is a union of conjugacy
classes.  Thus \(D_n-\langle t\rangle\) holds two conjugacy classes,
the conjugates of \(x\) and the conjugates of \(y\).  (So
\(D_n-\langle t\rangle\) consists entirely of involutions.)  Thus a
normal subgroup of \(D_n\) either is a subgroup of \(\langle
t\rangle\), or contains the conjugates of \(x\) or \(y\) or both.
If a normal subgroup contains conjugates of both \(x\) and \(y\),
then it is all of \(D_n\).  If it contains conjugates of \(x\), then
it contains \(x\) and \(yxy\) and thus \(xyxy=t^2\).  Thus the index
2 normal subgroup containing \(x\) consists of words in \(x\) and
\(y\) that either have odd length and center letter \(x\) or have
length a multiple of 4.  Similarly for \(y\).  Neither of these
normal subgroups is cyclic in \(D_n\) since neither \(x\) nor \(y\)
commute with \(t^2=xyxy\).  Thus the cyclic normal subgroups of
\(D_n\) are precisely the subgroups of \(\langle t\rangle\).  Note
that modding out by one of these cyclic normal subgroups means
adding the relation \(t^d=1\) to \ref{Dnxt} for some divisor \(d\)
of \(n\) or to \ref{Dinfxt} for any \(d\).  The quotient is the
dihedral group \[\langle x,t\mid x^2=t^d=(xt)^2=1\rangle \cong
\langle x,y\mid x^2=y^2=(xy)^d=1\rangle\] if \(d\ne1\) or
\(\ints_2\) if \(d=1\).

We have the following uniqueness result.

\begin{lemma}\label{CyclicNormal} The only non-trivial free product
with a non-trivial cyclic normal subgroup is \(D_{\infty}\).
\end{lemma}

\begin{demo}{Proof} Assume that \(G=A*B\) is a non-trivial free
product and that \(G\) has a non-trivial cyclic normal subroup
\(N\).  Every element of \(G\) is of the form \(a_1b_1\dots a_nb_n\)
where each \(a_i\) is in \(A\) and each \(b_i\) is in \(B\).  No
\(a_i=1\) except possibly \(a_1\) and no \(b_i=1\) except possible
\(b_n\).  Such an element is trivial if and only if \(n=1\) and
\(a_1=b_1=1\).  Assume such an element is the generator of \(N\).

If \(N\) is finite cyclic, then it is conjugate into one of \(A\) or
\(B\).  Snce it is normal, it would be in one of \(A\) or \(B\).
However, conjugating an element of \(A\) by an element of \(B\)
produces an element that is not in \(A\).  Thus \(N\) is not finite.

Since \(N\) is normal all rotations of the generator are in \(N\)
and are non-trivial.  Thus by passing to a subgroup of \(N\), we may
assume that the generator \(g\) of \(N\) is cyclically reduced and
has both \(a_1\ne1\) and \(b_n\ne1\).  Now positive powers of \(g\)
start with \(a_1\) and negative powers start with \(b^{-1}_n\).  The
length (measured in numbers of letters) of \(g^k\) is \(2n|k|\).
Conjugating \(g\) by \(a_1b_1\) produces \(a_2b_2\dots
a_nb_na_1a_2\) which can only equal \(g\).  Thus each \(a_i=a_1\)
and each \(b_i=b_1\).  Thus \(g=(a_1b_1)^n\) and every element of
\(N\) is of the form \((a_1b_1)^{nk}\).  Conjugating \(g\) by
\(a_1\) produces \((b_1a_1)^n\) which can only equal \(g^{-1}\) so
\(a_1=a^{-1}_1\) and \(b_1=b^{-1}_1\).  Assume that \(A\) has an
element \(c\) other than \(a_1\).  We know \(c^{-1}\ne a_1\) as
well.  Conjugating \(g\) by \(c\) produces \(a'(b_1a_1)^{n-1}b_1c\)
where \(a'=c^{-1}a_1\ne1\).  This is no element of \(N\).  Thus
\(a_1\) is the only element of \(A\).  Similarly, \(b_1\) is the
only element of \(B\).  Thus \(G\) is \(\ints_2*\ints_2\) or
\(D_{\infty}\).  \end{demo}

Next we take up a larger class of Fuchsian groups.

\subsection{The triangle groups}

Let \(p\), \(q\) and \(r\) be integers each greater than 1.  The
Fuchsian group \[\langle c_1,c_2,c_3\mid
c_1^p=c_2^q=c_3^r=c_1c_2c_3=1\rangle \] is isomorphic to \[\langle
c_1,c_2\mid c_1^p=c_2^q=(c_1c_2)^r=1\rangle.\] The presentations
make it clear that permutations of the subscripts \(\{123\}\) give
isomorphic groups so we may assume that \(p\le q\le r\).  These
groups are called \defit{triangle groups}.  The reason for this
terminology will become clear.  (Actually, these groups should
probably be called ``hexagonal groups.''  Another set of groups has
better claim to the title of triangle groups.  However, the notation
seems to be standard.)

Let \(\Gamma(p,q,r)\) have the presentation \(\langle a,b,c\mid
a^p=b^q=c^r=abc=1\rangle\) or equivalently \(\langle a,b\mid
a^p=b^q=(ab)^r=1\rangle\).  We will use \(\Gamma\) for short when
the specific values of \((p,q,r)\) are not important or are clear
from the context.  Let \(\Delta\) be a triangle with vertex angles
\(\pi/p\), \(\pi/q\) and \(\pi/r\).  Such a triangle can be realized
with geodesic sides on the sphere \(S^2\) if the sum of the angles
is greater than \(\pi\), on the Euclidean plane \(E^2\) if the sum
of the angles is exactly \(\pi\) and on the hyperbolic plane \(H^2\)
if the sum of the angles is less than \(\pi\).  This is determined
by comparing the sum
\begin{equation}\frac{1}{p}+\frac{1}{q}+\frac{1}{r}\label{TriangleSum}\end{equation}
to 1.  The triangle is realized on \(S^2\), \(E^2\) or \(H^2\)
respectively as \ref{TriangleSum} is greater than, equal to or less
than 1.  Since each of \(p\), \(q\) and \(r\) is at least 2, we can
easily describe the triples \((p,q,r)\), up to permutation, for
which \ref{TriangleSum} is no less than 1.  The sum is greater than
1 for \((2,2,r)\), \((2,3,3)\), \((2,3,4)\) and \((2,3,5)\).  The
sum is exactly 1 for \((2,3,6)\), \((2,4,4)\) and \((3,3,3)\).

\FIGURE{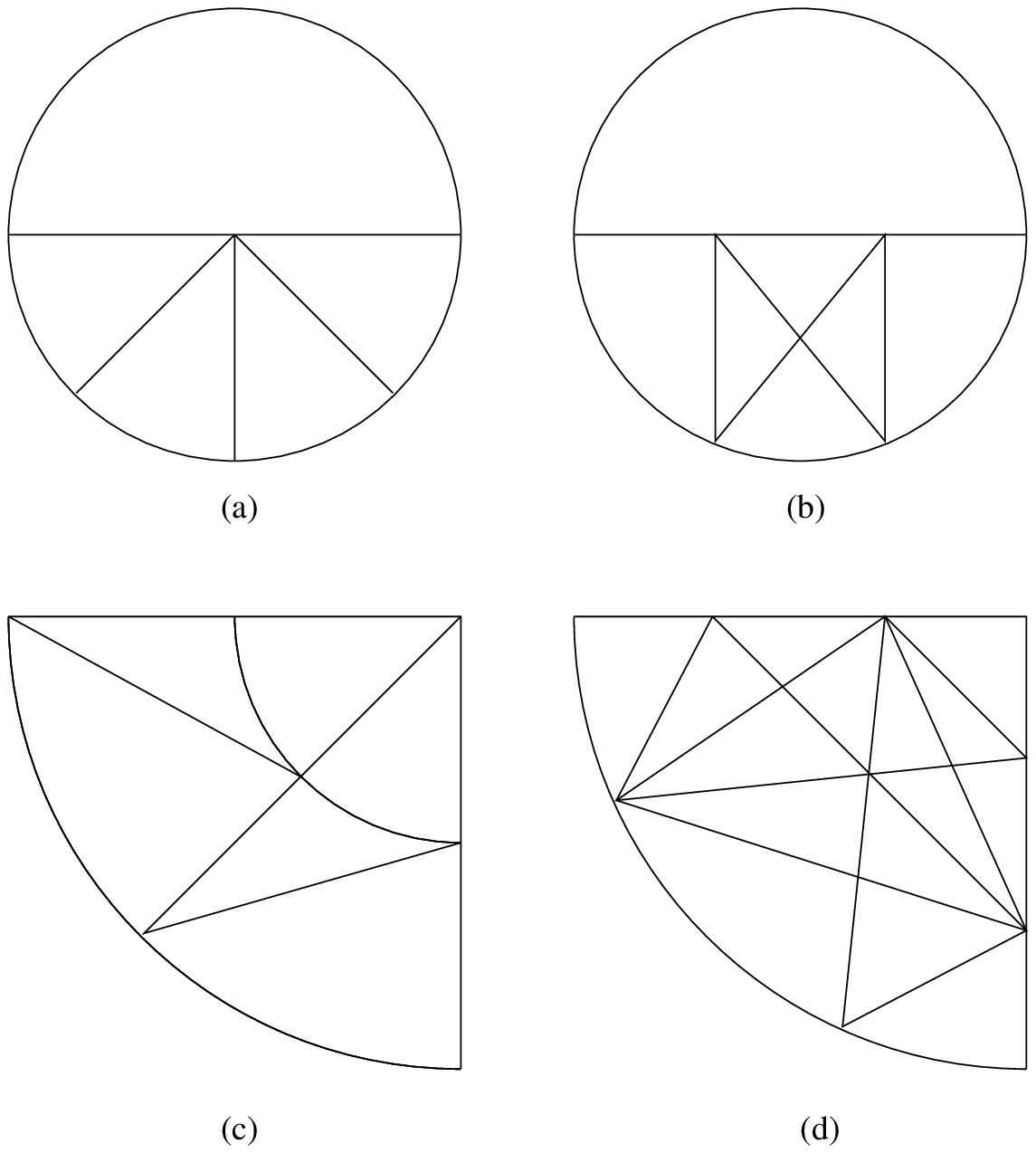}{Tessellations of \(S^2\)}

\FIGURE{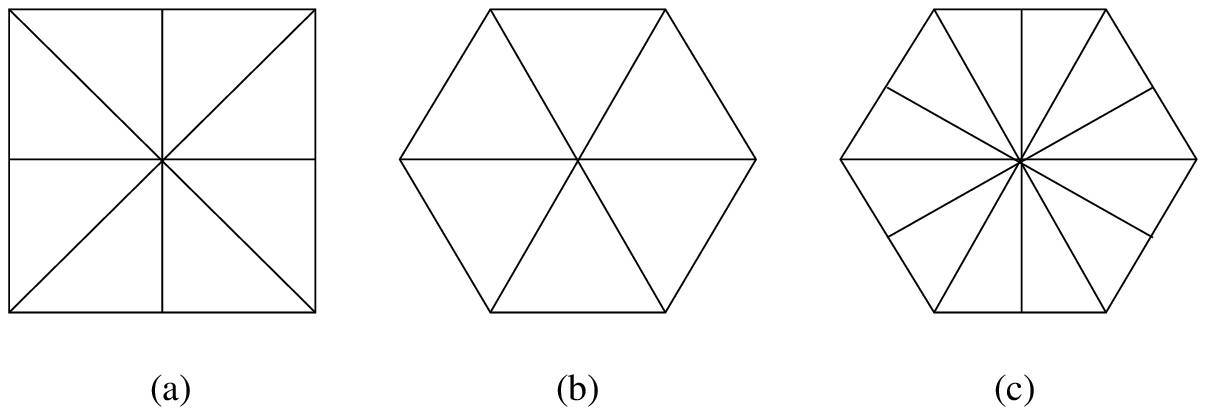}{Tessellations of \(E^2\)}

We will show that if we reflect \(\Delta\) across its edges
repeatedly, we will tessellate \(S^2\), \(E^2\) or \(H^2\) with
copies of \(\Delta\).  Figure \ref{SphereTessellations} shows how to
tessellate \(S^2\) for the values \((2,2,4)\), \((2,3,3)\),
\((2,3,4)\) and \((2,3,5)\) of \((p,q,r)\).  In Figure
\ref{SphereTessellations}, (a) and (b) show a quadrant of \(S^2\) as
viewed from above, while (c) and (d) show an octant of \(S^2\) as
viewed from above.  Figure \ref{PlaneTessellations} shows basic
units of tessellations of \(E^2\) for the allowable values of
\((p,q,r)\).  For examples of tessellations of \(H^2\) see various
etchings of M. C. Escher.

The original argument that I had for what follows was flawed.  The
argument below is mostly derived from B. Maskit, ``On Poincar\'e's
theorem for fundamental polygons,'' {\it Advances in Math.}, {\bf 7}
(1971), 219--230.

That the triangles successfully tessellate their respective spaces
is a covering argument.  When we do this we will be analyzing a
group that contains \(\Gamma\) as a subgroup of index 2.  (This
larger group is one that more naturally claims the title of a
triangle group.)  We could do the argument direclty with \(\Gamma\)
(using hexagons), but it is easier to work with the triangles, the
larger groups is interesting in its own right, the use of triangles
shows why \(\Gamma\) is called a triangle group, and once the
argument for the larger group is done, the argument for \(\Gamma\)
is seen as virtually identical.

We actually use two isomorphic groups in the argument.  At first it
is not clear that the groups are isomorphic.  One is defined
geometrically as as group of reflections, and the other by a
presentation.  We will get our tessellation by showing that an
abstract complex defined from the presentation is a covering space
for the space acted on by the geometric group.  We get as a further
consequence that the two groups are isomorphic.  This will be used
later to relate \(\Gamma\) to a subgroup of the geometric group.

We start with the group of reflections.  Let \(S\) be whichever of
\(S^2\), \(E^2\) or \(H^2\) contains \(\Delta\).  We start a pattern
of abusing notation and use \(x\), \(y\) and \(z\) to label the
edges of the triangle \(\Delta\).  We use \(x\) to label the edge
opposite the vertex of angle \(\pi/q\), we use \(y\) to lable the
edge opposite the vertex of angle \(\pi/r\), and we use \(z\) to
label the edge opposite the vertex of angle \(\pi/p\).  We let \(G\)
be the group of isometries of \(S\) generated by reflections in the
sides of \(\Delta\).  That is, reflections in the lines that contain
the sides of \(\Delta\).  Continuing our abuse of notation, we let
\(x\), \(y\) and \(z\) denote reflection of \(S\) across the lines
of \(S\) that contain the edges \(x\), \(y\) and \(z\) of \(\Delta\)
respectively.  Thus \(G\) is generated by the elements \(x\), \(y\)
and \(z\).  We let elements of \(G\) act on the left.  Since
elements of \(G\) are isometries, they are determined by their
actions on \(\Delta\).  If \(g\) is an element of \(G\), then we let
the triangle \(g(\Delta)\) have its edges labeled by \(x\), \(y\)
and \(z\) as carried over by the action of \(g\).  (This is not so
much a labeling of a subset of \(S\), but a labeling of an image of
a specific map.  We are not going to be concerned with the fact that
two elements of \(G\) might have the same images when restricted to
\(\Delta\).)

Note that each of \(x\), \(y\) and \(z\) is an involution in \(G\)
and that each is non-trivial.  Also, \(xy\) is a rotation about the
common point of the lines of reflection for \(x\) and \(y\).  That
is, \(xy\) is a rotation about the vertex of \(\Delta\) that is
common to the sides \(x\) and \(y\).  This is the vertex opposite
the side \(z\), and in our notation it is the vertex of angle
\(\pi/p\).  The element \(xy\) of \(G\) is seen to be a rotation
about this vertex by an angle \(2\pi/p\).  Thus \(xy\) has order
exactly \(p\) in \(G\).  Similarly, \(yz\) has order \(q\) and is a
rotation by \(2\pi/q\) about the vertex of \(\Delta\) of angle
\(\pi/q\), and \(zx\) has order \(r\) and is a rotation by
\(2\pi/r\) about the vertex of \(\Delta\) of angle \(\pi/r\).  We
use the relations that we have discovered about \(x\), \(y\) and
\(z\) in \(G\) to define an abstract group.

Let \(G_0\) be the group with presentation \[\langle x,y,z\mid
x^2=y^2=z^2=(xy)^p=(yz)^q=(zx)^r=1\rangle.\] The group \(G\) is a
homomorphic image of \(G_0\) by a homomorphism taking the generators
of \(G\) to the elements of the same letter in \(G\).  We will see
that the homomorphism just defined is actually an isomorhpism
between \(G_0\) and \(G\).  From what we know of \(G\), we know that
each of \(x\), \(y\) and \(z\) is a non-trivial involution in
\(G_0\) and that \(xy\), \(yz\), and \(zx\) have orders \(p\), \(q\)
and \(r\) respectively in \(G_0\).

Let \(\Omega\) be the graph of the group \(G_0\) with respect to the
generating set \(x\), \(y\) and \(z\).  Recall that the edges of the
graph \(\Omega\) are also labeled with the letters \(x\), \(y\) and
\(z\).  Also, recall that the edges of \(\Delta\) are labeled with
the letters \(x\), \(y\) and \(z\).  We use the graph \(\Omega\) to
define a simplical complex using copies of the triangle \(\Delta\).

Let \(K=G_0\x \Delta\).  Let \(g\) and \(gw\) be elements in \(G_0\)
where \(w\) is one of \(x\), \(y\) or \(z\).  Thus \(g\) and \(gw\)
are connected by an edge in \(\Omega\).  (Remember that each of
\(x\), \(y\) and \(z\) is its own inverse.)  In \(K\), we will now
identify the edges labeled \(w\) in the triangles \((g,\Delta)\) and
\((gw,\Delta)\).  Since \(\Delta\) is a geometric object, we can
insist that the identification is done by an isometry.  The only
thing that we have to specify is the orientation of the sewing.  We
require that the sewing be done so that the triangles \((g,\Delta)\)
and \((gw,\Delta)\) are seen as reflections of each other across the
common edge \(w\).  Another way of specifying this is to require
that for each endpoint of \(w\), the edges of \((g,\Delta)\) and
\((gw,\Delta)\) other than \(w\) that touch that endpoint (after
identification) have the same label.  See Figure
\ref{TriangleTessellation} for an example where \(w=x\).  We let
\(F\) represent the resulting complex.  Since the graph of \(G_0\)
is connected, we get that the complex \(F\) is connected.

\FIGURE{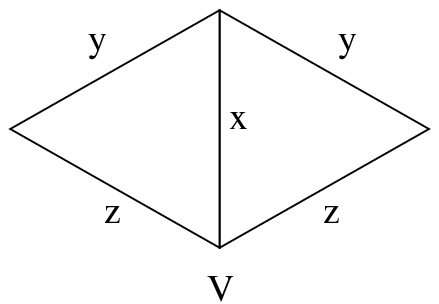}{Identified triangles}

Since every edge of a triangle in \(K\) is labeled with one of
\(x\), \(y\) and \(z\), and since each of \(x\), \(y\) and \(z\) is
a non-trivial involution in \(G_0\), we know that every edge of a
triangle in \(K\) is identified with exactly one other edge of a
triangle in \(K\).  Thus \(F\) is a surface except possibly at the
images of the vertices.  However, each of \(xy\), \(yz\) and \(zx\)
have finite order in \(G_0\).  This implies that each vertex has a
surface neighborhood.  (Consider the implications of the fact that
\(zx\) has finite order on the vertex labeled \(V\) in Figure
\ref{TriangleTessellation}.)  Thus \(F\) is a surface without
boundary.

If \(g\) is an element of \(G\), and \(w\) is one of \(x\), \(y\) or
\(z\), then the triangles \(g(\Delta)\) and \(gw(\Delta)\) are the
images under \(g\) of \(\Delta\) and \(w(\Delta)\) and are thus
reflections of each other in the edge \(w\).  If we now map
\((g,\Delta)\) in \(K\) to \(g(\Delta)\) and \((gw,\Delta)\) in
\(K\) to \(gw(\Delta)\) by isometries that preserve the labeled
edges, then these maps commute with the identifications that create
\(F\) from \(K\).

We now map each \((g,\Delta)\) in \(K\) to \(g(\Delta)\) in \(S\) by
an isometry that preserves the labeled edges.  This gives a map from
\(F\) to \(S\) that is an isometry on each triangle.  Further, this
map is a local isometry at each point in the interior of an edge
since the triangles that meet along an edge in \(F\) are carried to
triangles that are reflections of each other in \(S\).  Also, the
map is a local isometry at each vertex, since the number of
triangles meeting at a vertex of \(F\) is carried to the same number
of triangles meeting at a vertex in \(S\).  We thus have that the
map is a local isometry on the star of each vertex.

We now have a local homeomorphism \(f:F\into S\) and we wish to show
that it is a covering projection.  We must show that it evenly
covers.  In the surface \(F\), there is an \(\epsilon\) so that
every ball of radius \(\epsilon\) is contained in the star of some
vertex.  If \(f\) is not onto and \(a\) is a limit point of the
image, then there is a point \(b\) in \(F\) with \(f(b)\) less than
\(\epsilon\) from \(a\).  But the \(\epsilon\) ball about \(b\) maps
isometrically to \(S\) and will contain a point mapping to \(a\).
Thus the image of \(f\) is closed.  Since \(f\) is a local
homeomorphism, its image is also open and must be all of \(S\).  Now
if \(a\) is any point in \(S\), then the \(\epsilon\) ball about
\(a\) is evenly covered by \(f\).  Thus \(f\) is a covering
projection.

The space \(S\) is simply connected and the surface \(F\) is
connected.  Thus \(f\) must be a homeomorphism and \(S\) is
tessellated by the images of \(\Delta\) under the action of \(G\).

We know that the homomorphism from \(G_0\) to \(G\) is an
epimorphism since \(x\), \(y\) and \(z\) generate \(G\).  We see
that this homomorphism is an isomorphism by noting that different
elements \(g\) and \(g'\) of \(G_0\) correspond to different
triangles \((g,\Delta)\) and \((g',\Delta)\) in \(F\) which map to
different triangles \(g(\Delta)\) and \(g'(\Delta)\) of \(S\) (since
\(h\) is a homeomorphism) so that \(g\) and \(g'\) must be different
elements of \(G\).  We now drop the notation \(G_0\) and use only
\(G\) from now on.

We are now ready to consider \(\Gamma\).

The subgroup \(G_1\) of \(G\) generated by \(A=xy\), \(B=yz\) and
\(C=zx\) satisfies the relations \(A^p=B^q=C^r=ABC=1\) and so is a
homomorphic image of \(\Gamma\) by sending \(a\) to \(A\), \(b\) to
\(B\) and \(c\) to \(C\).  We must show that this gives an
isomorphism from \(\Gamma\) to \(G_1\).

The homomorphism already shows that \(a\), \(b\) and \(c\) have
exactly the orders \(p\), \(q\) and \(r\) respectively in \(\Gamma\)
because of the presentation of \(\Gamma\) and because these are
exactly the orders in \(G_1\).

Note that the action of \(A\) is to rotate by \(2\pi/p\) around the
vertex shared by the edges \(x\) and \(y\).  We again absue notation
and use \(A\) to represent the vertex shared by the edges \(x\) and
\(y\).  Similarly, \(B\) is the vertex shared by \(y\) and \(z\),
and \(C\) is the vertex shared bvy \(z\) and \(x\).  The action of
\(B\) in \(G_1\) is to rotate by \(2\pi/q\) about the vertex \(B\),
and the action of \(C\) in \(G_1\) is to rotate by \(2\pi/r\) about
the vertex \(C\).

Let \(\Delta\) and \(S\) be as above with the edges of \(\Delta\)
labeled \(x\), \(y\) and \(z\) as before.  For \(w\) one of \(x\),
\(y\) or \(z\), let \(\Delta_w\) be the union of the two triangles
in the barycentric subdivision of \(\Delta\) that share an edge with
\(w\).  We have that \(\Delta\) is the union of the three triangles
\(\Delta_x\), \(\Delta_y\) and \(\Delta_z\).  The three triangles
have disjoint interiors, and each \(\Delta_w\) is a neighborhood in
\(\Delta\) of the edge \(w\).  Recall that \(w\) is also a
reflection.  Thus \(w(\Delta_w)\) is a triangle in \(w(\Delta)\)
that is a neighborhood of \(w\) in \(w(\Delta)\) and \(\Delta_w\cap
w(\Delta_w)\) is a neighborhood of \(w\) in \(S\).

\FIGURE{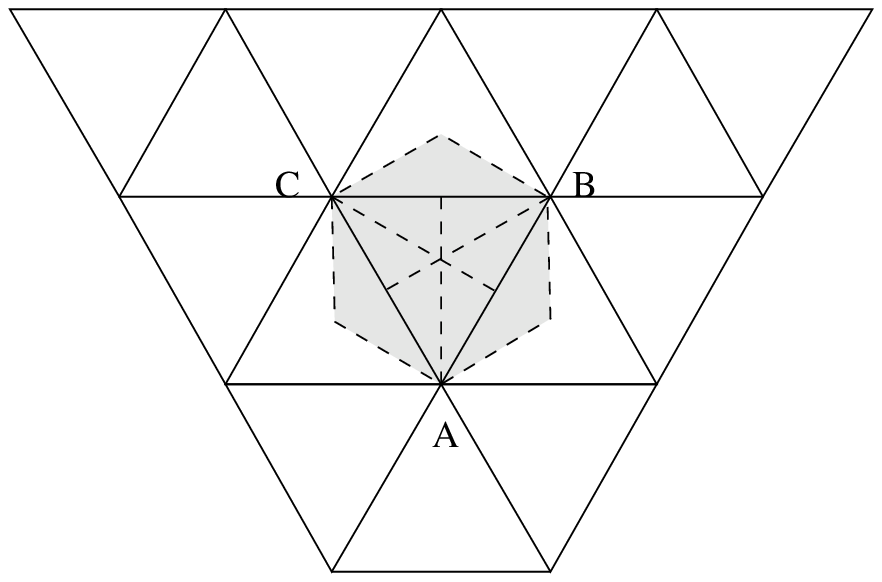}{The Hexagon \(\Xi\)}

We now let \(\Xi\) be the hexagon \[\Delta\cap x(\Delta_x)\cap
y(\Delta_y)\cap z(\Delta_z).\] The vertices \(A\), \(B\) and \(C\)
are now three of the six vertices of \(\Xi\).  See Figure
\ref{HexagonalTessellation}. The angle in \(\Xi\) at each of these
vertices is twice the angle found in \(\Delta\).  Thus the angle in
\(\Xi\) is \(2\pi/p\) at \(A\), \(2\pi/q\) at \(B\) and \(2\pi/r\)
at \(C\).

It is now possible to imitate the proof that \(\Delta\) tessellates
\(S\) under the action of \(G\) to show that \(\Xi\) tessellates
\(S\) under the action of \(G_1\).  One forms a surface from
\(\Gamma\x \Xi\) and shows that it is a cover of \(S\).  Slightly
more care is needed to define the identifications since the
directions of the rotations need to be taken into account.  The
important fact is that the action of \(G_1\) produces surface
neighborhoods of the edges of \(\Xi\) and the vertices \(A\), \(B\)
and \(C\).  We leave the details to the reader.  It also follows in
exactly the same way that \(\Gamma\) and \(G_1\) are isomorphic.  We
now drop the notatino \(G_1\) and use only \(\Gamma\) from now on.

We have now identified \(\Gamma(p,q,r)\) with a certain group of
symmetries of a triangular (actually hexagonal) tessellation of
\(S^2\), \(E^2\) or \(H^2\).  The group of symmetries is nice in
that it sets up a one to one correspondence between the elements of
\(\Gamma(p,q,r)\) and some of the triangles in the tessellation.
Namely, we can associate the element \(g\) in \(\Gamma\) with the
triangle \(g(\Delta)\) in \(S\).  The triangles used in this way are
``half'' the triangles in \(S\) (since the fundamental domain
\(\Xi\) for \(\Gamma\) is twice the size of the fundamental domain
\(\Delta\) for \(G\)) and are the triangles in \(S\) that are
obtained from \(\Delta\) by orientation preserving elements of
\(G\).  Thus \(\Gamma\) is the index two subgroup of \(G\)
consisting of all the orientation preserving elements of \(G\).

Note that the vertices of the triangular tessellation of \(S\) can
be given well defined labels from the alphabet \(\{ABC\}\).  This is
because each triangle in the tessellation in the image of \(\Delta\)
under a unique element of \(G\), and because when two triangles
\(g(\Delta)\) and \(g'(\Delta)\) in the tessellation share an edge
\(g(w)\), then \(g|w=g'|w\).  Thus for a vertex \(V\) in the
tessellation, we give it the same label as a vertex of \(\Delta\)
that is carried to \(V\) by any element of \(G\).

We can compute the orders of the triangle groups that are finite.
These are the groups \(\Gamma(2,2,r)\), \(\Gamma(2,3,3)\),
\(\Gamma(2,3,4)\) and \(\Gamma(2,3,5)\).  We will do an Euler
characteristic argument.  Let \(V\), \(E\) and \(F\) be the number
of vertices, edges and faces of the triangular tessellation.  The
vertices break into three classes, one label \(A\) and angles
\(\pi/p\), one with label \(B\) and angles \(\pi/q\) and one with
label \(C\) and angles \(\pi/r\).  We use the letters \(P\), \(Q\)
and \(R\) to represent the number of vertices in the respective
classes.  Thus \(V=P+Q+R\).  Each face is a triangle with exactly
one vertex of each class.  Thus \(2pP=2qQ=2rR=F\).  Also \(3F=2E\).
So
\[\chi(S^2)=2=V-E+f=\frac{F}{2}(\frac{1}{p}+\frac{1}{q}+\frac{1}{r}-1).\]
This allows us to get \(F\) for the various values of \((p,q,r)\).
The order of \(\Gamma\) is one half \(F\) so we get
\(|\Gamma(2,2,r)|=2r\), \(|\Gamma(2,3,3)|=12\),
\(|\Gamma(2,3,4)|=24\) and \(|\Gamma(2,3,5)|=60\).  The first is not
surprising since \(\Gamma(2,2,r)\) is the dihedral group of order
\(2r\).

The next paragraph is not needed.  I did not realize this when I put
it in, and I don't want to take it out.

We would like to argue that the finite groups are all different.
The only problem might be that \(\Gamma(2,3,3)\), \(\Gamma(2,3,4)\)
or \(\Gamma(2,3,5)\) might be dihedral groups.  Taking the
presentations as \(\langle a,b\mid a^2=b^q=(ab)^r=1\rangle\), we can
compute the abelianizations.  For \(\Gamma(2,2,r)\) we get \(2a=0\),
\(2b=0\) and \(ra+rb=0\).  If \(r\) is even, we get no new
information from \(ra+rb=0\) and the ableianization is
\(\ints_2\oplus\ints_2\).  If \(r\) is odd, then \(ra+rb=0\) reduces
to \(a+b=0\) so \(a=-b\) and the abelianization is \(\ints_2\).
Similar calculations show that \(\Gamma(2,3,3)\) abelianizes to
\(\ints_3\), \(\Gamma(2,3,4)\) abelianizes to \(\ints_2\) and
\(\Gamma(2,3,5)\) abelianizes to the trivial group.  The only
possibility for a dihedral group is \(\Gamma(2,3,4)\) but its order
is a multiple of 4 and the dihedral group of that order abelianizes
to \(\ints_2\oplus\ints_2\).  Thus the finite triangle groups are
all different.

We return to an arbitrary triangle group \(\Gamma(p,q,r)\).  Drawing
pictures shows that the rotations \(\lambda\) and \(\kappa\) do not
commute except when \(p=q=r=2\).  In this case there are four
triangles being permuted on \(S^2\) by a group of order four with
three involutions.  Thus \(\Gamma\) is not cyclic.

If \(\Gamma\) is a non-trivial free product, then either \(\Gamma\)
is \(\ints_2*\ints_2\) and \(\Gamma\) is two ended, or one of the
free factors has order at least 3 and \(\Gamma\) is infinitley
ended.  However, \(\Gamma\) acts as a group of covering
transformations on a simply connected space with compact quotient.
The simply connected space is either homeomorphic to \(S^2\) or
\(R^2\) and so \(\Gamma\) can have only zero ends or one end.  In no
case can \(\Gamma\) be a non-trivial free product.

\section{Properties of ``small'' of Seifert fibered spaces}

We start work on the excluded spaces in Lemma \ref{IrredClosedI}.
We will break the discussion into cases.

\subsection{Orbit surface \protect\(S^2\protect\) and less than four
exceptional fibers}

All Seifert fibered spaces with orbit surface \(S^2\) are orientable
and cannot contain two sided projective planes.  Thus showing
irreducibility is sufficient to show \(P^2\)-irreducibility.  Let
\(M\) be a closed Seifert fibered space with orbit surface \(S^2\)
and exactly three exceptional fibers.

\begin{lemma} The Heegaard genus of \(M\) is no more than 3.
\end{lemma}

\begin{demo}{Proof} Let \(A\) be the union of pairwise disjoint
fibered solid torus neighborhoods of the exceptional fibers, and let
\(B\) be the closure of the complement of \(A\).  We know that the
classifying homomorphism of \(M\) is trivial, so \(B\) is the
product of a disk with two holes with \(S^1\).  That is \(B\) is the
double along an annulus of a product of an annulus with \(S^1\).  A
hole can be drilled in an annulus cross \(S^1\) to connect the two
boundary components in such a way that a handlebody with of genus
two results.  The hole can miss a given annulus in the boundary.
The annulus will then become a \(\pi_1\) generator for one of the
handles.  Thus we can drill two holes in \(B\) to give the double of
two handlebodies of genus 2 doubled along this annulus.  The result
is a handlebody of genus 3.  The holes drilled in \(B\) are added to
\(A\) as 1-handles in such a way as to connect the boundary
components of \(A\).  This connects \(A\) and gives another
handlebody of genus 3.  \end{demo}

We are interested in whether \(M\) is irreducible.  Assume that it
is not.  Then it is a connected sum of two 3-manifolds, one of which
is genus 1, and \(\pi_1(M)\) is a non trivial free product, or
cyclic.  (Actually the only simply connected 3-manifold of Heegaard
genus 2 is \(S^3\), but this is hard to prove and we don't need it.)

A presentation for \(\pi_1(M)\) is \begin{align*}\langle
h,c_1,c_2,c_3\mid & c_1hc_1^{-1}=c_2hc_2^{-1}=c_3hc_3^{-1}=h, \\&
c_1^{\mu_1}h^{\beta_1}=c_2^{\mu_2}h^{\beta_2}=c_3^{\mu_3}h^{\beta_3}=1,
\\& c_1c_2c_3h^b=1\rangle.\end{align*} Each \(\mu_i\) is at least 2.
Since \((\mu_i,\beta_i)=1\), each \(\beta_i\) is not zero.  If we
mod out by the normal cyclic subgroup \(\langle h\rangle\), then we
obtain a triangle group \(\Gamma(\mu_1,\mu_2,\mu_3)\).  This is not
cyclic, so \(\pi_1(M)\) is not cyclic.  If \(\langle h\rangle\) is
trivial, then \(\pi_1(M)\) is isomorphic to
\(\Gamma(\mu_1,\mu_2,\mu_3)\).  But \(\Gamma(\mu_1,\mu_2,\mu_3)\) is
neither cyclic nor a non-trivial free product.  Thus \(\langle
h\rangle\) is non-trivial, \(\pi_1(M)\) is not cyclic, and
\(\pi_1(M)\) must be a non-trivial free product.

Since the only non-trivial free product with a non-trivial cyclic,
normal subgroup is \(\ints_2*\ints_2\), we have
\(\pi_1(M)\cong\ints_2*\ints_2\).  We know that \(\langle h\rangle\)
is contained in the unique maximal cyclic, normal subgroup of
\(\ints_2*\ints_2\cong D_{\infty}\) and that \(h\) has infinite
order.  Each element not in the maximal normal, cyclic subgroup of
\(\ints_2*\ints_2\) is an involution.  However, each generator
\(c_i\) has a power equal to a non-zero power of \(h\) and so each
\(c_i\) has infinite order and is not an involution.  Thus all of
the generators are in the maximal cyclic, normal subgroup and there
can be no involutions.  Since an irreducible, orientable 3-manifold
is \(P^2\)-irreducible, we have shown:

\begin{lemma} A closed Seifert fibered space with orbit surface
\(S^2\) and exactly three exceptional fibers is \(P^2\)-irreducible.
\BlackBox\end{lemma}

Note that we know which of these spaces have infinite fundamental
group.

\begin{lemma} Let \(M\) be a closed Seifert fibered space with orbit
surface \(S^2\) and exactly three exceptional fibers.  Let
\(h\in\pi_1(M)\) be represented by an ordinary fiber.  If the
indexes of the exceptional fibers are not one of \((2,2,r)\),
\((2,3,3)\), \((2,3,4)\) or \((2,3,5)\), then \(\pi_1(M)\),
\(\langle h\rangle\) and \(\pi_1(M)/\langle h\rangle\) are all
infinite.  Otherwise all three groups are finite.  \end{lemma}

\begin{demo}{Proof} We know that \(M\) is \(P^2\)-irreducible and
boundary irreducible.  The result follows immediately from Lemma
\ref{InfinitePiOne} and the fact that the only finite triangle
groups are the ones specified by the triples listed in the
hypothesis.  \end{demo}

We can also can give a criterion for having an incompressible
surface.

\begin{lemma} Let \(M\) be a closed Seifert fibered space with orbit
surface \(S^2\) and exactly three exceptional fibers.  Then \(M\)
has an incompressible surface if and only if \(H_1(M)\) is infinite.
\end{lemma}

\begin{demo}{Proof} \(M\) has a non-separating incompressible
surface if and only if \(H_1(M)\) is infinite.  Thus we must show
that \(M\) has no separating incompressible surface when \(H_1(M)\)
is finite.  What is actually true is that no closed Seifert fibered
space with orbit surface \(S^2\) and exactly three exceptional
fibers has a separating incompressible surface.

Assume that \(M\) has a separating incompressible surface.  Let
\(N_1\), \(N_2\) and \(N_3\) be pairwise disjoint fibered solid
torus neighborhoods of the exceptional fibers.  Let \(S\) be a
separating incompressible surface in \(M\) whose intersections with
the \(N_i\) are meridinal disks of the \(N_i\) and so that the total
number of these disks is minimal.  Let \(M'\) be \(M\) with the
interiors of the \(N_i\) removed and let \(S'\) be \(S\cap M'\).  If
\(S'\) compresses in \(M'\) along a non-separating curve in \(S'\),
then \(S\) is not incompressible in \(M\).  If \(S'\) compresses in
\(M'\) along a separating curve \(J\) in \(S'\), then \(J\) bounds a
disk in \(S\) which can be replaced by the compressing disk to lower
the number of boundary components of \(S'\).  Thus the minimality of
the number of boundary components of \(S'\) implies the
incompressibility of \(S'\) in\(M'\).

The space \(M'\) is the product of \(S^1\) with a disk with two
holes \(E\) and \(\pi_1(M')\) has non-trivial center.  Also,
\(\pi_1(M')\) is a free product amalgamated over \(\pi_1(S')\) so
either \(\pi_1(S')\) is abelian, or \(S'\) carries the fundamental
group of one of its complementary domains.  In the latter case,
\(S'\) would be boundary parallel and would be an annulus or a
torus.  In the former case, \(S'\) must be a sphere, annulus or a
torus.  Since \(M\) and \(M'\) are irreducible, a sphere is ruled
out.  If \(S'\) is an annulus, then \(S\) is a sphere and this case
is also ruled out.

If \(S'\) is a separating incompressible torus in \(M'\), it would
hit \(E\) in simple closed curves.  We can eliminate the curves that
are trivial on either \(E\) or \(S'\).  Elements of \(\pi_1(E)\)
that are conjugate in \(\pi_1(M')\) are conjugate in \(\pi_1(E)\),
so the curves must all be parallel.  Since \(E\) is planar, the
curves all separate \(E\).  Since \(E\) has only three boundary
components, the curves must be parallel to one of them.  However,
these curves demonstrate that \(S'=S\) is compressible in \(M\).
Thus \(S'\) is disjoint from \(E\).  This puts \(S'\) in a product
of a disk with three holes with an interval and puts
\(\ints\oplus\ints\) as a subgroup of the free group on two
generators.  This gives a contradiction.  \end{demo}

\begin{remark} From the presentation for \(\pi_1(M)\), we get that
\(H_1(M)\) is an abelian group with four generators with relation
matrix \[\left(\begin{matrix} \beta_1 & \mu_1 & 0 & 0 \\ \beta_2 & 0 &
\mu_2 & 0 \\ \beta_3 & 0 & 0 & \mu_3 \\ b & 1 & 1 & 1
\end{matrix}\right).\] If this matrix is of full rank over \(\ints\),
then \(H_1(M)\) is finite.  It is easy create triples
\((\mu_1,\mu_2,\mu_3)\) that make \(\pi_1(M)\) infinite and chose
the \(\beta_i\) and \(b\) so that \(H_1(M)\) is finite.  These were
the first examples discovered of closed, \(P^2\)-irreducible
3-manifolds with infinite \(\pi_1\) and no incompressible surfaces.
There are now many more known.  \end{remark}

We now consider a Seifert fibered space \(M\) with orbit surface
\(S^2\) and with fewer than three exceptional fibers.

\begin{lemma} The Heegaard genus of \(M\) is no more than 1.
\end{lemma}

\begin{demo}{Proof} There is a circle in the orbit surface whose
complementary domains are disks containing no more than one
exceptional point.  The preimages of these disks are fibered solid
tori in \(M\) and represent \(M\) as a closed 3-manifold with
Heegaard genus no more than 1.  \end{demo}

Thus \(M\) is a lens space.  If \(M\) is not irreducible, then
Haken's theorem gives us a sphere that intersects each solid torus
in the Heegaard splitting in a disk.  This shows that \(M\) is
\(S^2\x S^1\).  We know that \(S^2\x S^1\) is in fact realizable as
a Seifert fibered space with orbit surface \(S^2\) and with zero or
two exceptional fibers.  We have shown:

\begin{lemma} The only closed Seifert fibered space that is not
\(P^2\)-irreducible and that has an orientable orbit surface is
\(S^2\x S^1\).  \BlackBox\end{lemma}

We can also analyze the other properties that we have been
considering.

\begin{lemma} Let \(M\) be a closed, Seifert fiber space with orbit
surface \(S^2\) and with no more than two exceptional fibers.  Then
\(\pi_1(M)\) is cyclic, generated by some fiber (possibly
exceptional) and is finite except when \(M\) is \(S^2\x S^1\).  In
all cases, the associated Fuchsian group is finite.  The only case
in which \(M\) has an incompressible surface is when \(M\) is
\(S^2\x S^1\) in which case \(M\) has a 2-sphere that bounds no
3-cell.  \end{lemma}

\begin{demo}{Proof} This all follows from the structure of lens
spaces.  The generator of the fundamental group of a lens space is a
centerline of one of the solid tori in a Heegaard splitting.  An
ordinary fiber is a power of this centerline.  An incompressible
surface in a closed, orientable 3-manifold with finite fundamental
group must be a separating sphere.  The only reducible lens space is
\(S^2\x S^1\). \end{demo}

\subsection{Orbit surface \protect\(P^2\protect\) and less than two
exceptional fibers}

If \(M\) has orbit surface \(P^2\), then there are only two possible
classifying homomorphisms.  The trivial homomorphism yields a
non-orientable \(M\), and the non-trivial homomorphism yields an
orientable \(M\).  As in the case of orbit surface \(S^2\), we only
have to show that \(M\) is irreducible when \(M\) is orientable.
This will be easy because all of the orientable spaces that we will
be presented with except one will turn out to be spaces that we have
analyzed before.  The remaining space will turn out to be reducible.
We will also discover that there are only two non-orientable spaces
that arise.  These will turn out not to be \(P^2\)-irreducible.  We
will consider the cases of orientable \(M\) and non-orientable \(M\)
separately.

Let \(M\) be a Seifert fibered space with orbit surface \(P^2\).  If
\(M\) has no more than 1 exceptional fiber, then there is a disk on
\(P^2\) whose preimage in \(M\) is a fibered solid torus \(T\), and
whose closed complement on \(P^2\) is a \Mob\ \(A\) that has no
exceptional points.  The preimage \(B\) of \(A\) in \(M\) is a
circle bundle over \(A\).  The preimage \(S\) of the centerline of
\(A\) is either a torus or a Klein bottle depending on the
classifying homomorphism of \(M\).  If \(S\) is a torus, the
\(S^1\)-bundle over \(A\) is not twisted and \(M\) and \(B\) are
non-orientable.  If \(S\) is a Klein bottle, then the \(S^1\)-bundle
over \(A\) is twisted and \(M\) and \(B\) are orientable.

We can view \(B\) as an \(I\)-bundle over \(S\) and \(A\) is a
subbundle over a circle in \(S\).  Since \(A\) is a \Mob, the
\(I\)-bundle is twisted.  Thus the corresponding \(I\)-bundle over
\(S\) must be twisted.  There is only one class of classifying
homomorphisms over a torus that involve twisting.  There are two
that involve twisting over a Klein bottle, but only one will apply
here because when \(S\) is a Klein bottle, \(B\) must be orientable.

Assume that \(S\) is a torus.  Then there is no twisting along the
centerline of \(A\) and \(B\) is just the product of \(A\) with a
circle.  This is the case where \(M\) is non-orientable and the
classifying homomorphism for \(M\) is trivial.  However we
concentrate on the fact that \(B\) is an \(I\)-bundle over \(S\).
We let \(\phi\) be the classifying homomorphism for \(B\).  We will
have no reason to refer to the classifying homomorphism of the
Seifert fibered space \(M\) during the rest of this case.  Since
there is only one class of classifying homomorphism over a torus
that involves twisting, we can choose one representative for our
convenience.  Let \((J,K)\) be a pair of generating curves for
\(H_1(S)\) that intersect in a point and assume that \(\phi(J)=+1\)
and \(\phi(K)=-1\).  Any self homeomorphism of \(S\) that preserves
\(\phi\) extends to a fiber preserving, self homeomorphism of \(B\).
(In fact there are two extensions.  There is a homeomorphism of
\(B\) of period two that fixes \(S\) and reverses each of the
\(I\)-fibers.  This can be composed with any extension.)  Since the
value of \(\phi\) on \(aJ+bK\) depends only on the parity of \(b\),
we know that a homeomorphism whose action on \(H_1(S)\) is given by
\(\mtbt{a}{c}{b}{d}\) must have \(b\) even and \(d\) odd.  Since
\(\dtbt{a}{c}{b}{d}=\pm1\) it suffices to require that \(b\) be
even.

The boundary of \(B\) is a torus \(S'\) that double covers \(S\) by
pulling in along the fibers.  The fibers over \(J\) form an annulus.
One boundary of this annulus is a curve \(J'\) on \(S\).  The fibers
over \(K\) form a \Mob.  The boundary of this \Mob\ is a curve
\(K'\) that represents \(2K\) in the homology of \(B\).  The curves
\(J'\) and \(K'\) intersect in a point and give a generating pair
for \(H_1(S')\).

The action of the homeomorphism of \(B\) that reverses all the
\(I\)-fibers preserves the classes of \(J'\) and \(K'\).  Thus an
allowable homeomorphism of \(S\) induces a well defined automorphism
on \(H_1(S')\) by extending it to one of the two fiber preserving,
self homeomorphisms of \(B\).  We identify the homeomorphisms of
\(S\) with the matrices that give the effect on \(H_1(S)\).

Let a homeomorphism of \(S\) that preserves \(\phi\) be given by
\(\mtbt{a}{c}{2b}{d}\).  This homeomorphism takes \(J\) and the
annulus over it to \(\cv{a}{2b}\) and the annulus over it.  One
boundary component of this annulus runs \(a\) times around \(J\) and
\(2b\) times around \(K\).  This has it running \(a\) times in the
\(J'\) direction on \(S'\) and only \(b\) times around the \(K'\)
direction on \(S'\).  This can also be done with intersection
numbers since only half of the intersections of the curve with the
annulus over \(J\) are intersections with \(J'\) while all of the
intersections with the \Mob\ over \(K\) are intersections with
\(K'\).  Thus the image of \(J'\) is \(aJ'+bK'\) in \(H_1(S')\).

The same homeomorphism takes \(K\) and the \Mob\ over it to
\(\cv{c}{d}\) and the \Mob\ over it.  The boundary of this \Mob\ on
\(S'\) runs twice over the curve \(\cv{c}{d}\) on \(S\) and an
analysis similar to the image of \(J'\) shows that this curve is
\(2cJ'+dK'\) in \(H_1(S')\).  Thus the automorphism of \(H_1(S')\)
is given by \(\mtbt{a}{2c}{b}{d}\).  Any matrix
\(\mtbt{a}{2c}{b}{d}\) acting on \(H_1(S')\) is realized by an
extension of the homeomorphism of \(S\) that realizes
\(\mtbt{a}{c}{2b}{d}\) on \(H_1(S)\).

Let \(mJ'+nK'\) be a simple closed curve on \(S'\).  We know
\((m,n)=1\).  If \(m\) is odd, then its prime factorization does not
use the prime 2 and \((m,2n)=1\).  There are integers \(a\) and
\(c\) so that \(am+2cn=1\).  Thus \(\dtbt{a}{2c}{-n}{m}=1\) and
\(\mtbt{a}{2c}{-n}{m}\cv{m}{n}=\cv{1}{0}\).  If \(m\) is even, then
we just use the fact that \((m,n)=1\) to get integers \(b\) and
\(d\) with \(bm+dn=1\).  Then \(\dtbt{n}{-m}{b}{d}=1\) and
\(\mtbt{n}{-m}{b}{d}\cv{m}{n}=\cv{0}{1}\).  Both matrices are
realizable by extensions of homeomorphisms of \(S\).  Thus up to
homeomorphisms of \(B\), there are only two simple closed curves on
\(\bd B\) to which to attach a meridian of a fibered solid torus.

From the last paragraph, there are at most two topological spaces
that result from sewing a solid torus to \(B\).  There are many more
than two Seifert fibered spaces that result from sewing a fibered
solid torus to the trivial circle bundle over the \Mob.

We see what we get when we sew a solid torus with meridian sewn to
the curves from either of the two classes that we have found.  Note
that sewing a meridian along along the curve \(\cv{1}{0}\) does not
give a Seifert fibered space since the curve \(\cv{1}{0}\) is
homologous to a fiber.  We can use the curve \(\cv{1}{1}\) instead
which is in the same class when we wish to look at the space as a
Seifert fibered space.  The other class is represented by
\(\cv{0}{1}\).

Sewing a meridian of a solid torus to \(\cv{1}{0}\) sews a disk to
each of the boundary components of the annulus in \(B\) containing
\(J\).  This creates a non-separating \(S^2\) in \(M\).  Cutting
along this \(S^2\) cuts \(B\) into \(S^1\x I\x I\) from which \(B\)
is recovered by sewing \(S^1\x I\x\{0\}\) to \(S^1\x I\x\{1\}\) with
a homeomorphism that takes each \((x,t,0)\) to \((x,1-t,1)\).  The
cut along \(S^2\) cuts the solid torus into two cylinders that are
sewn to \(S^1\x I\x I\) along \(S^1\x\{0,1\}\x I\).  It is seen that
\(M\) is the non-orientable \(S^2\) bundle over \(S^1\).

To view \(M\) as a Seifert fibered space, we can sew the meridian of
the solid torus along \(\cv{1}{1}\) instead.  This is a crossing
curve for the fibers of \(M\) in \(\bd B\) so the resulting space
has no exceptional fibers.  The curve has intersection number one
with the boundary \(K'\) of the section of the projection from \(B\)
to its orbit surface \(A\) viewing \(B\) as a Seifert fibered space.
Thus the obstruction to a section for \(M\) is 1.  Since any curve
\(\cv{m}{n}\) can be used with \(m\) odd to get the same topological
space and \(n\) is the index of the fiber that results, we can
realize \(M\) with one exceptional fiber of any given index.  Note
that the presentation for \(\pi_1(M)\) gotten by using the curve
\(\cv{1}{1}\) is \(\langle h,x\mid xhx^{-1}h^{-1}=x^2h^1=1\rangle\)
or \(\langle h,x\mid x^2=h^{-1}\rangle\) which is just \(\ints\).

The other space that can be obtained is gotten by sewing a meridian
to \(\cv{0}{1}\).  This is a crossing curve for the fibers of \(M\)
in \(\bd B\), so this realizes \(M\) as a space with no exceptional
fibers.  It attaches a disk to the boundary of the section for \(B\)
so that it gives a section to all of \(M\).  Since there are no
exceptional fibers and \(B\) is \(A\x S^1\), this realizes \(M\) as
\(P^2\x S^1\).  A check of \(\pi_1(M)\) gives \(\langle h,x\mid
xhx^{-1}h^{-1}=x^2h^0=1\rangle\) or \(\ints\oplus\ints_2\) which is
consistent with our description of \(M\).  Since \(\pi_1(M)\) is not
a non-trivial free product, \(M\) is irreducible.  It is not
\(P^2\)-irreducible since it has a two sided copy of \(P^2\).  We
can also realize \(P^2\x S^1\) as a Seifert fibered space with a
single exceptional fiber of any index.  We have shown:

\begin{lemma} The only non-orientable toplogical spaces that can be
obtained as Seifert fibered spaces over \(P^2\) with less than two
exceptional fibers are \(P^2\x S^1\) and the twisted \(S^2\) bundle
over \(S^1\).  Niether space is aspherical as the first space is
irreducible but not \(P^2\) irreducible, and the second space is not
irreducible.  Both spaces can be realized with no exceptional fibers
or with one exceptional fiber of any integer index.
\BlackBox\end{lemma}

Note that both of these spaces have infinite fundamental group.  One
has fundamental group \(\ints\oplus\ints_2\) and the other has
fundamental group \(\ints\).  To see what the fate of an ordinary
fiber is in these spaces, we have to calculate the fundamental
groups using the curves \(J\), \(K\), \(J'\), \(K'\) described
above.  We can take \(J\) to be the fiber of \(B\) seen as an
\(S^1\)-bundle over the \Mob\ \(A\).  This is the fiber structure
inherited from \(M\).  The meridian of the fibered solid torus will
be attached along the simple closed curve \(C=mJ'+nK'\) on \(S'\)
with \((m,n)=1\).  Since \(J\) and \(K\) generate
\(H_1(B)=\ints\oplus\ints\), we want to express the curve in terms
of \(J\) and \(K\) and get \(C=mJ+2nK\).  If \(m\) is odd, then
\((m,2n)=1\) and the resulting fundamental group is \(\ints\).  If
\(m\) is even, then \((m,2n)=2\) and the resulting fundamental group
is \(\ints\oplus\ints_2\).  The only way for the fiber \(J\) to be
torsion in the result is for \(n\) to be zero.  This is not an
allowable sewing however, since this would put a meridian of the
attached fibered solid torus homologous to a fiber.  In all of the
above, modding out by the subgroup generated by the fiber yields the
cyclic group \(\ints_{2n}\).  We can summarize.

\begin{lemma} A non-orientable Seifert fibered space with orbit
surface \(P^2\) and with less than two exceptional fibers has
fundamental group \(\ints\) or \(\ints\oplus\ints_2\).  In all cases
the subgroup generated by an ordinary fiber is infinite cyclic and
the corresponding Fuchsian quotient is finite cyclic of non-zero
even order.  \BlackBox\end{lemma}

We now assume that \(S\) is a Klein bottle.  This means that the
classifying homomorphism of \(M\) is the only non-trivial
homomorphism available with \(P^2\).  Thus \(M\) and \(B\) are
orientable and \(B\) is the orientable \(I\)-bundle over the Klein
bottle.  This fibers as a Seifert fibered space in two ways.  One is
the fibration inherited from \(M\) which has the \Mob\ \(A\) as
orbit space.  The other uses curves parallel to the orientation
reversing curves of \(S\) that has a disk as orbit space with two
exceptional fibers of index two.  In spite of the fact that this is
not the fibering that we are presented with we will use it.  One
reason is that we are trying to extract topological information and
it does not matter which fibering we use.  The other is that it
leads to structures that we have already analyzed.

We must be careful since adding a fibered solid torus to \(B\) does
not allow a meridian of the of the fibered solid torus to be
homologous to a fiber of \(B\).  Thus the fiber of the new fibering
must be treated separately since it is not a valid sewing curve in
the new fibering but is allowed in the original.  Also, the fiber of
the original fibering must not be used even though it is allowed in
the new fibering.

Adding a fibered solid torus to the new fibering gives orbit surface
\(S^2\) with two or three exceptional fibers.  These spaces are
already known to be irreducible.  They have Fuchsian quotient a
triangle group \(\Gamma(2,2,r)\) if there is a third exceptional
fiber, and are lens spaces if there are only two exceptional fibers.

If we add a fibered solid torus with meridian sewn to the fiber of
the new fibering, then the meridian is a crossing curve of the
original fibering.  In fact the meridinal disk completes a section
of the orbit surface of \(M\).  Thus \(M\) is the orientable \(S^1\)
bundle over \(P^2\) which we have previously identified as
\(P^3\#P^3\).  This has 0 obstruction to a section and its
fundamental group has the standard presentation \(\langle x,h\mid
xhx^{-1}=h^{-1}, x^2=1\rangle\) of \(D_{\infty}\) or
\(\ints_2*\ints_2\).  The cyclic subgroup generated by an ordinary
fiber is infinite and the corresponding Fuchsian quotient is
\(\ints_2\oplus\ints_2\).

Note that only one sewing yields \(P^3\#P^3\).  When we note that
\(P^3\#P^3\) does not come up in any of the other sewing of Seifert
fibered spaces that are not \(P^2\)-irreducible, then we will know
that the fiber structure of \(P^3\#P^3\) is unique.

The forbidden sewing of the original fibering would add a solid
torus along a crossing curve in the new fibering.  Thus there would
be only two exceptional fibers and the space would be a lens space.
The curve would be one boundary component of an annulus which is an
\(I\)-bundle over an orientation preserving curve of the Klein
bottle \(S\).  Thus a meridinal disk of the added solid torus and a
parallel copy of this disk in the solid torus would create a sphere
with the annulus in \(B\).  Since the annulus does not separate
\(B\), the sphere would not separate the lens space.  Thus the
forbidden lens space is \(S^2\x S^1\).

We can show that \(S^2\x S^1\) does not arise from any other sewing.
The sewing must have no third exceptional fiber, so a crossing curve
must be used as the image of a meridian.  The sewing is then
determined by the invariant \(b\).  The fundamental group of a
sewing with no third exceptional fiber is presented as
\begin{align*}\langle h,c_1,c_2\mid & c_1hc_1^{-1}=c_2hc_2^{-1}=h,
\\& c_1^{\mu_1}h^{\beta_1}=c_2^{\mu_2}h^{\beta_2}=1, \\&
c_1c_2h^b=1\rangle.\end{align*} We have \(c_2=c_1^{-1}h^{-b}\) and
the fact that \(c_2\) commutes with \(h\) follows from this and the
fact that \(c_1\) commutes with \(h\).  We also know that the
indexes of the two exceptional fibers are both two, so
\(\mu_1=\mu_2=2\) and \(\beta_1=\beta_2=1\).  Thus we are left with
an abelian group generated by \(c_1\) and \(h\) with relations
\(2c_1+h=0\) and \(2c_1+(2b-1)h=0\).  The only way for the group to
be \(\ints\) is for the relations to be dependent.  This only
happens when \(b=1\).  Thus there is only one sewing that yields
\(S^2\x S^1\) and it is not one that agrees with a fibration over
\(P^2\) with no more than one exceptional fiber.

We can determine which lens spaces are obtainable.  First we see
which are possible, and then we show that they are realized.  Two
solid tori are sewn together by a homeomorphism of their boundaries.
We assume that the solid tori are oriented and have generators for
their first homologies so that the fibers are given as vectors
\(\cv{1}{2}\).  The action of the homeomorphism is given by
\(\mtbt{q}{r}{p}{s}\).  Since we want an orientable result, the
homeomorphism should reverse the orientations of the boundaries.
Thus \(\dtbt{q}{r}{p}{s}=-1\).  (This is not essential since fibered
solid tori determines by \(1/2\) admit afiber preserving,
orientation reversing self homeomorphism.)  We have
\[\cv{1}{2}=\mtbt{q}{r}{p}{s}\cv{1}{2}=\cv{q+2r}{p+2s}.\]
Multiplying both sides of \(q+2r=1\) by \(s\) and both sides of
\(p+2s=2\) by \(r\) and subtracting gives \(qs-pr=s-2r\) which
equals \(-1\) because of the restriction on the determinant.
Solving for \(s\) leads to \(p=4-4r\) and \(q=1-2r\).  Letting
\(n=1-r\) gives \(p=4n\) and \(q=2n-1\).  Thus the possible lens
spaces are the spaces \(L_{4n,2n-1}\).  We will show that they are
all realized by showing that all the possible fundamental groups are
realized.  This will suffice since each fundamental group is
represented by exactly one space among the possibilities.  We know
that the fundamental groups realized are abelian, generated by
\(c_1\) and \(h\) and have relations \(2c_1+h=0\) and
\(2c_1+(2b-1)h=0\).  So \(h=-2c_1\) and we get \(4c_1(1-b)=0\).
This realizes all cyclic groups of order a multiple of four and all
the \(L_{4n,2n-1}\) can be achieved.

\begin{lemma} If \(M\) is an orientable, Seifert fibered space with
orbit surface \(P^2\) and less than two exceptional fibers, then
\(M\) is homeomorphic either to a lens space of type
\(L_{4n,2n-1}\), or to a Seifert fibered space with orbit surface
\(S^2\) and three exceptional fibers with two of index 2, or to the
connected sum of two copies of \(P^3\).  All of the fundamental
groups are finite except for the fundamental group of \(P^3\#P^3\).
An ordinary fiber has infinite order in \(P^3\#P^3\) and the
corresponding Fuchsian quotient is \(\ints_2*\ints_2\). \end{lemma}

\subsection{\label{SmallSFS}Summary of ``small'' Seifert fibered
spaces}

We identify certain Seifert fibered spaces as ``small.''  The list
that follows will mention these spaces in overlapping groups.  The
``small'' Seifert fibered spaces are: fibered solid tori; lens
spaces which include proper lens spaces (all except \(S^3\) and
\(S^2\x S^1\)), the irreducible lens spaces (all except \(S^2\x
S^1\)), and the non-simply connected lens spaces (all except
\(S^3\)); the two \(S^2\)-bundles over \(S^1\); \(P^2\x S^1\);
\(P^3\#P^3\); and the platonic Seifert fibered spaces (Seifert
fibered spaces with orbit surface \(S^2\) and exactly three
exceptional fibers of indexes \((2,2,r)\), \((2,3,3)\), \((2,3,4)\)
or \((2,3,5)\)).

We have proven the following.

\begin{thm}\label{NotSmall} A compact, connected Seifert fibered
space \(M\) that is not ``small'' is \(P^2\)-irreducible, boundary
irreducible, aspherical, has torsion free \(\pi_1\), has an
incompressible surface with infinite fundamental group, and
\(\pi_1(M)\), \(\langle h\rangle\) and \(\pi_1(M)/\langle h\rangle\)
are all infinite where \(h\) is represented by some ordinary fiber.
\BlackBox\end{thm}

The exceptions to the various properties are listed below.

Not boundary irreducible: fibered solid tori.

Not \(P^2\)-irreducible: Both \(S^2\)-bundles over \(S^1\); \(P^2\x
S^1\); \(P^3\#P^3\).

Not aspherical: All ``small'' Seifert fibered spaces except the
fibered solid tori.

Torsion in \(\pi_1\): Platonic Seifert fibered spaces; irreducible
lens spaces; \(P^2\x S^1\); \(P^3\#P^3\).

Finite \(\pi_1\): Platonic Seifert fibered spaces; irreducible lens
spaces.

Finite \(\langle h\rangle\): Same as finite \(\pi_1\).

Finite Fuchsian quotient: All ``small'' Seifert fibered spaces.

No incompressible surfaces other than an \(S^2\) bounding no 3-cell
or a properly embedded non-boundary parallel disk: All ``small''
Seifert fibered spaces.

Thus an adequate definition of a ``small'' Seifert fibered space
could be that the subgroup of \(\pi_1\) generated by an ordinary
fiber is of finite index in \(\pi_1\).

Note that this consideration of cases has led to a result.

\begin{lemma} Let \(M\) be a compact, connected Seifert fibered
space and let \(\langle h\rangle\) be a cyclic subgroup of
\(\pi_1(M)\) generated by an ordinary fiber.  Then \(\langle
h\rangle\) is infinite if and only if \(\pi_1(M)\) is infinite.
\BlackBox\end{lemma}

\section{Relating the fundamental group to the topology}

In this section we try to discover how much of the fiber structure
of a Seifert fibered space \(M\) is determined by the topology of
the space.  The main problem is to try to identify a fiber from the
topology of \(M\).  The first step in this is to identify the
element of \(\pi_1(M)\) represented by a fiber.  To do this we have
to identify the behavior of an element of \(\pi_1(M)\) that is
represented by a fiber that is different from behavior of other
elements of \(\pi_1(M)\).

We know that an element \(h\) represented by an ordinary fiber
generates a cyclic normal subgroup.  From this it follows that the
centralizer of \(h\) has index no more than two in \(\pi_1(M)\).
Thus \(h\) commutes with at least half the elements in \(\pi_1(M)\).
However, \(\pi_1(M)\) is far from abelian for most Seifert fibered
spaces.  The quotient \(\pi_1(M)/\langle h\rangle\) is a Fuchsian
group which has as a quotient the fundamental group of a surface.
For most surfaces, centralizers of elements of their fundamental
groups are small.  This suggests that the cyclic normal subgroup
\(\langle h\rangle\) might be distinguished in \(\pi_1(M)\) by its
algebraic properties.

\subsection{Fuchsian groups}

Let \(H\) be a Fuchsian group.  It has presentation given by either
\ref{FuchsOrbl} or \ref{FuchsNOrbl}.  A complex \(X\) that has
\(\pi_1(X)\cong H\) can be obtained by modifying a surface.

If the presentation of \(H\) is given by \ref{FuchsOrbl}, then let
\(Y\) be an orientable surface of genus \(g\) and \(n+m\) boundary
components.  Let \(c_1,\dots,c_n\) denote the first \(n\) boundary
components, let \(E_1,\dots,E_n\) denote pairwise disjoint disks
that are disjoint from \(Y\), let \(h_i:\bd E_i\into c_i\) be a map
of degree \(\mu_i\), and let \(X\) be formed by attaching the
\(E_i\) to \(Y\) using the maps \(h_i\).  The standard calculation
of a presentation for \(\pi_1(X)\) gives the presentation
\ref{FuchsOrbl}.  The unused boundary components of \(Y\) correspond
to the elements \(d_i\) in the presentation.

If the presentation of \(X\) is given by \ref{FuchsNOrbl}, then the
same construction is done by starting with a non-orientable surface.

The surface \(Y\) in the above construction is compact and
connected.  We can enlarge the notion of a Fuchsian group to include
the fundamental group of any complex constructed as above, requiring
only that \(Y\) be a connected surface.

It pays to loosen the requirements even more by allowing any number
of disks to be attached to a given component of the boundary of the
surface.  We argue that this does not enlarge the class of groups
considered.  Let \(J\) be a boundary component of the surface to
which several disks are attached and let \(\alpha\) be the element
of the fundamental group that \(J\) represents.  Each disk attached
along \(J\) guarantees that some power \(d_i\) of \(\alpha\) is the
identity.  Let \(d\) be the greatest common divisor of the \(d_i\).
(This makes sense even if the number of attached disks is infinite.)
Then \(\alpha^d=1\) is a consequence of the relations
\(\alpha^{d_i}=1\), and each \(\alpha^{d_i}=1\) is a consequence of
the relation \(\alpha^d=1\).  Thus removing all the disks attached
to \(J\) and replacing them with one disk attached by a map of
degree \(d\), yields a space with the same fundamental group.  If
this is done for every boundary component of \(Y\), then we obtain a
Fuchsian complex with the same fundamental group as the original
group and with no more than one disk attached to each boundary
component of \(Y\).

We thus define a \defit{Fuchsian complex} to be the space obtained
from a connected surface by attaching disks to the surface using
maps of non-zero degree from the boundaries of the disks to the
boundary components of the surface.  (Of course the requirement that
the maps be of non-zero degree forces the disks to be attached to
compact boundary components of the surface.)  A \defit{Fuchsian
group} is the fundamental group of a Fuchsian complex.  This
definition makes the following immediate.

\begin{lemma} A cover of a Fuchsian complex is a Fuchsian complex,
and a subgroup of a Fuchsian group is a Fuchsian group.
\BlackBox\end{lemma}

We have shown that a Fuchsian group can be represented by a Fuchsian
complex \(X\) based on a surface \(Y\) in which each component of
\(\bd Y\) has no more than one disk attached to it.  Further, if
there is a disk attached by a map of degree one, then it can be
attached by a homeomorphism.  In this case, the disk can be regarded
as part of \(Y\) and the attaching curve removed from the set of
boundary components of \(Y\).  Thus every Fuchsian group can be
represented by a Fuchsian complex in which no more than one disk is
attached to any boundary component of the base surface and in which
every attaching map has degree at least two.  Such a complex is said
to be \defit{simplified}.

If \(X\) is a Fuchsian complex based on a connected surface \(Y\),
then we define the boundary of \(X\) to be all of the boundary of
\(Y\) to which no disks are attached.  Note that this is not an
invariant of the group since we can remove all of the boundary of
\(X\) and have the same fundamental group.

\begin{lemma} Let \(X\) be a finite Fuchsian complex with non-empty
boundary based on a surface \(Y\)..  Then \(\pi_1(X)\) is a free
product with one factor a (possibly trivial) free group and all
other factors finite cyclic.  If \(X\) is simplified, if \(J\) is a
component of \(\bd Y-\bd X\), and if \(d\) is the degree of the
attaching homeomorphism for the disk in \(X-Y\) attached to \(J\),
then the element of \(\pi_1(X)\) represented by \(J\) has order
exactly \(d\).  \end{lemma}

\begin{demo}{Proof} The complex can be assumed to be simplified for
all parts of the statemtent.  The surface \(Y\) has boundary
components \(J_1,\dots,J_n\) to which disks are attached, and
boundary components \(K_1,\dots,K_m\) to which no disks are
attached.  By hypothesis, \(m\ne0\).  There are pairwise disjoint,
properly embedded arcs \(\alpha_1,\dots,\alpha_n\) in \(Y\) with
each \(\bd\alpha_i\) having one point on \(K_1\) and one point on
\(J_i\).  A regular neighborhood \(A_i\) of \(J_i\cup\alpha_i\) in
\(Y\) is an annulus with \(J_i\) as one boundary component, and with
frontier in \(Y\) a properly embedded arc.  The closure \(Y'\) of
\(Y-\bigcup(A_i)\) is a surface with boundary.  Then \(X-Y'\) is a
disjoint union of Fuchsian complexes \(B_1,\dots,B_n\) where each
\(B_i\) is made from the annulus \(A_i\) by attaching to \(J_i\) the
disk attached to \(J_i\) in \(X\).  The results follow immediately.
\end{demo}

If \(X\) is a Fuchsian complex based on a surface \(Y\) and \(S\) is
a set of compact boundary components of \(Y\) to which disks are
attached, then we can form an identification space of \(X\) by
taking to a point each component in the collection \(S\) and all the
disks that are attached along that component.  We refer to this
identification as \defit{modding out the disks on \(S\)}.  If \(S\)
is the collection of all components of \(\bd Y\) that have disks
attached then we see that modding out the disks on \(S\) gives a
surface whose boundary is the boundary of \(X\).

\begin{lemma} Let \(X\) be a Fuchsian complex based on a surface
\(Y\) and let \(S\) be a set of components of \(\bd Y\) to which
disks have been attached.  Then the identification that mods out the
disks on \(S\) has image a Fuchsian complex and induces a surjection
on \(\pi_1\) whose kernel is normally generated by the loops in
\(S\).  If \(\pi_1(X)\) is torsion free, then the identification
induces an isomorphism.  \end{lemma}

\begin{demo}{Proof} The image of each collection of disks attached
to a curve in \(S\) is a point.  Thus the image of \(Y\) in the
identification is a surface \(Y'\) in which each curve in \(S\) has
become a point, and the image \(X'\) of \(X\) is a Fuchsian complex
based on the surface \(Y'\).  Let \(x_1,\dots,x_n\) be the image
points of the curves in \(S\).  A loop in \(X'\) can be homotoped to
miss the \(x_i\) and we see that the identification induces a
surjection on \(\pi_1\).  A loop in the kernel can be homotoped in
\(X\) to miss all attached disks and thus reside in \(Y\).  Its
image in \(Y'\) bounds a singular disk that may contain points
\(x_i\) in its interior.  Thus the original loop in \(Y\) bounds a
singular disk with holes whose other boundary components are
elements of \(S\).

If \(\pi_1(X)\) is torsion free, then each curve in \(S\) must
represent the trivial element in \(\pi_1(X)\).  The kernel of the
homomorphism on \(\pi_1\) must then be trivial.  \end{demo}

\begin{cor}\label{TorFreeFuchs} A torsion free Fuchsian group is a
surface group.  \end{cor}

\begin{demo}{Proof} Let \(X\) be a Fuchsian complex based on a
surface \(Y\), with \(\pi_1(X)\) torsion free.  If we mod out the
all the attached disks, then we do not change the fundamental group.
Modding out all the attached disks gives a surface.  \end{demo}

\begin{cor} If a finite Fuchsian complex \(X\) without boundary has
\(\pi_1(X)\) torsion free, then \(\pi_1(X)\) is the fundamental
group of a closed surface.  \end{cor}

\begin{demo}{Proof} Modding out the attached disks gives a compact
surface without boundary whose fundamental group is isomorphic to
\(\pi_1(X)\).  \end{demo}

\begin{cor} If a finite Fuchsian complex \(X\) with boundary has
\(\pi_1(X)\) torsion free, then \(\pi_1(X)\) is free.  \end{cor}

\begin{demo}{Proof} Modding out the attached disks gives a compact
surface with boundary whose fundamental group is isomorphic to
\(\pi_1(X)\).  \end{demo}

We say that a Fuchsian group is \defit{orientable} if it is the
fundamental group of a Fuchsian complex that is based on an
orientable surface \(Y\).  There are non-orientable Fuchsian groups.

\begin{lemma} A Fuchsian group with an infinite subgroup isomorphic
to the fundamental group of a closed, non-orientable surface is
non-orientable.  \end{lemma}

\begin{demo}{Proof} Let \(X\) be a Fuchsian complex based on a
surface \(Y\) and let \(\pi_1(X)\) have an infinite subgroup \(G\)
isomorphic to the fundamental group of a closed, non-orientable
surface \(S\).  Let \(\tilde X\) be the cover of \(X\) corresponding
to \(G\).  It is a Fuchsian complex based on the surface \(\tilde
Y\), the pre-image of \(Y\) in \(\tilde X\).  Since \(G\) is torsion
free, the curves in \(\bd\tilde Y\) with attached disks must
represent the trivial element in \(G\).  Modding out the attached
disks must result in a surface with fundamental group \(G\).  This
surface must be non-orientable, so \(\tilde Y\) and \(Y\) must be
non-orientable.  \end{demo}

\begin{lemma}\label{SeifFuchs} Let \(M\) be a Seifert fibered space,
let \(h\) be represented by an ordinary fiber, and let \(p:M\into
G\) be the projection to the orbit surface.  Then there is a
Fuchsian complex \(X\) with an isomorphism \(i:\pi_1(M)/\langle
h\rangle\into\pi_1(X)\) so that if \(q:X\into X'\) is the
identification map obtained by modding out the attached disks, then
there is a homeomorphism \(f:X'\into G\) taking the images of the
attached disks to the exceptional points of \(G\) so that the
following commutes. \[ \xymatrix{\pi_1(M) \ar[r] \ar[d]_{p_{\#}} &
\pi_1(M)/\langle h\rangle \ar[r]^i & \pi_1(X) \ar[dl]^{q_{\#}} \\
\pi_1(G) & \pi_1(X') \ar[l]_{f_{\#}}} \] \end{lemma}

\begin{demo}{Proof} The presentation \ref{PiOneOrbl} or
\ref{PiOneNOrbl} for \(\pi_1(M)\) is built from the orbit surface
\(G\) from which disk neighborhoods of the exceptional points and a
single disk neighborhood of an ordinary point have have been
removed.  Let the resulting surface be \(G'\), let \(c_1,\dots,c_n\)
be the boundary components of \(G'\) that are also the boundaries of
the disk neighborhoods of the exceptional points, and let \(e\) be
the boundary component of \(G'\) that is also the boundary of the
removed disk neighborhood of the ordinary point.  The presentation
of the Fuchsian quotient \ref{FuchsOrbl} or \ref{FuchsNOrbl} is
obtained if disks are attached to the curves \(c_i\) by maps of
degree \(\mu_i\) on the boundaries and a disk is attached to \(e\)
by a homeomorphism of the boundary.  The commutativity of the
diagram can now be checked.  \end{demo}

\begin{cor} If \(M\) is a Seifert fibered space, \(h\) is
represented in \(\pi_1(M)\) by an ordinary fiber, and
\(\pi_1(M)/\langle h\rangle\) is a non-orientable Fuchsian group,
then the orbit surface of \(M\) is non-orientable.  \end{cor}

\begin{demo}{Proof} This follows from the previous lemma.
\end{demo}

The assumption in the next lemma that \(\pi_1(X)\) be infinite is
seen necessary by considering a Fuchsian complex based on a \Mob\ to
which a single disk has been added to the boundary with a map of
degree \(d\).  The centerline of the \Mob\ represents an element of
order \(2d\) and is not conjugate to a power of the boundary of the
\Mob.

\begin{lemma}\label{TorInFuchs} Let \(X\) be a simplified Fuchsian
complex constructed from a surface \(Y\) and assume that
\(\pi_1(X)\) is infinite.  If \(\alpha\ne1\) is an element of
\(\pi_1(X)\) of finite order, then a conjugate of \(\alpha\) is
represented by a loop in \(\bd Y - \bd X\).  Two different
components of \(\bd Y-\bd X\) represent elements that are not
conjugate in \(\pi_1(X)\).  Also, an element of \(\pi_1(X)\)
represented by a boundary component of \(X\) has infinite order.
\end{lemma}

\begin{demo}{Proof} In proving the first two assertions we can
assume that the boundary of \(X\) is empty.  Let \(C\) be the cyclic
subgroup generated by \(\alpha\).  The cover \(\tilde X\)
corresponding to \(C\) is determined by the conjugacy class of
\(C\).  Thus we must show that some loop in the boundary of \(Y\)
lifts to \(\tilde X\).

Since \(C\) is finite and \(\pi_1(X)\) is infinite, \(\tilde X\) is
an infinite sheeted cover.  Thus \(\tilde Y\), the preimage of \(Y\)
is non-compact.  Modding out all attached disks gives a connected,
non-compact surface with empty boundary and whose fundamental group
is the image of a finite cyclic group.  Thus this surface is an open
disk and \(\tilde Y\) is planar.  We know that each component \(J\)
of \(\bd\tilde Y\) is torsion in the subspace of \(\tilde X\) formed
by \(J\) and the disks attached to it.  We thus see that
\(\pi_1(\tilde X)\) is the free product of cyclic groups.  But
\(\pi_1(\tilde X)\) is cyclic, so exactly one of the free factors is
non-trivial and exactly one component of \(\bd\tilde Y\) carries the
fundamental group of \(\tilde X\).  This finishes the first
assertion.

To consider the second assertion, let \(\alpha\) and \(\beta\) be
elements represented by two components \(J_1\) and \(J_2\)
respectively of \(\bd Y-\bd X\).  We again consider the cover
\(\tilde X\) corresponding to the cyclic group \(C\) generated by
\(\alpha\).  If \(\beta\) is conjugate to \(\alpha\), then both
\(J_1\) and \(J_2\) have lifts to \(\tilde X\).  But this would give
\(\pi_1(\tilde X)\) the structure of a non-trivial free product of
cyclic groups.

To consider the third assertion, we return to the original Fuchsian
complex \(X\) without discarding any of its boundary components.
Let \(J\) be a component of \(\bd X\). If we mod out all the
attached disks, we get a surface with boundary and \(J\) is one of
the boundary components.  If the quotient is a disk, then
\(\pi_1(X)\) is a free product of more than one cyclic group (since
\(\pi_1(X)\) is infinite), and \(J\) is the product of the
generators.  But this is an element of infinite order.  If the
quotient surface is not a disk, then its fundamental group is
infinite and \(J\) represents an element of infinite order in it.
Thus \(J\) represents an element of infinite order in \(\pi_1(X)\).
\end{demo}

Note that the ``small'' Seifert fibered spaces are those with finite
Fuchsian quotient.  This accounts for the hypothesis in the
following consequence of the above lemma.

\begin{lemma}\label{PowFib} Let \(S\) be a Siefert fibered space
that is not one of the ``small'' spaces.  Let \(\alpha\) be an
element of \(\pi_1(S)\) so that a power of \(\alpha\) is homotopic
to a power of an ordinary fiber.  Then a conjugate of \(\alpha\) is
represented by a power of an exceptional fiber.  If we assume in
addition that \(\alpha\) is represented by a loop \(J\) in a torus
boundary component \(C\) of \(S\), then \(J\) is homotopic in \(C\)
to a power of a fiber in \(C\).  \end{lemma}

\begin{demo}{Proof} As usual let \(H_i\) be the exceptional fibers
and \(N_i\) be pairwise disjoint fibered solid torus neighborhoods
of the \(H_i\) in \(S\).  Let \(S'\) be obtained from \(S\) by
removing the interiors of the \(N_i\).  Let \(h\) be represented by
an ordinary fiber.  The Fuchsian complex for \(\pi_i(S)/\langle
h\rangle\) is obtained from the orbit surface \(G'\) of \(S'\) by
attaching disks to the curves \(J_i\) that are the images of the
\(\bd N_i\).  The hypothesis says that \(\alpha\) maps to a torsion
element in \(\pi_1(S)/\langle h\rangle\).  Since the Fuchsian
quotient is infinite, Lemma \ref{TorInFuchs} says that the image of
\(\alpha\) is conjugate to a loop in one of the \(J_i\).  Thus
\(\alpha\) is conjugate mod \(\langle h\rangle\) to a loop in one of
the \(\bd N_i\).  But \(h\) can be represented by a loop in that
\(\bd N_i\) so \(\alpha\) is conjugate to a loop in that \(\bd
N_i\).  A loop in \(\bd N_i\) is either trivial or a power of the
centerline of \(N_i\).

If \(\alpha\) is represented by a loop in a boundary component \(C\)
of \(S\), then Lemma \ref{TorInFuchs} says that the image of
\(\alpha\) in \(\pi_1(S)/\langle h\rangle\) is either trivial in the
image of \(\pi_1(C)\), or has infinite order.  The hypothesis rules
out the latter.  The homomorphism from \(\Pi_1(C)\) is the
homomorphism from \(\ints\oplus\ints\) to \(\ints\) that kills the
factor generated by an ordinary fiber.  This completes the proof.
\end{demo}

The next result is proven by R. H. Fox in ``On Fenchel's conjecture
about F-groups,'' {\it Matematisk Tidsskrift,} B (1952), 61--65.  It
will be used to prove the lemmas that follow about torsion in
fuchsian groups.

\begin{thm} Let \(p\ge2\), \(q\ge2\) and \(r\ge2\) be three
integers.  Then there is a finite group with elements \(a\) and
\(b\) so that the orders of \(a\), \(b\) and \(ab\) are exactly
\(p\), \(q\) and \(r\) respectively.  \BlackBox\end{thm}

\begin{lemma} Let \(X\) be a simplified Fuchsian complex based on a
surface \(Y\) with \(\pi_1(X)\) infinite.  Let \(J\) be a component
of \(\bd Y-\bd X\) and let \(d\) be the degree of the attaching
homeomorphism for the disk in \(X-Y\) attached to \(J\).  Then the
element of \(\pi_1(X)\) represented by \(J\) has order exactly
\(d\).  \end{lemma}

\begin{demo}{Proof} If \(X\) has boundary, then the result is proven
above.  We assume that \(X\) has no boundary.

If there is a non-separating, two sided simple closed curve \(K\) in
\(Y\), then splitting \(X\) and \(Y\) along \(K\) gives a Fuchsian
complex \(X'\) based on a surface with at least two boundary
components and infinite fundamental group so that \(\pi_1(X')\) is a
free product of an infinite free group and a collection of finite
groups that includes \(\ints_d\) generated by \(J\).  The two copies
of \(K\) in \(X'\) would represent elements of infinite order in
\(\pi_1(X')\).  Now \(\pi_1(X)\) is an HNN extension into which
\(\pi_1(X')\) includes.

If there is a separating simple closed curve \(K\) in \(Y\) so that
the two complementary domains of \(K\) in \(X\) had infinite
\(\pi_1\), then a similar argument proceeds with a free product with
amalgamation.

If there is no such curve, then modding out the attached disks must
result in a 2-sphere or projective plane.  Thus \(Y\) is a disk with
holes or a \Mob\ with holes.  If the number of holes in the disk is
zero or one, or the number of holes in the \Mob\ is zero, then
\(\pi_1(X)\) is finite.  If the number of holes in the disk is
three, then the result follows from Fox's theorem and if the number
of holes is four or more, then a serpating curve \(K\) can be found
fitting the requirements of the paragraphs above.  If the number of
holes in the \Mob\ is one or more, then we can also find such a
curve.  (If an arc is run from the boundary of the \Mob\ to the
boundary of the hole, then the curve sought is the boundary of a
regular neighborhood of the union of the arc, the boundary of the
\Mob\ and the boundary of the hole.)  \end{demo}

\begin{cor} If \(M\) is a Seifert fibered space, \(h\) is
represented by an ordinary fiber, and \(\pi_1(M)/\langle h\rangle\)
is non-trivial and torsion free, then \(M\) has no exceptional
fibers.  \end{cor}

\begin{demo}{Proof} We know that \(\pi_1(M)/\langle h\rangle\) is
infinite.  We can construct a simplified Fuchsian complex with
fundamental group \(\pi_1(M)/\langle h\rangle\) from the
presentations \ref{FuchsOrbl} or \ref{FuchsNOrbl} as appropriate.
If there are exceptional fibers, then there will be disks attached
with homomeorphisms of degree above 1.  From the lemma, we know that
there will be torsion in the corresponding Fuchsian group.
\end{demo}

\begin{lemma}\label{TriangleSub} Every triangle group has a torsion
free, normal subgroup of finite index.  \end{lemma}

\begin{demo}{Proof} A finite group obviously has a torsion free,
normal subgroup of finite index, so we restrict our attention to
infinite groups.  Let \(\Gamma(p,q,r)\) be an infinite triangle
group presented by \(\langle x,y\mid x^p=y^q=(xy)^r=1\rangle\) and
let \(G\) be the group given by Fox's theorem.  There is a
homomorphism \(h\) taking \(x\) to \(a\) and \(y\) to \(b\).  The
only powers of \(x\), \(y\) and \(xy\) in the kernel of \(h\) are
the identity element.  Now \(\Gamma(p,q,r)\) is a Fuchsian group
with Fuchsian complex based on a compact, planar surface with three
boundary components corresponding to \(x\), \(y\) and \(xy\).  By
Lemma \ref{TorInFuchs}, we know that any torsion element of
\(\Gamma(p,q,r)\) is conjugate to a power of \(x\), \(y\) or \(xy\).
Thus no torsion elment of \(\Gamma(p,q,r)\) is in the kernel of
\(h\).  Since \(G\) is finite, we have that the kernel of \(h\) is a
torsion free subgroup of finite index in \(\Gamma(p,q,r)\).
\end{demo}

\begin{thm}\label{SurfInFuchs} Every Fuchsian group corresponding to
a finite Fuchsian complex contains the fundamental group of a
surface as a normal subgroup of finite index.  \end{thm}

\begin{demo}{Proof} By Corollary \ref{TorFreeFuchs}, we need only
find a normal, torsion free subgroup of \(\pi_1(X)\) of finite
index.  If we find a torsion free subgroup of finite index that is
not normal, then the intersection of its finitely many conjugates
will be torsion free, normal and of finite index.  Thus we are free
to replace \(X\) by a finite cover of \(X\) at any time.

Note that if \(\pi_1(X)\) is finite then the trivial subgroup
satisfies the conclusion.  We assume that \(\pi_1(X)\) is infinite.

Let \(X\) be a finite Fuchsian complex based on a surface \(Y\).  If
\(Y\) is not orientable, then we can consider the orientable double
cover of \(Y\).  The boundary components of \(Y\) are orientation
preserving, and lift to the double cover.  Thus we can construct a
double cover of \(X\) based on the orientable double cover of \(Y\).
This gives a subgroup of index two in \(\pi_1(X)\) whose Fuchsian
complex is based on an orientable surface.  Thus we assume that
\(Y\) is orientable.

We also assume that the complex \(X\) is simplified so that no more
than one disk is attached along the same curve, and every attaching
map is of degree at least two.

If we mod out by the attached disks, then we obtain an orientable
surface \(Y'\) and a finite number of points \(x_1,\dots,x_n\) that
are the images of the attached disks.  We consider cases based on
the genus of \(Y'\) and the number \(n\).

If \(n=0\), then \(X=Y=Y'\) and \(\pi_1(X)\) is an infinite surface
group and torsion free.

If \(n\ge3\), then there are a finite number of disks
\(D_1,\dots,D_m\) in \(Y'\) so that no \(x_i\) is in a \(\bd D_j\),
and so that each \(D_j\) contains exactly 3 of the \(x_i\).  The
\(D_j\) are not necessarily disjoint and some of the \(x_i\) will be
in more than one \(D_j\).  For each \(j\), let \(E_j\) be the
preimage of \(D_j\) in \(X\), let \(F_j\) be the closure of
\(X-E_j\), and let \(Z_j\) be obtained from \(X\) by identifying
\(F_j\) to a point.  Each \(Z_j\) is a Fuchsian complex based on a
disk with two holes and its fundamental group is a triangle group.

Each \(Z_j\) contains exactly three boundary components of \(Y\).
Let \(p\), \(q\) and \(r\) be the degrees of the maps for the disks
attached to these three boundary components.  We know that each
degree is at least two and we know that the corresponding power of
each boundary component is trivial in \(\pi_1(X)\).  In
\(\pi_1(Z_j)\), \(p\), \(q\) and \(r\) are exactly the orders of the
three boundary components.  Thus the only power of any of these
boundary components that is in in the kernel of the induced
homomorphism from \(\pi_1(X)\) to \(\pi_1(Z_j)\) is the identity
element.  By Lemma \ref{TriangleSub}, there is a torsion free,
normal subgroup of finite index in \(\pi_1(Z_j)\).  Let \(N_j\) be
the preimage of this subgroup in \(\pi_1(X)\).  The subgroup \(N_j\)
is of finite index in \(\pi_1(X)\).  Note that no power other than
the identity of the three relevant components of \(\bd Y\) are in
\(N_j\).

Let \(N\) be the intersection of the \(N_j\).  The subgroup \(N\) is
of finite index in \(\pi_1(X)\).  Let \(\alpha\) be a torsion
element of \(\pi_1(X)\).  By Lemma \ref{TorInFuchs}, \(\alpha\) is
conjugate to a power of a component of \(\bd Y\) to which a disk is
attached.  Thus \(\alpha\) is not in some \(N_j\) and not in \(N\).
Thus \(N\) is torsion free.

We now assume that \(0<n<3\) and consider the possibilities for
\(Y'\).  If \(Y'\) is a disk and \(n\le 1\), then \(\pi_1(X)\) is
finite cyclic.  If \(Y'\) is a disk and \(n=2\), then \(\pi_1(X)\)
is a free product of two cyclic groups.  We can either appeal to the
fact that a free product of cyclic groups has a free subgroup of
finite index, or we can use the followoing argument.  If both cyclic
groups are of order two, then \(\pi_1(X)\) is \(D_{\infty}\) and we
know that is has an infinite cyclic subgroup of index two.  Assume
one order is at least three.  We obtain \(X\) from a disk with two
holes by attaching a disk \(D\) to one hole with a map of some
degree \(d\) and attaching a disk \(E\) to the other hole by a map
of degree \(e\) with \(e\ge3\).  Now \(X-E\) has an \(e\)-fold cover
obtained from an annulus with \(e\) holes by attaching one disk to
each hole with a map of degree \(d\).  We obtain an \(e\)-fold cover
\(\tilde X\) of \(X\) by attaching \(e\) disks to one of the
remaining boundary components with degree one maps.  Now \(\tilde
X\) has at least \(e\) attaching sites of disks with attaching map
of degree more than one and we are done by referring to the case
where \(n\) was at least three.

If \(Y'\) is a sphere with \(n=1\) or \(n=2\), then \(\pi_1(X)\) is
finite.

Since \(Y'\) is compact and orientable, then any possibility for
\(Y'\) other than a disk or 2-sphere has \(\pi_1(Y')\) infinite.
Thus \(Y'\) has finite sheeted covers of arbitrarily high index.
Assuming \(n\ge1\), we can imitate our construction of \(\tilde X\)
above to obtain finite sheeted covers of \(X\) with any number of
attached disks of degree more than 1 and again refer to the case
where \(n\) was at least three.  \end{demo}

The sign of the Euler characteristic of the surface in Theorem
\ref{SurfInFuchs} is an invariant of the group.  If \(X\) is a
Fuchsian complex and \(\tilde X\) is a \(\lambda\)-sheeted cover of
\(X\) with fundamental group \(H\) isomorphic to the fundamental
group of the surface \(S\), then we can associate the rational
number \(\chi(S)/\lambda\) to \(\pi_1(X)\).  If we similarly obtain
another subgroup \(K\) of finite index \(\gamma\) in \(\pi_1(X)\)
and a surface \(S'\) with fundmaental group \(K\), then \(H\cap K\)
has index \(\gamma\) in \(H\) and index \(\lambda\) in \(K\).  There
is a surface \(T\) that is a \(\gamma\)-fold cover of \(S\) and a
surface \(T'\) that is a \(\lambda\)-fold cover of \(S'\) whose
fundamental groups correspond to \(H\cap K\).  Since
\(\pi_1(T)\cong\pi_1(T')\) and both are infinite, we have that \(T\)
and \(T'\) are homotopy equivalent and \(\chi(T)=\chi(T')\).  Since
Euler characteristics multiply with the sheetedness of covers,
\(\gamma(\chi(S))=\chi(T)=\chi(T')=\lambda(\chi(S'))\) and
\(\chi(S)/\lambda=\chi(S')/\gamma\).

We can compute the sign of the Euler characteristic from the
presentation of the Fuchsian group.  Let the Fuchsian complex \(X\)
be based on the surface \(Y\).  The surface \(Y\) has either \(g\)
handles or \(k\) crosscaps.  The surface \(Y\) also has \(n\)
boundary components to which disks are attached, and \(m\) boundary
components that remain boundary components of \(X\).  The Euler
charasteristic of \(Y\) is given by \(\chi(Y)=2-s-n-m\) where \(s\)
is either \(2g\) or \(k\).  Assume that \(\pi_1(X)\) has a torsion
free subgroup of index \(\lambda\).  Let \(\tilde X\) be the cover
of \(X\) determined by this subgroup.  Let \(\tilde Y\) be the
preimage of \(Y\) in \(\tilde X\).  We know that \(\tilde Y\) is a
connected \(\lambda\)-fold cover of \(Y\) and has \(\chi(\tilde
Y)=\lambda(2-s-n-m)\).  We get \(\tilde X\) from \(\tilde Y\) by
attaching disks to some of the boundary components of \(\tilde Y\).
Let \(J_1,\dots,J_n\) be the components of \(\bd Y\) to which disks
are attached in \(X\).  Let \(d_i,\dots,d_n\) be the degrees of the
repsective attaching maps.  (We assume that the Fuchsian complex
\(X\) is simplified.)  Each component of the preimage of \(J_i\) in
\(\tilde Y\) covers \(J_i\) by a \(d_i\)-fold covering map.  Thus
there are \(\lambda/d_i\) components of the preimage of \(J_i\).  To
complete \(\tilde Y\) to a surface \(S\) that carries the
fundamental group of \(\tilde X\), we must add exactly one disk (by
a homeomorphism of the boundary) to each preimage of each \(J_i\).
Thus we will add \(\lambda/d_i\) disks for each \(J_i\).  We get
\begin{align*}\chi(S)&=\lambda(2-s-n-m)+\sum\limits^n\frac{\lambda}{d_i}
\\&=\lambda\left[(2-s-m)-\sum\limits^n(1-\frac{1}{d_i})\right]\end{align*}
where again \(s\) is either \(2g\) or \(k\) depending on whether
\(Y\) is orientable or not.

We say that a Fuchsian group has \defit{positive, negative or zero
Euler characteristic} depending on the Euler characteristic of the
surface given by Theorem \ref{SurfInFuchs}.  Since surfaces of
positive Euler characteristic have finite fundamental groups, we
know that Fuchsian groups of positive Euler characteristic are
finite.  We can also decribe those Fuchsian groups of zero Euler
characteristic.

We let \(g\), \(k\), \(n\), \(m\) and \(d_1,\dots,d_n\) be as in the
previous paragraphs.  We let \(\Sigma\) denote
\[\sum\limits^n(1-\frac{1}{d_i})\] and let \(s\) denote \(2g\) or
\(k\) as appropriate.  We get Euler characteristic zero when
\(s+m+\Sigma=2\).  Since \(g=0\) and \(k=0\) both correspond to
2-spheres, we have the following possibilities where ``holes'' are
boundary components that have attached disks.

\begin{tabbing}\hspace{1in}\= \(S^2\) with 3 or 4 holes
\hspace{.25in}\= \(\Sigma=2\) \+ \\ disk with 2 holes \>
\(\Sigma=1\) \\ annulus \> \(\Sigma=0\) \\ torus \> \(\Sigma=0\) \\
\(P^2\) with 2 holes \> \(\Sigma=1\) \\ \Mob \> \(\Sigma=0\) \\
Klein bottle \> \(\Sigma=0\) \end{tabbing}

Ther are only a finite number of ways to get the values of
\(\Sigma\) above.  For \(\Sigma=2\), the values of the \(d_i\) must
be one of the following \((2,2,2,2)\), \((2,3,6)\), \((2,4,4)\) or
\((3,3,3)\).  For \(\Sigma=1\) the only combination possible is
\((2,2)\).  This is of interest because of the following.

\begin{lemma} If the fundamental group \(G\) of a finite Fuchsian
complex has a subgroup isomorphic to fundamental group of a torus or
a Klein bottle, then the Euler characteristic of \(G\) is zero.
\end{lemma}

\begin{demo}{Proof} Let \(H\) be the assumed subgroup.  A subgroup
of finite index in \(H\) is also one of the two possibilities for
\(H\).  This is seen by considering covering spaces of the
corresponding surfaces.  Thus the subgroup of finite index in \(G\)
that is the fundamental group of a surface \(S\) has a subgroup that
is one of the two possibilities for \(H\).  The cover of \(S\)
corresponding to \(H\) must be either a torus or Klein bottle.  This
forces the Euler characteristic of \(S\) to be zero.  \end{demo}

\begin{lemma} Every infinite Fuchsian group contains an element of
infinite order.  \end{lemma}

\begin{demo}{Proof} Let \(X\) be a Fuchsian complex based on a
surface \(Y\).  If \(\bd X\) is not empty, then we are done by Lemma
\ref{TorInFuchs}.  Let \(Y'\) be obtained by modding out all the
attached disks.  If \(Y'\) is not orientable, then we can pass to a
2-sheeted cover and assume that \(Y'\) is orientable.  If
\(\pi_1(Y')\) is infinite, then it has an element of infinite order
and so does \(\pi_1(X)\).  Thus we can assume that \(Y'\) is an
orientable, surface without boundary and finite fundamental group.
Thus \(Y'\) is an open disk or 2-sphere.

If \(Y'\) is a disk then \(\pi_1(X)\) is a free product of cyclic
groups and will be infinite if and only if at least two non-trivial
cyclic groups are involved.  If at least two non-trivial cyclic
groups are involved, then there will be an element of infinite
order.

If \(Y'\) is a 2-sphere, then \(X\) is a finite complex and we are
done by Theorem \ref{SurfInFuchs}.  \end{demo}

The next lemma discusses Fucshian groups that have \(\ints\) as a
subgroup of finite index.  Having \(\ints\) as a subgroup of finite
index is equivalent to having two ends.

\begin{lemma}\label{TwoEndFuchs} If the fundamental group \(G\) of a
finite Fuchsian complex has \(\ints\) as a subgroup of finite index,
then \(G\) is \(\ints\) or \(\ints_2*\ints_2\).  \end{lemma}

\begin{demo}{Proof} Let \(X\) be a finite, simplified Fuchsian
complex with fundamental group \(G\) based on a surface \(Y\).  Let
\(\tilde X\) be the cover corresponding to the \(\ints\) subgroup of
finite index.  If \(\bd\tilde X\) is empty, then
\(\ints=\pi_1(\tilde X)\) is the fundamental group of a closed
surface, so \(\bd\tilde X\) and \(\bd X\) must be non-empty.

The hypothesis implies that the number of ends of \(\pi_1(X)\) is
two.  Since \(\bd X\) is non empty, \(\pi_1(X)\) is a free product
of groups, one of which is free (and perhaps trivial) and all others
finite cyclic.  The number of ends will be infinite if the free
group factor has rank more than one, if there are more than two
non-trivial factors, or if there are two non-trivial factors and one
of the factors has order more than two.  Since \(\pi_1(X)\) is not
finite, the only possibilities left are those of the conclusion.
\end{demo}

\subsection{Centralizers in Seifert fibered spaces}

We know that the fundamental group of a Seifert fibered space has a
cyclic normal subgroup.  When this subgroup is infinite, the action
on it by conjugation can only carry a generator to itself or its
inverse.  Thus the cyclic normal subgroup is central in an index 1
or 2 subgroup in the fundamental group.  We can say exactly what
this subgroup is.

\begin{lemma} Let \(M\) be a Seifert fibered space with infinite
fundamental group, let \(h\) be represented by an ordinary fiber,
let \(p:M\into G\) be the projection to the orbit surface, and let
\(\phi:\pi_1(G)\into\ints_2\) be the classifying homomorphism.  Then
the centralizer of \(h\) in \(\pi_1(M)\) is the kernel of \(\phi
p_{\#}:\pi_1(M)\into\ints_2\).  \end{lemma}

\begin{demo}{Proof} Since \(\pi_1(M)\) is infinite, \(h\) has
infinite order in \(\pi_1(M)\) and is not equal to its inverse.
From the presentations \ref{PiOneOrbl} and \ref{PiOneNOrbl}, we see
that a generator \(a_i\), \(b_i\), or \(c_i\) is in \(C\), the
centralizer of \(h\) in \(\pi_1(M)\), if and only if it is in the
kernel \(N\) of \(\phi p_{\#}\).  Also, both \(N\) and \(C\) contain
\(h\).  We thus have two subgroups of index 1 or 2 that contain
exactly the same generators.  The two subgroups now yield the same
homomorphisms to \(\ints_2\) and are the same.  \end{demo}

We would like to know what other elements might have centralizers as
large.

We need one technical lemma.

\begin{lemma} If \(n\) is an integer greater than 1, then
\(\ints\oplus\ints_n\) is not a Fuchsian group.  \end{lemma}

\begin{demo}{Proof} Let \(X\) be a Fuchsian complex based on a
surface \(Y\) with \(\pi_1(X)\cong\ints\oplus\ints_n\).  We may
assume that \(\bd X\) is empty.  If we mod out the attached disks,
we get a surface without boundary whose fundamental group is a
quotient of \(\ints\oplus\ints_n\).  This surface must be \(S^2\),
\(P^2\) or an open disk, annulus or \Mob.

If the quotient surface is \(P^2\), then there is a subgroup of
index 2 with corresponding surface \(S^2\).  We will show that
neither \(\ints\) nor \(\ints\oplus\ints_n\) can occur if the
quotient surface is \(S^2\).  This will eliminate both \(S^2\) and
\(P^2\).  If the surface is \(S^2\), then as argued before,
\(\pi_1(X)\) is either finite cyclic or a triangle group or a free
product with amalgamation of two non-ableian groups (depending on
the number of attaching sites of disks in \(X\)).  However, all the
infinite triangle groups are non-abelian.

If the quotient surface is an open disk, then \(\pi_1(X)\) is a free
product of finite cyclic groups.  This cannot yield
\(\ints\oplus\ints_n\).  If the quotient surface is an open annulus
or \Mob, then \(\pi_1(X)\) is the free product of \(\ints\) and a
collection of finite cyclic groups.  Again, \(\ints\oplus\ints_n\)
cannot be obtained.  \end{demo}

We can now determine the infinite Fuchsian groups that have
non-trivial center.

\begin{lemma} The only infinite Fuchsian groups that have
non-trivial center are \(\ints\), \(\ints\oplus\ints\) or the
fundamental group of the Klein bottle.  \end{lemma}

\begin{demo}{Proof} Let \(G\) be an infinite Fuchsian group with
non-trivial center.  We start by showing that \(G\) is torsion free.
We know that \(G\) has an element of infinite order.  If it has
torsion, then it also has an element of finite order.  Let \(a\) be
an element of the center.  Let \(b\) be another element of \(G\)
chosen to have finite order if \(a\) has infinite order, and chosen
to have infinite order if \(a\) has finite order.  We know \(a\) and
\(b\) commute since \(a\) is in the center, so \(a\) and \(b\)
generate a subgroup of \(G\) isomorphic to \(\ints\oplus\ints_n\)
for some \(n>1\).  Since a subgroup of a Fuchsian group is a
Fuchsian group, this contradicts the previous lemma.

We now have that \(G\) is torsion free and must be the fundamental
group of a surface \(S\).  If \(S\) is a non closed surface, then
\(G\) is trivial or free.  The only infinite such group with center
is \(\ints\).  If \(S\) is a closed, orientable surface of genus at
least 2 or a closed, non-orientable surface with at least 3
crosscaps, then \(G\) is the free product with amalgamation of one
non-ableian free group with another free group of rank at least one
amalgamated over \(\ints\).  Since the generator of the amalgamating
subgroup does not commute with elements of the non-abelian free
group, the group \(G\) has trivial center.  The only remaining
surfaces have fundamental groups listed in the statement of the
theorem.  \end{demo}

The groups in the above lemma show up when looking at the Fuchsian
quotient of the fundamental group of the Seifert fibered space.  The
case in which we will apply the above, we will be modding out by a
cyclic central subgroup.  The next lemma analyzes what can happen.

\begin{lemma} Let \(G\) be a group and let \(N\) be an infinite,
cyclic, central subgroup of \(G\).  Assume that \(G/N\) is (i)
\(\ints\), (ii) \(\ints\oplus\ints\), or (iii) the fundamental group
of a Klein bottle.  Then in case (i), \(G\) is \(\ints\oplus\ints\),
in case (ii), \(G\) is \(\ints\oplus\ints\oplus\ints\) or \(N\) is
the entire center of \(G\), and in case (iii), \(G\) has an index 2
subgroup isomorphic to \(\ints\oplus\ints\oplus\ints\) which
contains \(N\).  \end{lemma}

\begin{demo}{Proof} The first case is elementary.  In the other
cases, \(G\) is generated by \(t\) a generator of \(N\), and two
other elements \(a\) and \(b\).  We have that \(t\) commutes with
\(a\) and \(b\) and that \(a^{-1}ba=b^{\epsilon}t^m\) for some
\(m\in\ints\) and \(\epsilon=1\) in case (ii) and \(\epsilon=-1\) in
case (iii).  If \(m=0\), then \(G\) is isomorphic to
\(\ints\oplus(G/N)\) which satisfies the conclusions of the two last
cases.  If \(m\ne0\) in case (iii), then
\(a^{-2}ba^2=a^{-1}b^{-1}at^m=bt^{-m}t^m=b\) and the index two
subgroup in question is generated by \(a^2\), \(b\) and \(t\).  We
now assume that \(m\ne0\) in case (ii) and show that every central
element of \(G\) is in \(N\).

Since \(t\) commutes with \(a\) and \(b\) and \(ba=abt^m\), we can
pass \(b\) over \(a\) at the expense of introducing a power of
\(t\).  Thus we can write any element of \(G\) as \(a^pb^qt^r\).  If
this element is central, then so is \(a^pb^q\) and \(a^p\) commutes
with \(b\) and \(b^q\) commutes with \(a\).  We have
\(a^{-1}b^qa=b^qt^{mq}\) so \(t^{mq}=b^{-q}a^{-1}b^qa=1\).
Similarly, \(t^{mp}=1\) is derived from \(bab^{-1}=at^m\).  But
\(t\) has infinite order and we assume \(m\ne0\), so both \(p\) and
\(q\) are 0.  Thus the central element \(a^pb^qt^r\) is in \(N\).
\end{demo}

The initial hypotheses in the next lemma could be stated by
requiring that \(S\) not be ``small.''

\begin{lemma}\label{TwoFibGrps} Let \(S\) be a compact, connected
Seifert fibered space with infinite, torsion free fundamental group
and infinte Fuchsian quotient.  Let \(G\) be the orbit surface of
\(S\), let \(h\) be represented by an ordinary fiber, let \(C\) be
the centralizer of \(h\) in \(\pi_1(S)\), and let \(\alpha\) not a
power of \(h\) generate a cyclic normal subgroup of \(\pi_1(S)\).
Then \(\ints\oplus\ints\) or \(\ints\oplus\ints\oplus\ints\) is a
subgroup of finite index in \(\pi_1(S)\) containing \(\langle
h\rangle\).  The possible indexes are given by the following table.
\newline\parbox{\textwidth}{\begin{tabbing} \qquad\(C=\pi_1(S)\) and
\(G\) is orientable\qquad\qquad \= 1, 2, or 4\qquad\qquad \= \kill
\> \(\underline{\ints\oplus\ints}\) \>
\(\underline{\ints\oplus\ints\oplus\ints}\) \\ no extra assumptions
\> 1, 2 or 4 \> 1, 2, 4 or 8 \\ \(C=\pi_1(S)\) \> 1 or 2 \> 1, 2 or
4 \\ \(G\) is orientable \> 1, 2 or 4 \> 1, 2 or 4 \\ \(C=\pi_1(S)\)
and \(G\) is orientable \> 1 or 2 \> 1 or 2
\end{tabbing}}\vspace{5pt}\newline If
\(\ints\oplus\ints\oplus\ints\) is not a subgroup of \(\pi_1(S)\)
with one of the indicated indexes, then \(\pi_1(S)/\langle
h\rangle\) is either \(\ints\) or \(\ints_2*\ints_2\).  \end{lemma}

\begin{demo}{Proof} Let \(N\) be generated by \(h\).  We know that
\(N\) is infinite, and that the index of \(C\) in \(\pi_1(S)\) is 1
or 2.  Similarly, the centralizer \(A\) of \(\alpha\) is of index 1
or 2 in \(\pi_1(S)\).  If \(\alpha\) is not in \(C\), then no power
of \(h\) commutes with \(\alpha\), no power of \(h\) is in \(A\),
and \(A\) is not of finite index in \(\pi_1(S)\).  Thus \(\alpha\)
is in \(C\).  By hypothesis, \(\alpha\) is not in \(N\).

Since \(\alpha\) is in \(C\), \(\alpha\) commutes with \(h\) and all
its powers, so \(N\) is in \(A\).  Also, \(N\) is in \(C\).  Let
\(\overline{\alpha}\), \(\overline{A}\) and \(\overline{C}\) be the
images of \(\alpha\), \(A\) and \(C\) in \(\pi_1(S)/N\).  Since
\(\pi_1(S)/N\) is infinite, both \(\overline{A}\) and
\(\overline{C}\) are infinite.  In particular, \(\overline{A}\) is
an infinite Fuchsian group with \(\overline{\alpha}\) as a
non-trivial central element.  Thus \(\overline{A}\) is one of
\(\ints\), \(\ints\oplus\ints\), or the fundamental group of a Klein
bottle.

Let \(H=A\cap C\) and let its image in \(\pi_1(S)/N\) be
\(\overline{H}\).  We know that \(N\) is in \(H\) and is central in
it.  Also, \(H\) is of index 1 or 2 in \(A\), and \(\overline{H}\)
is of index 1 or 2 in \(\overline{A}\).  So \(\overline{H}\) is also
one of \(\ints\), \(\ints\oplus\ints\) or the fundamental group of a
Klein bottle.  We also know that \(N\) is not the entire center of
\(H\) since \(\alpha\notin N\) is central in \(H\).  By the previous
lemma, \(H\) is \(\ints\oplus\ints\) or has
\(\ints\oplus\ints\oplus\ints\) as an index 1 or 2 subgroup
containing \(N\).  The possibility that \(H=\ints\oplus\ints\) only
occurs when \(\overline{H}\) is \(\ints\).  But \(\overline{H}\) is
a subgroup of finite index in \(\pi_1(S)/\langle h\rangle\).  This,
and lemma \ref{TwoEndFuchs} give the last conclusion.

We must calculate indexes.  All indexes discussed are powers of 2.
The index of \(C\) in \(\pi_1(S)\) is 1 or 2.  The index of \(H\) in
\(A\) is the index of \(C\) in \(\pi_1(S)\) and the index of \(A\)
in \(\pi_1(S)\) is 1 or 2.  Thus the index of \(H\) in \(\pi_1(S)\)
is no more than 4 and if \(C=\pi_1(S)\), then it is no more than 2.
This accounts for the possible indexes if \(H=\ints\oplus\ints\) or
\(H=\ints\oplus\ints\oplus\ints\).  The only way for \(H\) not to be
\(\ints\oplus\ints\oplus\ints\) but have
\(\ints\oplus\ints\oplus\ints\) as a subgroup of index 2 in \(H\) is
for \(\overline{H}\) to be the fundamental group of a Klein bottle.
Then \(\pi_1(S)/\langle h\rangle\) is a non-orientable Fuchsian
group and the orbit surface of \(S\) is non-orientable.  This
accounts for the remaining possibilities.  \end{demo}

The last lemma says that for ``large'' Seifert fibered space, the
absence of \(\ints\oplus\ints\) and \(\ints\oplus\ints\oplus\ints\)
as subgroups of the fundamental group implies that the cyclic normal
subgroup generated by an ordinary fiber is the unique maximal,
cyclic normal subgroup of the fundamental group.  The presence of
\(\ints\oplus\ints\) or \(\ints\oplus\ints\oplus\ints\) as subgroups
limits the Fuchsian quotient to have Euler characteristic zero.
This limits the structure of the Seifert fibered space.  The fact
that these subgroups are of finite index, says that there is a
finite cover of the space having these subgroups as fundamental
group.  This is even more of a restriction.  In the next section we
will consider the structure of finite covers of Seifert fibered
spaces.  Before that we extract a little more information from our
knowledge of the fundamental group.

\begin{lemma}\label{TwoFibCombos} Let \(S\) be a compact, connected
Seifert fibered space with infinite, torsion free fundamental group
and infinite Fuchsian quotient.  Let \(h\) be represented by an
ordinary fiber,and let \(\alpha\) not a power of \(h\) generate a
cyclic normal subgroup of \(\pi_1(S)\).  Then the orbit surface of
\(S\) and the exceptional fibers are limited to the following
combinations.\newline
\parbox{\textwidth}{\begin{tabbing}\hspace{1in} \= orbit
surface\hspace{1in} \= exceptional fibers \+ \\ \vspace{10pt} \(S^2\)
\> four of index 2 \\ \(D^2\) \> two of index 2 \\ \(S^1\x I\) \>
none \\ \(S^1\x S^1\) \> none \\ \(P^2\) \> two of index 2 \\ \Mob
\> none \\ Klein bottle \> none \end{tabbing}} \end{lemma}

\begin{demo}{Proof} We know that the Fuchsian group
\(\pi_1(S)/\langle h\rangle\) is built from the orbit surface \(G\)
of \(S\) by removing neighborhoods of the exceptional points and
attaching disks with maps of various degrees to the resulting
boundaries.  The list in the conclusion is the list giving the
Fuchsian groups of zero Euler characteritic with some omissions.  We
must show that \(\pi_1(S)/\langle h \rangle\) has zero Euler
characteristic and that the omissions are legitimate.

We know that \(\pi_1(S)\) has a subgroup of finite index containing
\(\langle h\rangle\) isomorphic to \(\ints\oplus\ints\) or
\(\ints\oplus\ints\oplus\ints\).  In the latter case,
\(\ints\oplus\ints\) is a subgroup of \(\pi_1(S)/\langle h\rangle\)
and \(\pi_1(S)/\langle h\rangle\) has zero Euler characteristic.  In
the former case, \(\pi_1(S)/\langle h\rangle\) is \(\ints\) or
\(\ints_2*\ints_2\).  These Fuchsian groups can only be built from
an annulus or \Mob\ with no exceptional points, or from a disk two
exceptional points of index two.  These are also included on the
list of possibilities with zero Euler characteristic.

We must show that the items omitted from the list are not possible
under the hypotheses.  The omitted items all have orbit surface the
2-sphere, and three exceptional fibers with indexes \((2,3,6)\),
\((2,4,4)\) or \((3,3,3)\).  Since the orbit surface is \(S^2\), we
are not in the situation analyzed in the previous paragraph where
\(\ints\oplus\ints\) is a subgroup of finite index in \(\pi_1(S)\).
The classifiying homomorphism is trivial, so \(h\) is central in
\(\pi_1(S)\).  Also, \(S^2\) is orientable.  Thus
\(K=\ints\oplus\ints\oplus\ints\) is a subgroup of index 1 or 2 in
\(\pi_1(S)\).  Since \(h\) is in \(K\), the image \(\overline{K}\)
in \(\pi_1(S)/\langle h\rangle\) is of index 1 or 2 in
\(\pi_1(S)\langle h\rangle\).  We cannot have
\(\overline{K}=\ints\oplus\ints\oplus\ints_n\) since
\(\ints\oplus\ints_n\) is not a Fuchsian group.  Thus
\(\ints\oplus\ints\) is a subgroup of index 1 or 2 in
\(\pi_1(S)/\langle h\rangle\).  However, \(\pi_1(S)/\langle
h\rangle\) is the triangle group \(\Gamma(2,3,6)\),
\(\Gamma(2,4,4)\) or \(\Gamma(3,3,3)\).  These all have torsion
elements of order more than 2 and cannot have a torsion free
subgroup of index 1 or 2.  \end{demo}

\subsection{Covers of Seifert fibered spaces}

%checked to here next page is 140

In the previous section, we saw that certain assumptions about a
Seifert fibered space implied that the fundamental group would have
certain groups as subgroups of finite index.  This has strong
implications for the structure of the covers of the Seifert fibered
spaces corresponding to these subgroups.  This motivates the topic
of this section --- the structure of covers of Seifert fibered
spaces.

Let \(M\) be a Seifert fibered space.  Finite sheeted covers of
\(M\) turn out to inherit a Seifert fibered structure from \(M\).
The invariants of this fiber structure are computable, and in some
cases easily computable.

In this section we will have to be careful and distinguish between
lifts and pre-images.  If \(p:X\into Y\) is a covering projection
and \(J\) is a circle in \(Y\), then \(p^{-1}(J)\) is the pre-image
of \(J\).  This may or may not be connected and a given component of
\(p^{-1}(J)\) may or may not be a circle.  We get a lift of \(J\) by
chosing a point in \(J\) as a basepoint and an orientation so that
we can regard \(J\) as a closed path.  A lift of \(J\) is a path in
\(X\) starting at a pre-image of the basepoint that covers \(J\)
once as a path.  It is possible that no lift of \(J\) is a circle
while some or all components of the pre-image of \(J\) are circles.
The fact that a component of the pre-image of \(J\) is a circle
implies that some power of \(J\) (regarded as an element of
\(\pi_1(Y)\)) lifts to a closed path in \(X\).

Let \(p:\tilde M\into M\) be a covering projection and let \(H\) be
a fiber in \(M\).  Then \(p^{-1}(H)\) is a 1-manifold and each
component is a line or a circle.  This fibers \(\tilde M\) with
lines and circles.  Let \(N\) be a fibered solid torus neighborhood
of \(H\), let \(H'\) be a component of \(p^{-1}(H)\) and let \(N'\)
be the component of \(p^{-1}(N)\) that contains \(H'\).  The degree
of \(p|{N'}\) equals the degree of \(p|{H'}\).  This degree is
finite precisely when \(H'\) is a circle, or equivalently, when
\(N'\) is a solid torus.  Let \(d(H')\) denote the degree of
\(p|{N'}\).

If the degree of \(p|{N'}\) is finite, then for any fiber of \(M\)
in \(N\), a component of the pre-image of the fiber that intersects
\(N'\) lies in \(N'\) and is a circle.  Thus a component of the
pre-image of a fiber that is a circle has a neighborhood of
pre-images of fibers that are circles.  Similarly, a component of
the pre-image of a fiber that is a line has a neighborhood of
pre-images of fibers that are lines.  Thus in \(\tilde M\), the
components of pre-images that are circles and the components of
pre-images that are lines form disjoint open sets that cover
\(\tilde M\).  If \(\tilde M\) is connected, one of these sets is
empty.  We have shown:

\begin{lemma} If \(M\) is a Seifert fibered space and \(p:\tilde
M\into M\) is a covering projection with \(\tilde M\) connected for
which one fiber \(H\) of \(M\) has a component of \(p^{-1}(H)\) a
circle, then \(\tilde M\) inherits a Seifert fiber structure from
\(M\).  \BlackBox\end{lemma}

From now on we only consider covers of \(M\) where all components of
the pre-images of fibers are circles.  This happens in finite
sheeted covers of Seifert fibered spaces, but can happen in other
ways as well.  Note that \(G\x S^1\) with \(G\) a surface with
infinite fundamental group is infinitely covered by \(R^2\x S^1\).

Let \(H\) now be an ordinary fiber of \(M\) and let \(H'\), \(N\)
and \(N'\) be as above.  The pre-images of fibers in \(N'\) fiber
\(N'\) as an ordinary solid torus and the restriction of \(p\) to
each fiber in \(N'\) has the same degree as \(p|{H'}\).  Thus the
fuction \(d\) is constant over all the fibers in \(N'\).  We also
have that the pre-images in \(\tilde M\) of the exceptional fibers
are isolated.  Thus the pre-images in \(\tilde M\) of the ordinary
fibers of \(M\) form a connected subset of \(\tilde M\) which can be
covered by pre-images of ordinary solid tori in \(M\).  This implies
that the fuction \(d\) is constant over all the pre-images in
\(\tilde M\) of ordinary fibers in \(M\).  Let this common value be
called the \defit{degree} of \(p\) and be denoted \(d(p)\).

Now let \(H\) be an exceptional fiber of \(M\) and let \(H'\), \(N\)
and \(N'\) be as above.  Let the fiber structure of \(N\) be
determined by \(\nu/\mu\).  Let \(\sigma=d(H')\).  The fiber
structure of \(N\) is derived from a rotation of the unit disk by
\(2\pi(\nu/\mu)\).  The fiber structure of \(N'\) is derived from a
rotation of the unit disk by \(2\pi\sigma(\nu/\mu)\).  Assuming that
\(\nu/\mu\) is in reduced terms, the rational number in reduced
terms that determines \(N'\) is
\[\frac{\sigma\nu/(\sigma,\mu)}{\mu/(\sigma,\mu)}\] where
\((\sigma,\mu)\) is the greatest common divisor of \(\sigma\) and
\(\mu\).  Let \(g=(\sigma,\mu)\), so the invariants of \(N'\) are
\((\sigma\nu/g)/(\mu/g)\).

Let \(J\) be a simple closed curve on \(\bd N\).  In terms of a
meridian-longitude pair \((m,l)\) we can write \(J=\alpha m-\beta
l\).  We write it this way because with
\(\hint(m,l)=-\hint(l,m)=1\), we have \(J=\hint(J,l)m-\hint(J,m)l\)
and \(\hint(J,m)=\beta\).

If we think of \(N\) as made from one copy of \(D\x I\), then we can
think of \(N'\) as made from \(\sigma\) copies of \(D\x I\).  We can
declare the ``length'' of a curve \(J\) in \(\bd N\) to be
\(\hint(J,m)\), the algebraic intersection of the curve with \(m\).
We can declare the ``length'' of a curve in \(\bd N'\) to be the
algebraic intersection of the curve with the full pre-image in
\(N'\) of \(m\).  A full pre-image in \(N'\) of \(m\) will consist
of \(\sigma\) disjoint, parallel meridians of \(N'\).  This gives
\(J\) length \(\beta\) in \(N\).  The full pre-image of \(J\) in
\(N'\) has length \(\sigma\beta\).

The subgroup of \(\pi_1(N)\) corresponding to the covering map from
\(N'\) to \(N\) is \(\sigma\ints\).  Since \(J\) represents
\(\beta\) times a generator of \(\pi_1(N)\), the smallest power of
\(J\) that lies in \(\sigma\ints\) is \([\sigma,\beta]/\beta\) where
\([\,\,,\,]\) is the least common multiple.  Since
\([\sigma,\beta]=\sigma\beta/(\sigma,\beta)\), we have that the
smallest power of \(J\) that lies in \(\sigma\ints\) is
\(\sigma/(\sigma,\beta)\).  Such a power of \(J\) will have length
\(\sigma\beta/(\sigma,\beta)\) in \(N\) and a lift of this power of
\(J\) will be a simple closed curve of the same length in \(N'\).
Thus each component of the pre-image of \(J\) in \(N'\) has length
\(\sigma\beta/(\sigma,\beta)\), the full pre-image of \(J\) in
\(N'\) has length \(\sigma\beta\), and \(J\) is covered by
\begin{equation}(\sigma,\beta)\label{CurveLiftCount}\end{equation}
components.

Let \(J'\) be a component of the pre-image of \(J\) in \(N'\).  The
intersection of \(J'\) with the full pre-image of \(m\) is the
length of \(J'\) and is \(\sigma\beta/(\sigma,\beta)\).  But \(m\)
lifts to \(\sigma\) copies of \(m\).  Thus the intersection of
\(J'\) with a meridian \(m'\) of \(N'\) is \(\beta/(\sigma,\beta)\).

The pre-image \(l'\) of \(l\) in \(N'\) is connected, is a longitude
of \(N'\) and has length \(\sigma\).  The intersection of \(J\) with
\(l\) is \(\alpha\), and a calculation similar to the one just done
gives that the intersection of \(J'\) with \(l'\) is
\(\alpha\sigma/(\sigma,\beta)\).

Thus the homology class of \(J'\) is determined as
\begin{equation}J'=\hint(J',l')m'-\hint(J',m')l'=
\frac{\alpha\sigma}{(\sigma,\beta)}m'-
\frac{\beta}{(\sigma,\beta)}l'.\label{CurveLiftForm}\end{equation}

If we apply \ref{CurveLiftCount} and \ref{CurveLiftForm} to an
ordinary fiber of form \(\nu m+\mu l\) in \(\bd N\), we get that it
is covered by \((\sigma,\mu)\) components (which are ordinary fibers
in \(\bd N'\)) and that each component is of the form \[
\frac{\nu\sigma}{(\sigma,\mu)}m'+ \frac{\mu}{(\sigma,\mu)}l'.\] This
agrees with our analysis of the type of \(N'\) above.

Recall that \(g=(\sigma,\mu)\).  Ordinary fibers in \(N\) are
covered by \(g\) ordinary fibers in \(N'\).  The restriction of
\(p\) to \(N'\) induces a map of the orbit surface of \(N'\) to the
orbit surface of \(N\).  Both surfaces are disks, and off the center
point, the map is \(g\) to one.  Thus the map is a branched cover of
degree \(g\) with branch point the center of the disk.

Since each ordinary fiber in \(N\) is covered by \(g\) ordinary
fibers in \(N'\), and each ordinary fiber in \(M\) is covered
\(d(p)\) times by each component of its pre-image in \(M\), we have
that each ordinary fiber in \(N\) is covered \(gd(p)\) times by its
full pre-image in \(N'\).  We get that \(gd(p)=\sigma\) and
\(\sigma\) is a multiple of \(d(p)\).

Let \(A\) be the subgroup of \(\pi_1(M)\) that corresponds to the
covering \(p:\tilde M\into M\), and let \(h\) be an element of
\(\pi_1(M)\) determined by an ordinary fiber.  Letting \(d=d(p)\),
we know that \(h^d\) is in \(A\) and no smaller power is.  In other
words, \(A\cap\langle h\rangle\) is a subgroup of index \(d\) in
\(\langle h\rangle\).  If \(\langle h\rangle\) is in \(A\), then
\(d=1\) and every ordinary fiber lifts to a closed path.  Because of
the results in the previus section, we will be interested in covers
where ordinary fibers lift to closed paths.  Let us call such a
cover a \defit{primary} cover of \(M\).

If \(p\) is a primary covering map, then \(d=d(p)=1\) and for a
component of the pre-image \(H'\) of an exceptional fiber with
\(\sigma\) as above, we have \(g=gd=\sigma\) where
\(g=(\mu,\sigma)\).  Thus \(\sigma=(\mu,\sigma)\) and
\(\sigma|\mu\).  The induced map of orbit surfaces is a branched
cover of degree \(g=\sigma\) near the corresponding exceptional
point.

We are primarily interested in primary covering maps.  This is
because the subgroups discussed in the previous section included the
cyclic normal subroup generated by an ordinary fiber.  The covers
corresponding to such subgroups are primary covers.

Let \(H\) be an exceptional fiber with invariant \(\nu/\mu\), and
let \(N\), \(N'\), \(m\), \(\sigma\) and \(g=(\sigma,\mu)\) be as
above.  Let \(Q\) be a crossing curve on \(\bd N\).  Now the
meridian \(m=\mu Q+\beta H\) is a simple closed curve, so
\((\mu,\beta)=1\).  Since the cover is primary, \(\sigma|\mu\) and
\((\sigma,\beta)=1\).  We have \(Q=\alpha m+\beta l\) for some
\(\alpha\) and from \ref{CurveLiftCount} and \ref{CurveLiftForm} we
get that the pre-image of \(Q\) in \(N'\) is connected and has
intersection \(\beta\) with a meridian of \(N'\).

We are going to look at two types of covers.  One will be a special
type of primary cover in which not only the ordinary fibers lift to
circles, but all fibers lift to circles.  This makes the handling of
the crossing curves determined by a fibered solid torus easy.  We
will see examples of this type of cover when we give two
applications.  One will be the orientable double cover of a
non-orientable Seifert fibered space.  The second will be a double
cover of a space with a non-orientable orbit surface in which the
cover has an orientable orbit surface.  The second type of cover
will be one in which there are no exceptional fibers in the covering
space.  The crossing curves are easy in this case as well.  That
there are such covers is shown by the following.

\begin{lemma} Let \(M\) be a Seifert fibered space, let
\(h\in\pi_1(M)\) be represented by an ordinary fiber and assume that
\(\pi_1(M)/\langle h\rangle\) is infinite.  Then \(M\) has a primary
finite sheeted cover with no exceptional fibers in the induced
fibration.  \end{lemma}

\begin{demo}{Proof} Let \(h\in\pi_1(M)\) be represented by an
ordinary fiber.  The Fuchsian group \(\pi_1(M)/\langle h\rangle\)
has a torsion free subgroup \(K\) of finite index.  We write
\(\overline{K}\) to represent the pre-image in \(\pi_1(M)\) of
\(K\).  Then \(\overline{K}\) has finite index in \(\pi_1(M)\).  The
cover \(\tilde M\) of \(M\) corresponding to \(\overline{K}\) is a
Seifert fibered space and each ordinary fiber of \(M\) lifts to an
ordinary fiber in \(\tilde M\).  To avoid extra notation with
induced homomorphisms, we regard \(\pi_1(\tilde M)\) as the subgroup
\(\overline{K}\) of \(\pi_1(M)\).  If \(\tilde h\in\pi_1(\tilde M)\)
is represented by an ordinary fiber of \(\tilde M\), then
\(\langle\tilde h\rangle=\langle h\rangle\).  The Fuchsian quotient
\(\pi_1(\tilde M)/\langle\tilde h\rangle\) is just
\(\overline{K}/\langle h\rangle=K\) and is torsion free and
infinite.  Thus \(\tilde M\) has no exceptional fibers.  \end{demo}

Seifert fiberings of the 3-sphere with exceptional fibers show that
the hypotheses of the previous lemma are necessary.

\begin{cor} The cover \(\tilde M\) of the previous lemma can be
assumed to be orientable.  \end{cor}

\begin{demo}{Proof} Fibers are orientation preserving curves in a
Seifert fibered space.  The orientable double cover of a Seifert
fibered space is a primary cover and a primary cover of a primary
cover is a primary cover.  Lastly, covers cannot introduce
exceptional fibers where there were none before.  Taking the
orientable double cover of the result of the previous lemma gives
the result.  \end{demo}

We now consider a cover of a Seifert fibered space in which all
lifts of all fibers are circles.  We first assume that the base is
oriented.  We consider the non-orientable case later when we look
only at the orientable double cover.

Let \(M\) be a connected, oriented Seifert fibered space and let
\(p:\tilde M\into M\) be a covering projection as described.  We
assume that \(M\) and \(\tilde M\) are compact so that the number of
sheets \(\lambda\) of the cover is finite.  The number \(\sigma\)
defined above which gives the degree of the projection on a fiber in
the cover is 1 for all fibers in \(\tilde M\).  We recall more
notation for the structure of \(M\).  We let \(H_1,\dots,H_n\) be
the exceptional fibers of \(M\), we let \(N_1,\dots,N_n\) be
pairwise disjoint fibered solid torus neighborhoods of the \(H_i\)
with meridians \(m_1,\dots,m_n\), we let
\((\mu_1,\beta_1),\dots,(\mu_n,\beta_n)\) be the crossing
invariants, and let \(Q_1,\dots,Q_n\) be the crossing curves
determined by the \(N_i\).  We let \(N_0\) with meridian \(m_0\) be
an ordinary solid torus neighborhood of some ordinary fiber in \(M\)
with \(N_0\) disjoint from the \(N_i\).  We let \(Y\) be the image
of a section of the orbit surface of the closure \(M_0\) of
\(M-(N_0\cup N_1 \cup \dots \cup N_n)\) with \(\bd Y\cap\bd
N_i=Q_i\) for \(i\ge1\).  We let \(Q_0=\bd Y\cap\bd N_0\).  We know
that \(m_0=Q_0+bH\) and \(m_i=\mu_i Q_i+\beta_i H\) where \(H\)
represents an ordinary fiber in the appropriate torus.  We let
\(\tilde M_0\) be the preimage of \(M_0\) in \(\tilde M\).  It is
connected since it is obtained from \(\tilde M\) by removing the
pre-images of the \(N_i\).  These are solid tori with connected
boundary and if \(\tilde M_0\) is not connected, then \(\tilde M\)
cannot be connected.  The preimage of \(Y\) in \(\tilde M_0\) is a
surface \(\tilde Y\) which hits each fiber of \(\tilde M_0\) in a
single point.  Thus \(\tilde Y\) is connected.

Each \(N_i\), \(i\ge1\), lifts to \(\lambda\) copies of itself that
cover \(N_i\) once each.  Thus each component of \(p^{-1}(N_i)\)
determines a crossing curve that is a lift of \(Q_i\).  This is just
a boundary component of \(\tilde Y\).  Each component of a lift of
\(N_0\) has a meridian with intersection \(b\) with the appropriate
boundary component of \(\tilde Y\).

If we remove the lifts of the \(N_i\), \(i\ge1\), from \(\tilde M\)
and replace them with ordinary solid tori with meridians going to
the lifts of the \(Q_i\), then we would be in a position to
calculate the obstruction to the section.  The partial section would
contain the surface \(\tilde Y\) together with the meridians bounded
by the \(Q_i\).  If \(N_0\) lifted to a single solid torus, we would
be in good position.  Unfortunately it lifts to \(\lambda\)
different tori.  We need a lemma.

\begin{lemma} Let \(M\) be a closed, connected, oriented Seifert
fibered space with no exceptional fibers.  Let
\(N_1,\dots,N_\lambda\) be pairwise disjoint saturated solid tori in
\(M\), let \(M_0\) be the closure of \(M-(N_1\cup\dots\cup
N_\lambda)\), and let \(Y\) be the image of a section for \(M_0\)
with boundary components \(Q_i=\bd Y\cap\bd N_i\).  Let
\(b_1,\dots,b_\lambda\) be the intersections of the meridians
\(m_i\) of the \(N_i\) with the \(Q_i\).  Then the obstruction to
the section for \(M\) is \(\sum b_i\).  \end{lemma}

\begin{demo}{Proof} If all the \(b_i\), \(i>1\), are zero, then we
are done since \(Y\) extends through the \(N_i\), \(i>1\).  We will
be done by induction when we show that we can replace \(b_1\) by
\(b_1+b_2\) and \(b_2\) by zero while leaving the \(b_i\), \(i>2\)
the same and obtain a manifold of the same fiber type as \(M\).  Let
\(\alpha\) be an arc in \(Y\) from \(Q_1\) to \(Q_2\).  The
saturated annulus \(A\) over \(\alpha\) has two fibers for its
boundary.  We can orient these fibers to have the ``same'' direction
in \(A\).  If we orient \(\bd N_1\) and \(\bd N_2\) consistently
with the orientations that \(N_1\) and \(N_2\) inherit from \(M\)
(the convention from the first chapter), then \(Q_1\) and \(Q_2\)
have intersections of opposite sign with the fibers in \(\bd A\).
There is a self homeomorphism of \(M_0\) fixed off a neighborhood of
\(A\) that rotates one ``side'' of \(A\) through \(b_2\) full
rotations.  The direction can be chosen to take \(Q_2\) to \(Q_2+b_2
H=m_2\).  This will take \(Q_1\) to \(Q'_1=Q_1-b_2 H\).  This makes
\(m_2\) and \(Q'_1\) the new boundary curves of reference, and the
new intersections are calculated from \(m_2=m_2+0H\) and
\(m_1=Q_1+b_1 H=Q_1-b_2 H+b_2 H+b_1 H=Q'_1+(b_1+b_2)H\).  This
completes the proof.  \end{demo}

We can now prove the following.

\begin{lemma} Let \(M\) be a compact, connected, oriented Seifert
fibered space and let \(p:\tilde M\into M\) be a primary,
\(\lambda\)-sheeted cover with \(\lambda\) finite.  Let \(\tilde M\)
have the induced fiber structure.  If all fibers of \(M\) lift to
circles in \(\tilde M\) then each exceptional fiber \(H\) is covered
by \(\lambda\) exceptional fibers of \(\tilde M\) with the same
crossing invariants as \(H\).  If \(M\) is closed with obstruction
to section \(b\), then the obstruction to section for \(\tilde M\)
is \(\lambda b\).  The orbit surface of \(\tilde M\) is a
\(\lambda\)-sheeted cover of the orbit surface of \(M\), and the
classifying homomorphism for \(\tilde M\) is obtained by composing
the homomorphism induced by the projection with the classifying
homomorphism for \(M\).  \end{lemma}

\begin{demo}{Proof} The previous lemma gives the statement about the
obstruction to the section.  We have demonstrated above the facts
about the crossing invariants.  The remaining statements are
straightforward.  \end{demo}

We now consider a non-orientable Seifert fibered space \(M\) with
\(n\) exceptional fibers.  We look at the orientable double cover
\(\tilde M\) of \(M\).  Since every fiber in \(M\) is orientation
preserving in \(M\), it lifts only to circles in \(\tilde M\).  We
adopt the notation developed for the case of oriented \(M\).  Each
exceptional fiber \(H_i\) with fibered solid torus neighborhood
\(N_i\) has crossing invariants \((\mu_i,\beta_i)\) with
\(0\le\beta_i\le\mu_i/2\).  There will be two components of the
pre-image of \(N_i\) in \(\tilde M\) so \(\tilde M\) will have
\(2n\) exceptional fibers.  The covering translation of \(\tilde M\)
is an orientation reversing involution which will carry one
component of the pre-image of \(N_i\) to the other.  Thus one
component will have crossing invariants \((\mu_i,\beta_i)\) and the
other will have crossing invariants \((\mu_i,\mu_i-\beta_i)\).  This
determines the exceptional fibers of \(\tilde M\) and the collection
of crossing invariants.  It turns out that the obstruction to the
section \(\tilde b\) for \(\tilde M\) is remarkably easy to
calculate.  The existence of a fiber preserving, orientation
reversing homeomorphism says that \(\tilde M\) is fiber equivalent
to \(\tilde M\) with the opposite orientation.  The result at the
end of chapter 1 says that \(\tilde b\) must be equal to
\(-2n-\tilde b\) since there are \(2n\) exceptional fibers in
\(\tilde M\).  Thus \(\tilde b=-n\) and is independent of the
crossing invariants.  We have proven.

\begin{lemma} Let \(M\) be a compact, connected, non-orientable
Seifert fibered space with \(n\) exceptional fibers and let \(\tilde
M\) be the orientable couble cover of \(M\).  Then \(\tilde M\) has
\(2n\) exceptional fibers.  Each crossing invariant pair
\((\mu_i,\beta_i)\) for \(M\) is replace by two pairs
\((\mu_i,\beta_i)\) and \((\mu_i,\mu_i-\beta_i)\) for \(\tilde M\)
and if \(M\) is closed, then the obstruction to a section for
\(\tilde M\) is \(-n\).  The orbit surface for \(\tilde M\) is a
double cover of the orbit surface for \(M\) and the classifying
homomorphism for \(\tilde M\) must be the orientation homomorphism
for its orbit surface.  \BlackBox\end{lemma}

We now consider primary covers with no exceptional fibers.  We
restrict ourselves to oriented base.  We continue our previous
notation and keep the same meanings for \(M\), \(M_0\), \(\tilde
M\), \(\tilde M_0\), \(H_i\) and \((\mu_i,\beta_i)\) for \(1\le i\le
n\), \(N_i\), \(Q_i\) and \(m_i\), for \(0\le i\le n\), \(Y\),
\(\tilde Y\) and \(\lambda\).  Each \(H_i\) has a set of numbers
\(\sigma_{ij}\) where \(H_{ij}\) are the components of the pre-image
of \(H_i\) and \(\sigma_{ij}\) is the degree of the map from
\(H_{ij}\) to \(H_i\).  Since the cover is a \(\lambda\)-sheeted
map, we must have \(\sum\limits_j\sigma_{ij}=\lambda\) for each
\(i\).  Each component of the pre-image of \(N_i\) is an ordinary
solid torus.  We let \(N_{ij}\) denote the component containing
\(H_{ij}\).  Since we have a primary cover, the analysis done at the
beginning of this section shows that each crossing curve \(Q_i\) has
one component in its pre-image in each \(N_{ij}\) and that it has
intersection \(\beta_i\) with a meridian of \(N_{ij}\).

Since we have a primary cover, we know that each
\(\sigma_{ij}|\mu_i\).  However, the index of the fiber \(H_{ij}\)
is \(\mu_i/(\mu_i,\sigma_{ij})\).  Since the cover has no
exceptional fibers, this must be 1 and \((\mu_i,\sigma_{ij})\) must
be \(\mu_i\).  Thus \(\mu_i|\sigma_{ij}\) and \(\sigma_{ij}=\mu_i\)
for all \(j\).  Thus every \(N_{ij}\) is a \(\mu_i\)-fold cover of
\(N_i\) and there are \(\lambda/\mu_i\) of them.  The total
intersection of the meridians of the \(N_{ij}\) with the lifts of
\(Q_i\) (boundary components of a section for \(\tilde M_0\)) is
\(\lambda\beta_i/\mu_i\).  We have previously shown that the total
intersection of the meridians of the lifts of \(N_0\) with the
boundaries of a section for \(\tilde M_0\) is \(\lambda b\).  From
the lemma above on the obstruction to a section, we have that the
obstruction to a section for \(\tilde M\) is \[\lambda
b+\sum\frac{\lambda\beta_i}{\mu_i}=
\lambda\left(b+\sum\frac{\beta_i}{\mu_i}\right).\] We can now extend
one of the lemmas above.

\begin{lemma} Let \(M\) be a compact, connected, oriented seifert
fibered space with \(n\) exceptional fibers with crossing invariants
\((\mu_i,\beta_i)\), \(1\le i\le n\).  If \(M\) is closed, let the
obstruction to a section for \(M\) be \(b\).  Let \(h\in\pi_1(M)\)
be represented by an ordinary fiber.  Let \(\pi_1(M)/\langle
h\rangle\) be infinite and have a torsion free subgroup if finite
index \(\lambda\).  Then \(M\) has a \(\lambda\)-sheeted cover
\(\tilde M\) with no exceptional fibers and, if \(M\) is closed, it
has obstruction to a section given by
\[\lambda\left(b+\sum\frac{\beta_i}{\mu_i}\right).\]
\BlackBox\end{lemma}

\begin{remark} The quantity \[b+\sum\frac{\beta_i}{\mu_i}\] is, in
fact, more natural than \(b\) itself.  This has been exploited in
various papers (e.g., Neumann and Raymond, Lecture Notes in Math.
664).  \end{remark}

\subsection{Uniqueness of fibers}

We can combine the results of the previous sections to give a
statement about the uniqueness of the elements of the fundamental
group represented by the ordinary fibers of \defit{any} Seifert
fibration of a Seifert fibered space.  Such a statement needs
restrictions since fiber structures are not unique for certain
spaces.  Once we understand the uniqueness of the elements of the
fundamental group that can be represented by an ordinary fiber, we
can give results about the uniqueness of the fiber structure of
certain Seifert fibered spaces.  We will only do this under the
assumption that the spaces have non-empty boundary even though it is
possible to prove similar statements about closed manifolds.

To help with the next statement we define certain Seifert fibered
spaces as ``flat.''  The include some closed manifolds and some with
boundary.  The closed manifolds are those determined by the
following: \[(O,o,0\mid-2,(2,1),(2,1),(2,1),(2,1)),\]
\[(O,o,1\mid0),\] \[(N,o,1\mid\epsilon),\]
\[(O,n,1\mid-1,(2,1),(2,1)),\] \[(N,n,I,1\mid-,(2,1),(2,1)),\]
\[(O,n,2\mid0),\] \[(N,n,I\hbox{--}II,2\mid\epsilon)\] where
\(\epsilon\) is either 0 or 1.  The manifolds with boundary are the
unique space with orbit surface a disk and two expectional fibers of
index 2, and the four spaces (orientable and non-orientable) with no
exceptional fibers with orbit surface either an annulus or a \Mob.

Recall that certain Seifert fibered spaces were identified in
\ref{SmallSFS} as ``small.''  These are the spaces with finite
Fuchsian quotient.  Thus the hypotheses in the next statement
require that the space be neither ``small'' nor ``flat.''

\begin{lemma}\label{HomFiber} Let \(M\) be a compact, connected
Seifert fibered space, let \(h\in\pi_1(M)\) be represented by an
ordinary fiber and assume that \(\pi_1(M)/\langle h\rangle\) is
infinite.  Further assume that \(M\) is not ``flat.''  Then any
cyclic normal subgroup of \(\pi_1(M)\) is contained in \(\langle
h\rangle\).  \end{lemma}

\begin{demo}{Proof} If the conclusion is false, then the hypotheses
of Lemma \ref{TwoFibCombos} are satisfied and \(M\) is one of the
spaces listed in the conclusion of that lemma.  We also know that
\(\pi_1(M)\) contains \(\ints\oplus\ints\) or
\(\ints\oplus\ints\oplus\ints\) as a subgroup \(N\) of finite index
as described in Lemma \ref{TwoFibGrps}.

We note that the ``flat'' Seifert fibered spaces are some of the
spaces in the conclusion of Lemma \ref{TwoFibCombos}.  We must
explain why not all of the spaces in the conclusion of Lemma
\ref{TwoFibCombos} are listed among the ``flat'' spaces.

All of the spaces with boundary in the conclusion of Lemma
\ref{TwoFibCombos} are listed as ``flat.''  All non-orientable
spaces in the conclusion of Lemma \ref{TwoFibCombos} are listed as
``flat.''  The remaining spaces in Lemma \ref{TwoFibCombos} only
show up among the ``flat'' spaces with specific values for the
obstruction to a section.  We must show why these are the only
possible values under the hypotheses of the present lemma.

We assume that \(M\) is closed and orientable.  Let \(N\) be the
subgroup of finite index guaranteed by Lemma \ref{TwoFibGrps}, and
let \(\tilde M\) be the cover of \(M\) corresponding to \(N\).  We
first argue that \(N\) is \(\ints\oplus\ints\oplus\ints\).  If not,
then \(N\) is \(\ints\oplus\ints\).  Since \(N\) is of finite index,
and \(M\) is closed, \(\tilde M\) is also closed.  Since \(M\) is
not ``small,'' it is aspherical as is \(\tilde M\).  This would make
a closed 3-manifold have the homotopy type of \(S^1\x S^1\) which is
impossible.  Thus \(N\) is \(\ints\oplus\ints\oplus\ints\).

We let \(h\in\pi_1(M)\) be represented by an ordinary fiber.  Lemma
\ref{TwoFibGrps} says that \(\langle h\rangle\) is in \(N\), so
\(\tilde M\) is a primary cover of \(M\).  We can regard
\(\pi_1(\tilde M)\) as a subgroup of \(\pi_1(M)\) and thus regard
\(h\) as an element of \(\pi_1(\tilde M)\).  We know that \(h\) is
non-trivial.  (In fact it has infinite order.)  Thus \(\pi_1(\tilde
M)/\langle h\rangle\) is a proper quotient of
\(\ints\oplus\ints\oplus\ints\).  Since \(\ints\oplus\ints_n\) is
not a Fuchsian group, the quotient must be exactly
\(\ints\oplus\ints\).

From Lemma \ref{SeifFuchs}, we know that the orbit surface of
\(\tilde M\) is a torus.  Since \(\pi_1(\tilde M)\) is torsion free,
we know that \(\tilde M\) has no exceptional fibers.  The only
information needed to determine \(\tilde M\) is \(\tilde b\) the
obstruction to a section.  A presentation for \(\pi_1(\tilde M)\) is
\(\langle h,x,y\mid x^{-1}y^{-1}xyh^{\tilde b}=1\rangle\).  Since
\(h\) has infinite order in \(\pi_1(\tilde M)\) we get that \(a\)
and \(b\) commute only if \(\tilde b=0\).  This makes \(\tilde M\)
homeomorphic to \(S^1\x S^1\x S^1\).

We now consider the structure of \(M\).  If \(M\) has \(n\)
exceptional fibers, we let \((\mu_1,\beta_1),\dots,(\mu_n,\beta_n)\)
be the crossing invariants of the exceptional fibers.  Let \(b\) be
the obstruction to a section.  From the previous section we know
that \(b+\sum(\mu_i/\beta_i)\) must be a divisor of \(\tilde b\) and
therefore must be zero.  Thus each space has
\(b=-\sum(\mu_i/\beta_i)\).  This accounts for the list of closed
orientable spaces that appear among the ``flat'' Seifert fibered
spaces.  \end{demo}

We take up the question of the uniqueness of the fiber structure of
a Seifert fibered space.  We restrict ourselves to Seifert fibered
spaces with boundary.  This vastly simplifies the arguments.

\begin{lemma}\label{OneBdFib} Let \(M\) be a compact, connected
3-manifold with non-empty boundary with \(M\) neither a solid torus
nor an \(I\)-bundle over a torus or Klein bottle.  Let \(C\) be a
component of \(\bd M\) and let \(J\) and \(K\) be simple closed
curves in \(C\) that are fibers in two Seifert fiberings of \(M\).
Then \(J\) and \(K\) are isotopic in \(C\).  \end{lemma}

\begin{demo}{Proof} A fibered solid torus is the only ``small''
Seifert fibered space with non-empty boundary.  The remaining spaces
excluded by the hypothesis are the bounded ``flat'' Seifert fibered
spaces.  Note that the Seifert fibered space with orbit surface a
disk and two exceptional fibers of index 2 is homeomorpihic to an
orientable \(I\)-bundle over a Klein bottle.

By Lemma \ref{HomFiber}, \(J\) and \(K\) generate the maximal cyclic
normal subgroup of \(\pi_1(M)\).  Thus (up to reversal) they are
homotopic in \(M\).  By Lemma \ref{PowFib}, \(J\) and \(K\) are
homotopic in \(C\) to powers of each other.  Since \(\pi_1(C)\) is
\(\ints\oplus\ints\), \(J\) and \(K\) are homotopic in \(C\) and
thus isotopic in \(C\).  The isotopy can be extended to all of
\(M\).  \end{demo}

Before we go further, we need some techniques to deal with
incompressible surfaces in a Seifert fibered space.  We start with a
description of some canonical incompressible surfaces in a Seifert
fibered space with boundary.

Let \(M\) be a compact, connected Seifert fibered space with
non-empty boundary and assume that \(M\) is not a fibered solid
torus.  Let \(p:M\into G\) be the projection to the orbit surface.
We know that either \(G\) is not a disk, or \(G\) has more than one
exceptional point.  There is a finite set
\(\{\alpha_1,\dots,\alpha_n\}\) of pairwise disjoint arcs properly
embedded in \(G\) minus the exceptional points so that when \(G\) is
split along the \(\alpha_i\), then each component is a disk, no
component has more than one exceptional point, and no more than one
component has no exceptional point.  We call the saturated annuli
\(A_i=p^{-1}(\alpha_i)\) a \defit{canonical system} of saturated
annuli in \(M\).  Splitting \(M\) along the \(A_i\) represents \(M\)
as a union of fibered solid tori and no more than one ordinary solid
torus sewn together along saturated annuli in their boundaries.

\begin{lemma} Let \(M\) be a fibered solid torus and let \(A\) be a
properly embedded annulus in \(M\).  Assume that one boundary
component \(H\) of \(A\) is a fiber in \(\bd M\).  Then there is an
ambient isotopy of \(M\) rel \(H\) that carries \(A\) to a saturated
annulus in \(M\).  If \(X\) is a union of saturated annuli in \(\bd
M\) that is disjoint from \((\bd A)-H\), then the isotopy can be
taken rel \(X\).  If both components of \(A\) are fibers in \(\bd
M\), then the isotopy can be taken rel \(\bd M\).  \end{lemma}

\begin{demo}{Proof} A fiber is not trivial in \(M\).  Thus the
boundary components of \(A\) are parallel non-trivial curves in
\(\bd M\).  If \(H'\) is the component of \(\bd A\) that is not
\(H\) and \(H'\) is not a fiber, then there is an isotopy taking
\(H'\) to a fiber.  The isotopy can be taken rel \(H\) and rel \(X\)
if \(X\) is as specified in the hypothesis.  All further isotopies
will be rel \(\bd M\).

Since a fiber is not trivial in \(M\), we know that \(A\) is
incompressible in \(M\).  This expresses \(\pi_1(M)\) as a free
product with amalgamation over \(\ints\), or an HNN extension over
\(\ints\).  An HNN extension over \(\ints\) cannot yield \(\ints\)
which is \(\pi_1(M)\).  Also, since \(\pi_1(M)\) is \(\ints\), the
amalgamating subgroup must surject onto one of the factors.  Thus
one of the components \(T\) of \(M\) split along \(A\) is a product
and \(A\) is boundary parallel.  We have that \(T\) is a solid torus
with boundary made of two annuli \(A\) and \(A'\) with \(A'\) in
\(\bd M\).  Let \(A''\) be a saturated annulus parallel to \(A'\)
and properly embedded in \(T\) with \(\bd A''=\bd A'=\bd A\).  Now
\(A\) and \(A''\) are parallel and the result follows by creating an
isotopy carrying \(A\) to \(A''\).  This isotopy will be rel \(\bd
M\).  \end{demo}

\begin{lemma} Let \(M\) be a compact, connected Seifert fibered
space with non-empty boundary and let \(A\) be a properly embedded
annulus in \(M\).  Assume that one boundary component \(H\) of \(A\)
is a fiber of \(\bd M\).  Then there is an ambient isotopy of \(M\)
rel \(H\) that carries \(A\) to a saturated annulus in \(M\).  If
\(X\) is a union of saturated annuli in \(\bd M\) that is disjoint
from \((\bd A)-H\), then the isotopy can be taken rel \(X\).  If
both components of \(\bd A\) are fibers, then the isotopy can be
taken rel \(\bd M\).  \end{lemma}

\begin{demo}{Proof} If \(M\) is a fibered solid torus, then we are
done by the previous lemma so we assume that it is not.  If the
boundary component \(H'\) of \(A\) that is not \(H\) is not a fiber
then it is homotopic to a fiber.  By Lemma \ref{PowFib}, \(H'\) is
homotopic in \(\bd M\) to a fiber in \(\bd M\).  (Here we use the
fact that the only ``small'' Seifert fibered space with boundary is
a fibered solid torus.)  Since \(H'\) is a simple closed curve, it
is isotopic to a fiber in \(\bd M\).  If \(H'\) is in the same
component of \(\bd M\) as \(H\), then we can assure that the isotopy
does not disturb \(H\).  If \(X\) is as specified in the hypothesis,
then we can also assure that the isotopy does not disturb \(X\).

From now on, all isotopies will be rel \(\bd M\).  Since \(M\) is
not a fibered solid torus, it has a canonical system
\(\{A_1,\dots,A_n\}\) of saturated annuli.  Since both components of
\(\bd A\) are fibers, we can assume that neither component of \(\bd
A\) is in any of the \(A_i\).  We can isotop \(A\) to make the
intersections of \(A\) with the \(A_i\) a set of simple closed
curves that are non-trivial in \(A\) and the \(A_i\).  These are
isotopic to fibers in the \(A_i\) and so we can isotop \(A\) so that
the intersections are fibers.  We now use the previous lemma to fix
up the parts of \(A\) that remain when we split \(M\) along the
\(A_i\).  The fix will be done relative to the boundaries of the
pieces of \(M\) that are created by the splitting.  \end{demo}

\begin{lemma}\label{IsoFiber} Let \(M\) and \(N\) be compact,
connected Siefert fibered spaces with boundary and let \(f:M\into
N\) be a homeomorhpism that takes a fiber \(H\) in \(\bd M\) to a
fiber in \(\bd N\).  Then \(f\) is isotopic (rel \(H\)) to a fiber
preserving homeomorphism.  \end{lemma}

\begin{demo}{Proof} If the orbit surface of one is a disk with no
more than one exceptional point, then both \(M\) and \(N\) are
fibered solid tori.  The hypotheses say that they have the same
invariants and have the same fiber type.  An isotopy can be built in
a manner similar to previous exercises.

We now assume that \(G\), the orbit surface of \(M\) is not a disk
or it contains at least two exceptional fibers.  Let \(\alpha\) be a
properly embedded arc in \(G\) that misses all exceptional points,
that has the image of \(H\) in \(G\) as one of its boundary points,
and so that \(\alpha\) is not boundary parallel in \(G\) minus the
exceptional points.  Then the pre-image of \(\alpha\) in \(M\) is a
saturated annulus \(A\) containing \(H\) as one boundary component.
The image of \(A\) in \(N\) is an annulus that is properly embedded
in \(N\) with one boundary component a fiber of \(N\).  By the
previous lemma, there is an ambient isotopy of \(N\) carrying the
image of \(A\) to a saturated annulus in \(N\).  Thus we may assume
that the image of \(A\) is a saturated annulus in \(N\).  If we
split \(M\) along \(A\) and \(N\) along the image of \(A\), then we
obtain spaces with ``simpler'' orbit surfaces and can repeat the
process.  We can assure that the annuli used for later steps are
disjoint from the copies of \(A\) and that the isotopies do not
disturb the copies of \(A\).

We can do this for each element of a canonical system of \(M\).  At
the end we are left with fibered solid tori with the map already
fiber preserving along certain annuli in the boundary.  The lemma is
finished by fixing the maps on these fibered solid tori rel the
annuli.  This is similar to previous exercises.  \end{demo}

\begin{thm} Let \(M\) and \(N\) be compact, connected Seifert
fibered spaces with \(f:M\into N\) a homeomorphism.  Assume that
\(M\) has non-empty boundary and is homeomorphic to neither a solid
torus nor an \(I\)-bundle over a torus or Klein bottle.  Then \(f\)
is isotopic to a fiber preserving homeomorphism.  \end{thm}

\begin{demo}{Proof} By Lemma \ref{OneBdFib}, a fiber in \(\bd M\)
has image that is istopic to a fiber in \(\bd N\).  Thus we can
assume that \(f\) takes a fiber in \(\bd M\) to a fiber in \(\bd
N\).  We now apply Lemma \ref{IsoFiber}.  \end{demo}

\end{document}